\documentclass[3p,times]{elsarticle}

\usepackage{graphicx}
\usepackage[dvipsnames]{xcolor}

\usepackage{amssymb}
\usepackage{amsmath}
\usepackage{dsfont}
\usepackage{rotating}
\usepackage{multirow}

\DeclareMathAlphabet\mathbfcal{OMS}{cmsy}{b}{n}  

\newcommand{\tA}{\tilde{A}}
\newcommand{\ta}{\tilde{a}}
\newcommand{\tw}{\tilde{w}}
\newcommand{\tc}{\tilde{c}}

\newcommand{\D}{\mathbf{D}}

\def\nuu{s}

\newcommand{\Ng}{\mathcal{M}}

\newcommand{\F}{\mathbf{F}}

\newtheorem{remark}{Remark}

\def\R{{\cal R}}

\def\RR{\mathbb R}

\def\be{\begin{equation}}
	\def\ee{\end{equation}}
\def\bea{\begin{eqnarray}}
	\def\eea{\end{eqnarray}}
\def\ba{\begin{array}{l}\displaystyle}
	\def\ea{\end{array}}
\def\E{{\cal E}}

\newcommand{\SSS}{\mathbb{S}}

\newcommand{\ens}[1]{\mathbb{#1}}
\def\Ball{ {\cal B}}

\def\gf{\breve f}

\begin{document}
	
	\begin{frontmatter}
		
		\journal{Computers \& Fluids} 
		
		
		
		\title{High order modal Discontinuous Galerkin Implicit-Explicit Runge Kutta and Linear Multistep schemes for the Boltzmann model on general polygonal meshes} 
		
		\cortext[mycorrespondingauthor]{Corresponding author}
		
		\author[UNIFE,Centro]{Walter Boscheri\corref{mycorrespondingauthor}}
		\ead{walter.boscheri@unife.it}
		\author[UNIFE,Centro]{Giacomo Dimarco}
		\ead{giacomo.dimarco@unife.it}
		\address[UNIFE,Centro]{Department of Mathematics and Computer Science, University of Ferrara, Via Machiavelli 30, 44121 Ferrara, Italy}
		\address[Centro]{Center for Modeling, Computing and Statistic CMCS, University of Ferrara, Via Muratori 9, 44121 Ferrara, Italy}

%

\begin{abstract}
Deterministic solutions of the Boltzmann equation represent a real challenge due to the enormous computational effort which is required to produce such simulations and often stochastic methods such as Direct Simulation Monte Carlo (DSMC) are used instead due to their lower computational cost. In this work, we show that combining different technologies for the discretization of the velocity space and of the physical space coupled with suitable time integration techniques, it is possible to compute very precise deterministic approximate solutions of the Boltzmann model in different regimes, from extremely rarefied to dense fluids, with CFL conditions only driven by the hyperbolic transport term. To that aim, we develop modal Discontinuous Galerkin (DG) Implicit-Explicit Runge Kutta schemes (DG-IMEX-RK) and Implicit-Explicit Linear Multistep Methods based on Backward-Finite-Differences (DG-IMEX-BDF) for solving the Boltzmann model on multidimensional unstructured meshes. The solution of the Boltzmann collision operator is obtained through fast spectral methods, while the transport term in the governing equations is discretized relying on an explicit shock-capturing DG method on polygonal tessellations in the physical space. A novel class of WENO-type limiters, based on a shifting of the moments of inertia for each zone of the mesh, is used to control spurious oscillations of the DG solution across discontinuities.  
The use of Linear Multistep Methods (LMM) allows the Boltzmann solutions to be consistent not only with the compressible Euler limit but also with the Navier-Stokes asymptotic regime. In addition, as numerically proven, they also permit to strongly reduce the computational effort compared to Runge-Kutta approaches while maintaining the same or even larger accuracy. The performances of these different time discretization techniques are measured comparing both precision and efficiency. At the same time, comparisons against simpler relaxation type kinetic models such as the BGK model are proposed. The order of convergence is numerically measured for different regimes and found to agree with the theoretical findings. The new methods are validated considering two-dimensional benchmark test cases typically used in the fluid dynamics community. A prototype engineering problem consisting of a supersonic flow around a NACA 0012 airfoil with space-time-dependent boundary conditions is also presented for which the pressure coefficients are measured.
\end{abstract}
%
\end{frontmatter}

%
%
\section{Introduction}
The kinetic theory describes systems made of a large number of particles. The key idea which this theory relies on is that, instead of following the details of motion of each single particle, one describes the collective behavior of the system at a probabilistic level \cite{cercignani}. This statistical description is handled by introducing the so-called distribution function $f(x,v,t)$ defined in $\RR^{d_x}\times\RR^{d_v}\times \RR^{+}$ where $d_x,d_v$ are the physical and velocity space dimensions, respectively, and $fdxdv$ represents the probability distribution of particles to be in a given state at a given time. Thanks to this point of view, one is able to characterize a very large number of different physical phenomena ranging from classical fluids, rarefied gases and plasma \cite{ACTA,cercignani,bird} to, more recently, biological, social or economic phenomena \cite{pareschi2013BOOK}. 

One of the powerful aspects of kinetic theory is that, from the Boltzmann model, the standard compressible Euler and Navier-Stokes equations can be recovered through a suitable expansion of the distribution function in terms of the so-called Knudsen number, which is a parameter for measuring the typical time scale of the collisional dynamics \cite{ACTA,cercignani}. This means that the Boltzmann equation contains a richer spectrum of dynamics compared to the standard fluid mechanics theory. However, as well known, the counterpart of this richness is the large computational effort related to the numerical simulation of such model which has prevented for many years its applicability to realistic scenarios, especially for the case of deterministic solvers \cite{ACTA}. In fact, compared to the standard fluid model, the kinetic equations live in a larger space depending additionally on the velocity variable $v\in\RR^{d_v}$. 
Thus, their applicability is limited to the difficulties related to their numerical discretization, the so-called curse of dimensionality. Typically, Direct Simulation Monte Carlo (DSMC) methods \cite{cf, bird} are used for simulating realistic scenarios since they are the only schemes capable of describing the different spatial and time scales involved and the high dimensionality at a reasonable price in terms of computational time and memory requirements. Nevertheless, also these techniques exhibit limitations which are caused by the large statistical noise typical of the Monte Carlo methods, even if some remedies have been proposed in the past to relieve these difficulties \cite{dimarco1,dimarco2,Hadji}.

Concerning deterministic methods, it is clear that the multidimensional nature of the Boltzmann equation combined with the five fold Boltzmann operator integral causes a real challenge for the development of numerical schemes. To overcome the issue related to the numerical integration of the Boltzmann model, several authors introduced simplified models, typically of relaxation type, such as the very famous BGK and ES-BGK models \cite{Gross, Holway} or the Shakhov model \cite{Titarev1,Titarev2,Titarev4} which results can be shown to be quite effective for flows not too far from the thermodynamic equilibrium state \cite{Mieussens,perth}. We also quote the series of works on the so-called Unified-Gas-Kinetic-Schemes dealing with a special treatment of the transport part of the kinetic equations \cite{Xu1,Xu2,Xu3}. Related to this transport part, we recall some research directions based on Lagrangian techniques which have proven to be able to handle very efficiently multidimensional kinetic equations \cite{FKS,FKS_HO,FKS_GPU,FKS_Bolt}. Although very efficient, this approach was only first order in time and space and it was based on a regular Cartesian mesh. In this direction we also recall some very recent results in which high order semi-Lagrangian methods have been applied to the solution of kinetic equations even if only simple relaxation operators were considered \cite{Russo,Russo2}.

Regarding the approximation of the multidimensional full Boltzmann model by high order in space and time methods, up to our knowledge, there exist only few examples in literature. Finite volume methods both for the BGK equation as well as for the Boltzmann equation have been recently presented in \cite{BosDim,BosDim2}, which employ Implicit-Explicit Runge Kutta time discretizations. Alternatively, Discontinuous Galerkin (DG) schemes might be adopted for achieving high order of accuracy in space. These discretizations were first applied to neutron transport equations \cite{RedHill73} and later extended to general nonlinear systems of hyperbolic
conservation laws in one and multiple space dimensions in a series of papers \cite{CS-DG-1,CS-DG-2,CS-DG-3,CS-DG-4,CS-DG-5}. In the DG context, the numerical solution is approximated by polynomials within each control volume, hence leading to
a natural piecewise high order data representation. The numerical solution is represented by a polynomial expansion in terms of a set of basis functions, that can be either of nodal or modal type. Thus, DG schemes do not need any reconstruction procedure, unlike high order finite volume schemes, and therefore exhibit a much more compact stencil that is also suitable for the development of parallel solvers with MPI on several cores. Focused on kinetic models of Boltzmann type, high order nodal DG methods on unstructured meshes have been proposed in \cite{Hu3,Hu4,WU}. However, in these cases, either explicit time discretizations or low order (first order) implicit schemes have been designed. Here, instead we focus on high order in time implicit methods which are of paramount importance in kinetic equations due to the different temporal scale involved in the process under study \cite{ACTA}.

In details, in this work, we deal with high order modal Discontinuous Galerkin methods on general polygonal unstructured meshes. The methods here discussed are clearly more demanding in terms of computational resources compared to \cite{FKS} for instance. Nevertheless, differently from \cite{BosDim2}, where Runge-Kutta methods were employed, here we have a strong reduction of the computational effort thanks to the use of Linear Multistep Methods for the time integration. Moreover, as well known, DG methods are more versatile, more accurate compared to finite volume schemes, and they present a very good linear scaling on parallel machines although typically smaller CFL conditions are required for stability. With the aim of producing a high order time and space numerical method for the Boltzmann equation, we introduce in the sequel four different levels of discretization. First, the unbounded velocity space is truncated and discretized on a uniform Cartesian mesh, realizing the so-called discrete ordinate approximation or Discrete Velocity Model (DVM) \cite{bobylev} of a kinetic equation. Second, the Boltzmann integral is solved by using fast spectral methods based on the FFT technology \cite{MoPa:2006,FiMoPa:2006} which have spectral accuracy and complexity of the order of $\mathcal O(N_M^{d_v} \log(N_M^{d_v}))$ where $N_M$ is the number of modes taken for the velocity space in one direction and $d_v$ the dimension of the velocity space. We underline that the research field in this direction is very active and many recent contributions are present in the literature (see \cite{wu2015fast,wu2013deterministic,gamba:2010,Hu3,Rey3} and the references therein). Third, the transport part of the DVM model is solved by a modal DG method on arbitrarily polygonal meshes in the physical space \cite{CBS-book,CBS-convection-dominated,BosDG,BosDG2}. Emphasis is given in the construction of the scheme to the limiter procedure of the fluxes in the DG framework and in the unstructured setting \cite{Kuzmin2013,Kuzmin2014,DGCWENO}. In fact, to stabilize the DG solution, an \textit{a priori} WENO-type limiting technique is employed, which permits to find a less oscillatory polynomial approximation. Then, by using a modal basis, the evaluation of the oscillation indicators is carried out allowing higher order modes to be properly limited, while leaving the zeroth order mode untouched for ensuring conservation \cite{BosDG3}. This approach is possible because a \textit{conservative} modal basis is used and the moments of inertia can be easily shifted from one zone to another one as proposed in \cite{WENOAOunstr}. Fourth, the time discretization is handled by using two different Implicit-Explicit methods: Runge-Kutta (RK) \cite{Dimarco_stiff2} and Linear Multistep Methods (LMM) \cite{Dimarco_stiff4}. In fact as discussed in detail in \cite{ACTA}, kinetic equations are in particular hard to solve due to their multiscale nature in which collision and transport scales coexist. Close to the fluid limit, the collision rate grows exponentially, while the fluid dynamics time scale conserves a much lower pace \cite{cercignani}. Several authors have tackled the above problem in the recent past \cite{Jin2, Jin_review,Filbet_Jin, BLM,LemouAP, degondrev, DegondAP, Dimarco_stiff1, Dimarco_stiff3, Qin,Hu1,Hu2, Puppo,Rey} by developing and studying the class of methods known as Asymptotic Preserving. These techniques allow the full problem to be solved for time steps which are independent from the collisional fast scale identified by the Knudsen number. Here, we are particularly interested in comparing the RK and LMM approaches which are both stable and consistent with the limit model identified by the compressible Euler equation in the limit in which the collisional scale grows to infinity. Moreover, LMM methods can be shown to be also consistent with the compressible Navier-Stokes asymptotic regime, less computationally demanding and more accurate with respect to RK methods. Specifically, we rely on Backward-Finite-Difference (BDF) methods of LMM, because they only involve one integration point at the next (unknown) time level. Concerning this time integration technique applied to the study of kinetic equation we are aware of only two works in which they have been successfully employed \cite{Hulmm,Dimarco_stiff4}. However, these works were about the analytical properties of these schemes, with no evidences in a more applied setting. Here, we aim in showing that LMM represent a very interesting direction of research for what concerns the approximate solution of kinetic equations with respect to the more traditional RK methods.

The methods proposed in this work are successively tested on benchmark multidimensional rarefied gas dynamics problems comparing relaxation type kinetic equations with the Boltzmann model and the DG discretization with a finite volume (FV) approach. The theoretical accuracy is numerically verified for different regimes and for all methods developed hereafter. A more realistic test case involving a flow around a NACA 0012 airfoil for different values of the Knudsen number is also performed in the supersonic regime with a discontinuous change in space and time of the inflow boundary condition. A MPI parallelization is realized distributing the space variable on different threads to improve the overall computational efficiency. 

The article is organized as follows. In Section \ref{sec_Boltzmann}, we introduce the Boltzmann model and a simplified relaxation model. The properties and their fluid dynamics limits are also very briefly recalled. In Section \ref{sec_num_meth}, we present the four adopted discretizations: the discrete ordinate discretization, the Fourier approximation of the Boltzmann operator, the modal DG method in space and the IMEX Runge-Kutta and Linear Multistep schemes applied to the phase-space discretization of the kinetic equations. Section \ref{sec.validation} is devoted to test our method against some benchmark problems and one prototype engineering application involving a wing profile. Up to third order accuracy in space and time is proven and the performances of the different methods are measured.  Conclusions and future investigations are drawn in a final section.

%
%

\section{The Boltzmann model and its relaxation approximation}
\label{sec_Boltzmann}
In the rest of the paper, we consider kinetic equations in non-dimensional form of type~\cite{cercignani}
\be
\partial_t f + v\cdot\nabla_{x}f
=\frac1{\varepsilon}Q(f), \label{eq:1b}\ee
where $\varepsilon$ is the Knudsen number which is a non-dimensional quantity directly proportional to the mean free path between particles.
As already stated, $f(x,v,t)$ is a non-negative function giving the time evolution of the distribution of particles with velocity $v \in \R^{d_v}$, $d_v \geq 1$ and position $x \in \Omega \subset \R^{d_x}$, $d_x \geq 1$
at time $ t > 0$ and describing the state of the system. 

The operator $Q(f)$ characterizes the particles
interactions, which is assumed in this work to be the quadratic Boltzmann collision operator of rarefied gas dynamics, thus $Q(f)=Q_B(f,f)$ for $d_v\geq 2$. This collision term is then given by
\be Q_{B}(f,f)=\int_{\RR^{d_v}\times S^{d_v-1}} B(|v-v_*|,n_v)
[f(v')f(v'_*)-f(v)f(v_*)]\,dv_*\,dn_v, \label{eq:Q} \ee where $n_v$ is a vector of the unitary sphere $\SSS^{d_v-1} \subset \R^{d_v}$ and \be
v'=v+\frac12(v-v_*)+\frac12|v-v_*|n_v,\quad
v'_*=v+\frac12(v-v_*)-\frac12|v-v_*|n_v. \ee 
The term $B(|v-v_*|,n_v)$ is a
non-negative collision kernel characterizing the details of the
collision. In the sequel, we will assume Maxwellian type of molecules and then  $B(|v-v_*|,n_v)=b_0(\cos(\chi))$ with  $\chi=n_v\cdot(v-v_*)/|v-v_*|$. 

The local conservation properties are satisfied by the Boltzmann model: \be\int_{\R^{d_v}} \phi(v) Q(f)\,
dv=:\langle \phi Q(f)\rangle=0, \label{eq:QC}\ee where
$\phi(v)=\left(1,v,{|v|^2}/{2}\right)^T$ are the collision invariants, i.e. mass, momentum and energy are preserved by the collision dynamics. In addition, the entropy inequality \be
\frac{d}{dt}H(f) = \int_{\R^{d_v}} Q(f)\log f dv
\leq 0,\qquad H(f)=\int_{\R^{d_v}}f\log f\,dv, \label{eq:entropy} \ee
is also satisfied. This inequality permits to infer that equilibrium distributions for this problem, i.e. functions $f$ for which $Q(f)=0$, are local Maxwellian of the form
\be M[f]=M(\rho,u,T)=\frac{\rho}{(2\pi
	T)^{d_v/2}}\exp\left(\frac{-|u-v|^{2}}{2T}\right), \label{eq:M}\ee
where $\rho$, $u$, $T$ are the density, mean velocity and
temperature of the gas in the $x$-position and at time $t$ defined as
\be (\rho,\rho u,E)^T=\langle \phi f \rangle, \qquad
T=\frac1{d_v\rho}E-\frac{1}{2}\rho|u|^2, \ee 
where we have assumed the gas constant $R=1$ as effect of the adimensionalization of the equations. The high computational complexity of the Boltzmann collision operator
$Q_{B}(f,f)$ causes often its replacement by simpler relaxation operators. A leading example is given by the BGK operator, i.e. $Q(f)=Q_{BGK}(f)$, which substitutes the binary interactions with a
relaxation towards the equilibrium of the form~\cite{Gross,Holway} \be\label{BGK}
Q_{BGK}(f)=\mu(M[f]-f),\ee
where $\mu > 0$ is the collision frequency that is assumed to be equal to the macroscopic density $\rho(x,t)$. The validity of this operator, with particular care about the definition of $\mu$, in
describing the physics of non equilibrium phenomena has been the
subject of many investigations, see~\cite{cercignani,Mieussens,Gross} and references therein. Here we do not discuss this aspect and in the numerical experiments we consider a simple relaxation operator with $\mu=\rho$. Closer agreement between Boltzmann and BGK or the Ellipsoidal BGK (ES-BGK) \cite{Holway} model can be obtained by choosing more sophisticated collision frequencies. 

In the following, for simplicity, we will fix $d_x=d_v=d$, i.e. we consider the same dimension in the physical and in the velocity space. Furthermore, all numerical simulations are performed in the $d_x=d_v=2$ case which means in a four dimensional space plus the time. Moreover, for keeping notations as simple as possible, we embed the Knudsen number in the definition of the collision operator, thus we first define $\tilde Q(f)=\frac{1}{\varepsilon}Q(f)$ and we successively call $Q(f)$ this new operator $\tilde Q(f)$. 
Integration of (\ref{eq:1b}) against the collision
invariants in the velocity space leads to the following set of conservation laws \be
\partial_t \langle \phi f\rangle+{\rm div}_x
\langle v\otimes\phi f\rangle=0,\label{eq:macr}\ee
which still need a closure relation. Close to fluid regimes, the mean free path between two collisions is
very small. In this situation, passing to the limit $\varepsilon\rightarrow 0$ we formally obtain $Q(f)=0$ from (\ref{eq:1b}) and therefore $f=M[f]$. Thus, at least formally, we
recover the closed hyperbolic system of the compressible Euler equations
\be
\partial_t U+{\rm div}_x \F(U)=0,
\label{eq:Euler} \ee with
$U=\langle \phi M[f]\rangle = (\rho,\rho u,E)^T$ and
\[
\F(U)=\langle v\otimes\phi M[f]\rangle=\left(
\begin{array}{c}
	\rho u  \\
	\rho u \otimes u+pI  \\
	Eu+pu     
\end{array}
\right)
,\quad p=\rho T,\]
where $I$ is the identity matrix. Note that the above conclusions are
independent from the particular choice of $Q(f)$, provided the collision operator satisfies (\ref{eq:QC}) and admits Maxwellian of the form (\ref{eq:M}) as local equilibrium functions. For small but non zero values of the Knudsen number, the evolution equation for the moments can be derived by the so-called Chapman-Enskog expansion~\cite{cercignani}.
This originates the compressible Navier-Stokes equations
as a second order approximation with respect to $\varepsilon$ to the
solution of the Boltzmann equation
\be
\partial_t U+{\rm div}_x \F(U)=\varepsilon\, {\rm div}_x \D(\nabla_x U),
\label{eq:NavierStokes} \ee with
\be
\D(\nabla_x U)=
\left(
\begin{array}{c}
	0 \\
	\nu_{NS}\sigma(u)\\
	\kappa\nabla_x T+\nu_{NS}\sigma(u)\cdot u     
\end{array}
\right),
\quad 
\sigma(u)=\frac12\left(\nabla_x u +(\nabla_x u)^T -\frac2{3}{\rm div}_x uI\right),
\label{eq:NSsigma}
\ee
where the viscosity $\nu_{NS}$ and the thermal conductivity $\kappa$ are defined according to the chosen collisional operator.

Finally, we need to define suitable boundary conditions in physical space in order to determine the solutions for kinetic equations of type \eqref{eq:1b}. In this context, the only conditions that need to be discussed are the ones required when the fluid meets an object (the surface of a wing for instance) or equivalently it hits a wall and particles interact with the atoms of the surface before being reflected backward. For $v \cdot n \le 0$ and $x\in\partial \Omega$, where {$n$} denotes the unit normal pointing outside the domain {$\Omega$} and $\partial \Omega$ is the boundary of the domain, such boundary conditions are modeled by 
\begin{equation}
	|v \cdot n| f(x,v,t)
	= \int_{v_{\ast} \cdot n>0} |v_{\ast}\cdot n(x)| K(v_{\ast}
	\rightarrow v, x,t) f(x,v_{\ast},t)\,dv_{\ast}. \label{eq:KE} 
\end{equation}
The above expression means that the ingoing flux is defined in terms of the outgoing flux modified by a given boundary kernel {$K$}. The only condition imposed is that positivity and mass conservation at the
boundaries are guaranteed. Typically, outgoing molecules are partially reflected back into the domain and partially thermalized. When thermalized, molecules are re-emitted in the domain accordingly to a Maxwellian distribution with velocity and temperature of the object that have been encountered by the molecules (see \cite{ACTA,cercignani} for details).

%
%
\section{The numerical method}\label{sec_num_meth}
We first introduce the discrete velocity space and the discrete ordinate approach. Then, the spectral method for the Boltzmann integral is detailed. In the third part, we present the discrete spatial computational domain and the details of the unstructured mesh. On this discrete phase-space setting, we subsequently define the modal DG scheme. The Implicit-Explicit Runge-Kutta and Linear Multistep techniques are given in the last part of this section. 

\subsection{The Discrete Velocity Models (DVM)} \label{sec_DVM}
We consider the following discrete space in velocity:
\be\label{disc_space}
\mathcal{V}=\left\{ v_{(k_1,k_2)}=(k_1\Delta v_{x_1}+v_{x_1,min},k_2\Delta v_{x_2}+v_{x_2,min}), \ \Delta v_{x_1}=\frac{(v_{x_1,max}-v_{x_1,min})}{N_{x_1}}, \ \Delta v_{x_2}=\frac{(v_{x_2,max}-v_{x_2,min})}{N_{x_2}}\right\},\ee 
where $v_{min}$ is the vector indicating the lower admissible velocity and $v_{max}$ the maximum one. In the following, for simplicity and to shorten the notation, we introduce the mono-index $k$, which spans all the discrete space of $N=N_{x_1}\times N_{x_2}$ total elements, i.e. the discrete velocity will be denoted by $v_k, \ k=1,\ldots,N$. Furthermore, a regular Cartesian grid is assumed, hence implying that the mesh spacing is uniform, i.e. $\Delta v_{x_1}=\Delta v_{x_2}:=\Delta v$. In this setting, the discrete space is fixed once for all, i.e. fixed for all times, and independent of the space location. Both requirements can be relaxed in order to better describe the physics of the problem. However, this research direction is rather unexplored, and we recall a very interesting recent contribution addressing this topic \cite{Mieussens3,Rey2}. On our side, we refer to future investigations to take into account these delicate aspects. 

The continuous distribution function $f$ is then replaced by the vector $f_{\mathcal{K}}(x,t)$ of size $N$. Each component of this vector is assumed to be an approximation of the distribution function $f$ at location $v_{k}$: \be\label{disc_ord}
f_{\mathcal{K}}(x,t)=(f_{k}(x,t))_{k},\qquad f_{k}(x,t) \approx
f(x,v_{k},t). \ee  
The consequence of the truncation of the velocity space is that the exact conservation of mass, momentum and energy is impossible, because in general the support of the distribution function is non compact, as for the case of the Maxwellian equilibrium state $M[f]$. Nevertheless, taking sufficiently large bounds in velocity makes the loss of conservation negligible. If exact conservation of macroscopic quantities is demanded, $L^2$ minimization techniques can be set into place in order to locally restore the exact values of mass, momentum and energy \cite{gamba}. These are based on the search of the closest possible distribution function $\tilde f(x,v_k,t)$, in the $L^2$ sense, starting from $f(x,v_k,t)$ under the constraint that $\langle \phi \tilde f(x,v_k,t)\rangle$ has to match a given set of moments $\tilde U$, namely the ones to be conserved. However, it has to be noticed that the truncation of the velocity space is not the only responsible for the loss of conservation, since also the spectral approach largely employed in the past to solve the Boltzmann integral causes a loss of energy conservation and thus these minimization techniques have to be carried out at each time step and in each spatial cell after the collision takes place, increasing the computational effort. Here, we do not pursue this direction since the loss of conservation is assumed to be negligible for the simulations shown in this work. Anyway, we stress that the rest of the scheme which is detailed in the sequel can be adapted to exactly match the moments of the distribution, if required.

The discrete ordinate (or velocity) kinetic model consists then in the following linear hyperbolic system of $N$ equations coupled through a strongly nonlinear term on the right hand side characterizing the collisions among particles 
\be
\partial_t f_{k} + v_{k} \cdot\nabla_{x}f_{k} = Q_k(f_{\mathcal{K}}), \ k=1,\ldots,N.
\label{eq:DM1_gen} \ee
For the two cases here considered we have either $ Q_k(f_{\mathcal{K}})=Q_{k,BGK}(f_{\mathcal{K}})=\mu(\E_{k}[U]-f_{k})/\varepsilon$ for the BGK model, or $Q_k(f_{\mathcal{K}})=Q_{k,B}(f_{\mathcal{K}})$ for the Boltzmann model. In the first case, $\E_{k}[U]$ represents a suitable approximation of $M[f]$, e.g. $\E_{k}[U]=M[f](x,v_k,t)$, and $U=(\rho,\rho u ,E)^T$ is the vector of the macroscopic quantities recovered from the knowledge of the distribution function $f_{\mathcal{K}}(x,t)$ thanks to discrete summations on the discrete velocity space, that is 
\be\label{disc_summ}
U=\sum_{k=1}^{N}\Delta v^2 \, \phi_kf_k=\langle\phi_kf_k\rangle.
\ee
We claim that different type of discrete integration in the velocity space can be performed instead of \eqref{disc_summ} if desired. In the second case, $Q_{k,B}(f_{\mathcal{K}})$ corresponds to the solution given by the spectral approximation of the collision integral \eqref{eq:Q} projected over the discrete space $\mathcal{V}$. We describe in the following the fast spectral scheme used for the computation of $Q_{k,B}(f_{\mathcal{K}})$.

\subsection{Spectral discretization of the Boltzmann operator}
\label{sec:spectral}
Since collisions act only at a local in space level, we limit us in this part to deal with the approximation of $Q_B(f)$ omitting the dependence of the distribution function $f$ on space and time: $f=f(v)$. The same computation of the collision integral will be repeated for all space locations and time intervals in the full discretization as detailed later. Furthermore, the description of the method is carried out for the specific two-dimensional velocity space, that is $d_v=2$.

We summarize the method along the lines of \cite{PRspectral2000,FiMoPa:2006}. As a consequence of the discrete ordinate approach, the distribution function $f$ has compact support on the ball $\Ball_0(R)$ of radius $R$ centered in the origin. To simplify and fix the notation, we set $R=2\pi/(3+\sqrt{2})$. This choice permits to consider then $f(v)$ to be periodic on $[-\pi,\pi]^{2}$, while assuming $f(v)=0$ outside the ball $\Ball_0(2\pi/(3+\sqrt{2}))$. Periodicity is needed in order to properly define the Fourier transform of the distribution function used in the approximation of the Boltzmann operator. Moreover, the above choices are done in order to assure that, in a single computation of the collision integral, aliasing is avoided (see \cite{FiMoPa:2006} for details). Using now a single index to denote the two-dimensional sums, we have that the approximate function, denoted by $f_{N_M}$, can be continuously represented as the truncated Fourier series given by

\begin{equation}
	f_{N_M}(v)= \sum_{k=-N_m}^{N_M} \gf_k e^{i k \cdot v},
	\label{eq:FU}
\end{equation}
\begin{equation}
	\gf_k = \frac{1}{(2\pi)^{2}}\int_{[-\pi,\pi]^{2}} f(v)
	e^{-i k \cdot v }\,dv,
	\label{eq:FC}
\end{equation}
where the same index $k$, adopted for the discrete ordinate representation, is used here to span the discrete Fourier space. We stress that this representation of the solution is continuous in $v$, and $N_M$ is related to the number of modes which are employed. A spectral quadrature of the collision operator is then obtained by projecting the truncated Boltzmann operator onto the space of trigonometric polynomials of degree less or equal to $N_M$, i.e.
\begin{equation}
	{\breve Q}_k={\frac{1}{\varepsilon(2\pi)^2}}\int_{[-\pi,\pi]^{2}}
	Q_B(f_{N_M})
	e^{-i k \cdot v}\,dv, \quad {k=-N_M,\ldots,N_M}. 
	\label{eq:VAR}
\end{equation}
By substituting expression (\ref{eq:FU}) in (\ref{eq:VAR}) one gets
\begin{equation}
	{\breve Q}_k = {\sum_{\substack{l,m=-N_M\\l+m=k}}^{N_M} \gf_l\,\gf_m
		\breve\beta(l,m),\quad k=-N_M,\ldots,N_M}.
	\label{eq:CF1}
\end{equation}
The terms $\breve\beta(l,m)$ are defined as $\breve\beta(l,m)=\breve B(l,m)-\breve B(m,m)$, where $\breve B(l,m)$ are given by
\begin{equation}
	\breve B(l,m) =\frac{1}{\varepsilon} \int_{\Ball_0(2\pi/(3+\sqrt{2}))}\int_{\ens{S}^{1}} 
	|q| b_0(|q|, \cos\chi) e^{-i(l\cdot q^++m\cdot q^-)}\,d\chi\,dq, \label{eq:KM}
\end{equation}
with \begin{equation}
	q^{+} = \frac12(q+\vert q\vert \chi), \quad
	q^{-} = \frac12(q-\vert q\vert \chi).
	\label{eq:VV2}
\end{equation}
The usage of (\ref{eq:CF1}) as approximation of the collision operator requires $O(N_M^4)$ operations. In order to reduce the number of operations, one can express the Boltzmann operator adopting another representation \cite{Carl:EB:32}. Omitting the details, this second representation reads 
\begin{equation}
	\label{defQBCarleman}
	Q_B (f)=\frac{1}{\varepsilon} \int_{\R^{2}} \int_{\R^{2}} {\tilde B}(y_1,y_2) 
	\delta(y_1 \cdot y_2) 
	\left[ f(v + y_2) \, f(v+ y_1) - f(v+y_1+y_2) \, f(v) \right] \, dy_1 \,
	dy_2,
\end{equation} 
with 
\begin{equation}\label{eq:Btilde}
	\tilde{B}(|y_1|,|y_2|) =
	2 \, \sigma\left(\sqrt{|y_1|^2+|y_2|^2}, \frac{|y_1|}{\sqrt{|y_1|^2+|y_2|^2}} \right). 
\end{equation}
This transformation yields the following new spectral quadrature formula 
\begin{equation}\label{eq:ode}
	\breve{Q}_k  =
	\sum_{\underset{l+m=k}{l,m=-{N_M}}}^{{N_M}} {\breve{\beta}}_F(l,m) \, \breve{f}_l \, \breve{f}_m, \ \ \
	{k=-N_M,\ldots,N_M},
\end{equation}
where the terms $\breve B_F(l,m)$ in the definition ${\breve{\beta}}_F(l,m)=\breve B_F(l,m)-\breve B_F(m,m)$ are now 
given by
\begin{equation}
	\breve B_F(l,m) = \frac{1}{\varepsilon}\int_{\Ball_0(2\pi/(3+\sqrt{2}))} \int_{\Ball_0(2\pi/(3+\sqrt{2}))}
	\tilde{B}(y_1,y_2) \, \delta(y_1 \cdot y_2) \, 
	e^{i (l \cdot y_1+ m \cdot y_2)} \, dy_1 \, dy_2.
	\label{eq:FKM}
\end{equation}
Now, expression \eqref{eq:ode} can be evaluated relying on a convolution structure by approximating each
${\breve{\beta}}_F(l,m)$ with the sum
\[ {\breve{\beta}}_F(l,m) \simeq \sum_{p=1} ^{A} \alpha_p (l) \, \alpha' _p (m), \]
where $A$ denotes the number of finite possible directions of collisions. This gives a sum of $A$ discrete convolutions and, consequently, the algorithm can be computed in $O(A \, {N_M^2 \log_2 N_M^2})$ operations by means of
standard FFT technique (see \cite{FiMoPa:2006} for details). Typically, for the simulations shown in this work, the number of discrete angle used is very low ($A=8$). In fact, due to the use of the spectral method, this is sufficient to get enough accurate results, as demonstrated by numerical evidences. Once that $Q_B(f)$ is computed, it is projected onto the discrete ordinate space by setting 
\be
Q_{k,B}(f)\approx Q_B(f)(v_k), \quad k=1,\ldots,N,
\ee
where now $k$ is the index spanning the discrete velocity space as in the previous section.
%
%
%

\subsection{Unstructured mesh in physical space}
\label{sec_num_approx}
\begin{figure}[!htbp]
	\begin{center}
		\begin{tabular}{cc}
			\includegraphics[width=0.47\textwidth]{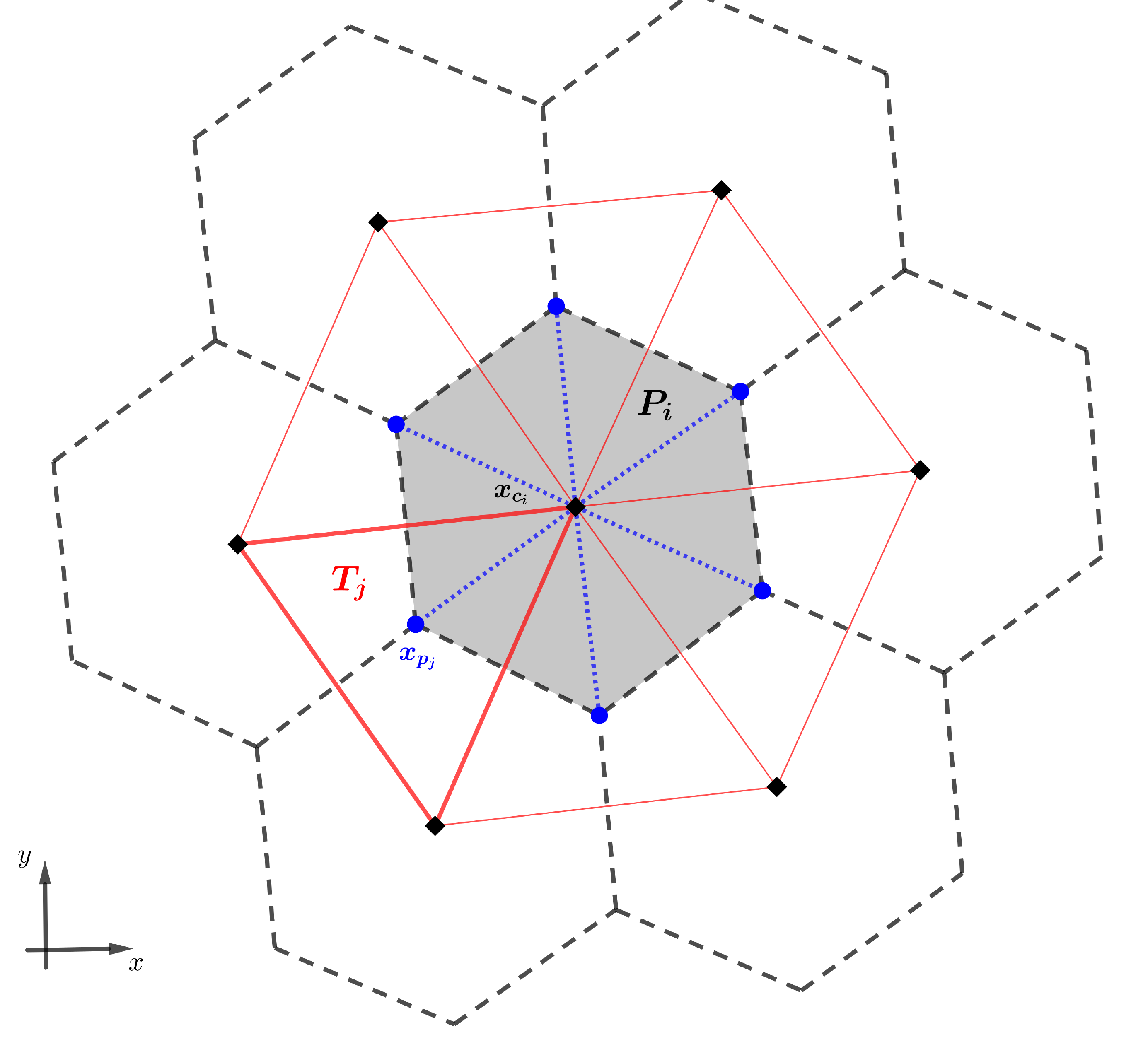}  &          
			\includegraphics[width=0.47\textwidth]{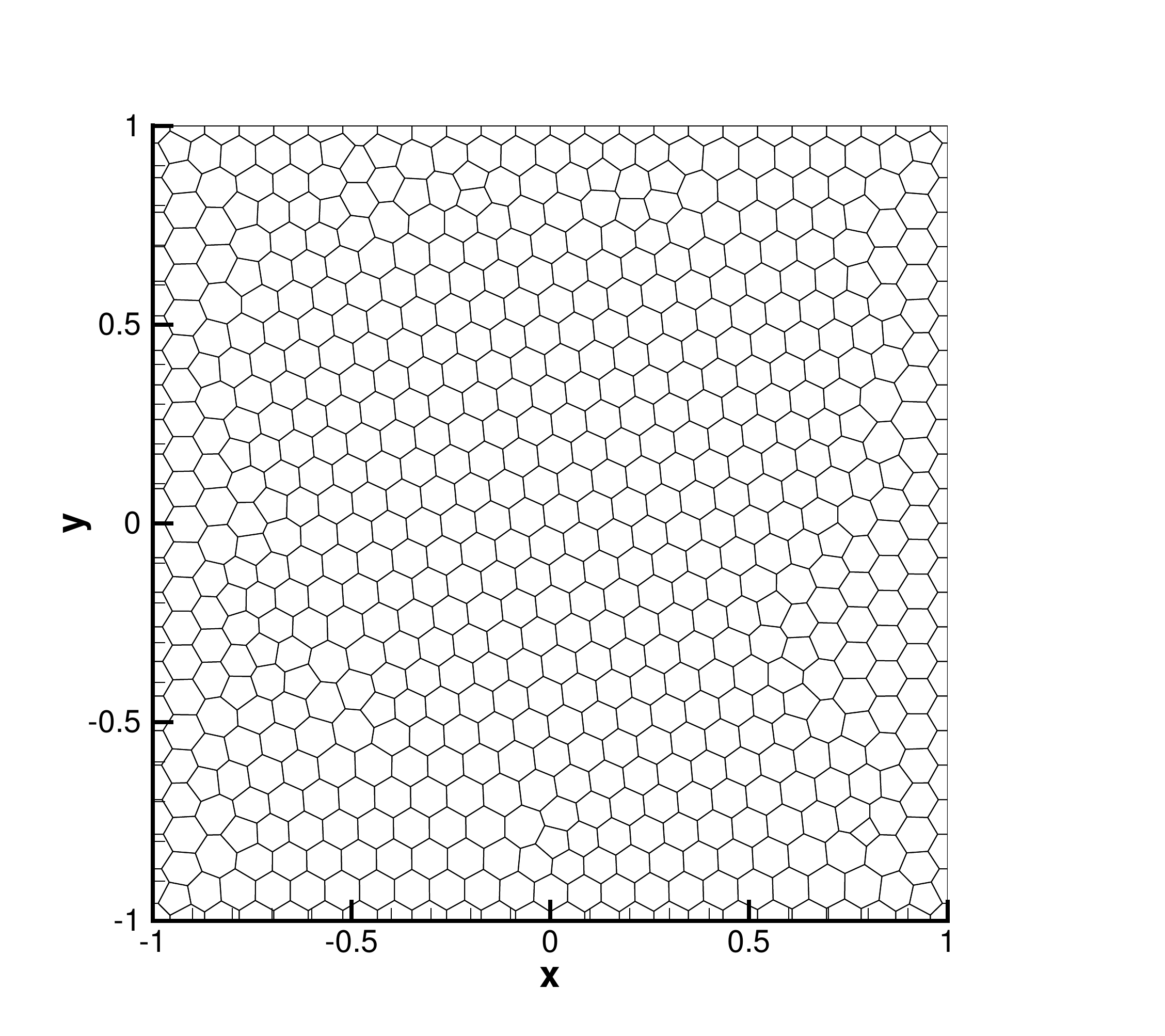} \\ 
		\end{tabular}
		\caption{Left: example of the construction of the polygonal tessellation (black dotted line) for an element $P_i$ (highlighted in gray) starting from the Delaunay triangulation (red solid line) composed of triangles $T_j$. The generator points $x_{c_i}$ are identified by black diamond markers, while the set of points lying within each primary element $T_j$ are highlighted by a blue dot and labelled with $x_{p_j}$. These points are then used to build the polygonal tessellation. The sub-triangulation used to effectively integrate the unknowns is given by the dotted blue lines. Right: example of an arbitrarily polygonal unstructured mesh.}
		\label{fig.grid}
	\end{center}
\end{figure}

We detail here how the construction of the mesh in the physical space is realized. We restrict to two-dimensional and arbitrarily shaped centroid based Voronoi-type meshes. The tessellation is made of $N_P$ non-overlapping polygons $P_i, i=1, \dots, N_P$ obtained from an underlying primary Delaunay triangulation with $N_T$ triangles. The vertexes of these triangular elements provide the position $x_{c_i}, i=1,\dots, N_P$ of the generator points of the real polygonal mesh employed in the computation. Then, inside each triangular element $T_j, j=1,\dots, N_T$, we arbitrarily choose a point of coordinates $x_{p_j}$ which will represent one of the vertices of the polygonal control volumes. Once the set of points $x_{p_j}$ is fixed, then each element $P_i$ is constructed by connecting the surrounding points $x_{p_j}\in I_{c_i}$ having the generator point $x_{c_i}$ as a vertex, see Figure \ref{fig.grid} on the left for a visual explanation. In other words, $I_{c_i}$ is the set of the Voronoi neighborhood on the triangular mesh for the vertex located at $x_{c_i}$. For example, if $x_{p_j}$ is chosen as the center of mass of element $T_j$, then a standard Voronoi tessellation is obtained, as depicted in Figure \ref{fig.grid} on the right. The center of mass $x_{m_i}$ of the obtained polygon $P_i$ is successively evaluated by
\be
x_{m_i} = \frac{1}{|P_i|} \int_{P_i} x \, dx,
\ee
This point will be used next for the definition of the basis functions in the Discontinuous Galerkin formulation. 

Let now $N_{V_i}$ denote the number of vertexes of polygon $P_i$, i.e the number of points $x_{p_j}$ belonging to each set $I_{c_i}, c_i=1,\ldots,N_p$. In order to be able to numerically integrate any quantity inside the cell $P_i$, we connect the centroid $x_{m_i}$ with each vertex of $I_{c_i}$ and we subdivide the polygon $P_i$ in $N_{V_i}$ sub-triangles. Thus integration is performed by Gauss formulae over each sub-triangles of suitable order of accuracy \cite{stroud}. This sub-triangulation is referred to as $\mathcal{T}(P_{i})$ in the following, and it is highlighted in blue in Figure \ref{fig.grid}. Before concluding, we also introduce the set $\mathcal{D}(P_i)$, which is the set of the Neumann neighborhood of element $P_i$ used also in the sequel of the paper.

\subsection{Discontinuous Galerkin approach}
The numerical solution is obtained through a space and time approximation of the discrete ordinate method which looks for a solution of $f(x,v_k,t)$ for each discrete velocity $v_k, \ k=1,\ldots,N$ (see \eqref{disc_ord}). In particular, we search an approximation of the distribution function in each element $P_i, \ i =1,\ldots,N_p$ of the spatial mesh by resorting to a linear combination of basis functions $\varphi_l, \ l=1,\ldots,\Ng$. These functions are used to span the space $\mathcal{U}_h$ of piecewise polynomials of degree less or equal to $M$ and give then the following approximate distribution function in phase-space 
\be\label{DG}
f(x,v_k,t)=f_{k}(x,t)\approx f_{h}(x,v_k,t)=f_{h,k}(x,t)=\sum_{l=1}^{\Ng}\varphi_l(x)\hat f_{k,l}(t), \ x\in P_i,
\ee
where $\hat f_{k,l}(t)$ are the unknown expansion coefficients and $\Ng$ represents the number of unknowns in the spatial element $P_i$ as an effect of the DG approximation. As classically done in the DG framework, these unknowns are referred to as the degrees of freedom, or \textit{dof}, of the numerical solution. In the following, we will also denote with
\be\label{DG_vec}
f_\mathcal{K}(x,t)=(f_k(x,t))_k\approx f_{h,\mathcal{K}}(x,t)=\left(\sum_{l=1}^{\Ng}\varphi_l(x)\hat f_{k,l}(t)\right)_k, \ x\in P_i,
\ee
the vector with components given by \eqref{DG}. Thus, the problem of finding an approximate solution to equation \eqref{eq:1b} reduces to the problem of finding at each time step $t^n$ the coefficients $\hat f_{k,l}(t^n)$ of the expansion for all mesh elements $P_i, i=1,\ldots,N_p$ and discrete velocities $v_k, \ k=1,\ldots, N$. We also recall that the relation occurring between the total number $\Ng$ of the expansion coefficients and the local polynomial degree $M$, for the two-dimensional case, is given by
\be
\Ng=\frac{1}{2}(M+1)(M+2).
\ee
A set of \textit{conservative} Taylor functions are used as \textit{modal} basis functions $\varphi_l$ for the DG approximation. Specifically, they are obtained by the direct definition of a truncated Taylor series of degree $M$ around the center of mass $x_{m_i}$ of the physical element $P_i$. This expansion is additionally normalized by the characteristic length $h_i$ of the element defined as the square root of the element surface, i.e. $h_i=\sqrt{|P_i|}$ with $|P_i|$ being the surface of polygon $P_i$. This gives, defining $x=(x_1,x_2)$,
	\be
	\label{eq.Dubiner_phi_spatial}
	\varphi_l(x_1,x_2)|_{P_i} = \frac{(x_1 - x_{m_i,1})^{r_\ell}}{h_i^{r_\ell}} \, \frac{(x_2 - x_{m_i,2})^{q_\ell}}{h_i^{q_\ell}} - \frac{1}{|P_i|}\int_{P_i}\frac{(x_1 - x_{m_i,1})^{r_\ell}}{h_i^{r_\ell}} \, \frac{(x_2 - x_{m_i,2})^{q_\ell}}{h_i^{q_\ell}} \, dx, \qquad 0 \leq {r_\ell}+{q_\ell}\leq M,
	\ee
	with $l = 1, \dots, \Ng$ denoting the same mono-index of \eqref{DG} counting the total \textit{dof} of the expansion. Let us notice that from \eqref{eq.Dubiner_phi_spatial}, it is clear that the basis functions $\varphi_l$ depends additionally on the spatial location $x$ and thus a more precise notation for them would require the indication of the mesh element $P_i$ which they belong to, i.e. one would need to use $\varphi_l^i$ to denote the basis functions. However, since the notation is already quite heavy, in the rest of the paper we will simply use $\varphi_l, \ l=1,\ldots,\Ng$, being clear the additional dependence of these functions on the mesh element. Furthermore, let us also notice that the basis functions \eqref{eq.Dubiner_phi_spatial} yield a \textit{conservative} expansion, thus it holds true that
	\be
	\frac{1}{|P_i|}\int_{P_i} \sum_{l=1}^{\Ng}\varphi_l(x)|_{P_i} \, dx = 1,
	\label{eqn.phi_cons}
	\ee
	meaning that the first degree of freedom $\hat f_{k,1}(t)$ of each element $P_i$ (i.e. the one identified by $l=1$) represents the cell average value of the unknown, in the finite volume sense. Likewise, all high order modes for $l=2, \ldots, \Ng$ automatically vanish when integrating over the control volumes $P_i$, according to the definition \eqref{eq.Dubiner_phi_spatial}. This is a crucial property for the design of suitable limiters in the DG framework, as explained later on.
		
	A weak formulation of \eqref{eq:DM1_gen} is now enforced in order to compute the approximate DG solution. In other words, we multiply the kinetic equation \eqref{eq:DM1_gen} by a test function $\varphi_m$ and we integrate over element $P_i$. This has to be done for each mesh element $P_i$ with $i=1,\ldots,N_p$ and for each velocity $v_k$ with $k=1,\ldots,N$, hence leading to
	\be\label{weak}
	\int_{P_i} \partial_t f_{h,k}(x,t)\varphi_m dx+\int_{P_i} \nabla_x \cdot\left( v_k f_{h,k}(x,t)\right)\varphi_m dx=\int_{P_i} Q_k(f_{h,\mathcal{K}})\varphi_m dx, \qquad  m=1,\ldots,\Ng.
	\ee 
	Then, by using Gauss theorem and by substituting the expansion \eqref{DG} into \eqref{weak}, we obtain
	\begin{eqnarray}
	&&\sum_{l=1}^{\Ng}\left(\int_{P_i}\varphi_l\varphi_m \, dx\right)  \partial_t \hat f_{k,l} +\int_{\partial P_i}\varphi_m\left(\mathcal{\bar F}_k^i\cdot \hat n^i \right)\, dS -\sum_{l=1}^{\Ng}\hat f_{k,l}v_k\left(\int_{P_i} \varphi_l\nabla_x \cdot( \varphi_m) \, dx\right)=\nonumber \\ &&=\int_{P_i} Q_k\left(\left(\sum_{l=1}^{\Ng}\varphi_l(x)\hat f_{k,l}\right)_k\right)\varphi_m \, dx, \qquad  m=1,\ldots,\Ng, \label{DGscheme}
	\end{eqnarray}
	where $\mathcal{\bar F}_k^i$ is the numerical flux relative to the discrete velocity $v_k$ and the mesh element $P_i$, while $\hat n^i$ is the local outward normal unit vector. The numerical flux, also known as Riemann solver, can be written as
	\be\label{num_flux}
	\int_{\partial P_i}\varphi_m\left(\mathcal{\bar F}_k^i\cdot \hat n^i \right)\, dS= \sum_{j=1}^{N_{V_i}} \int_{\partial P_{ij}} \varphi_m\left( 
	\mathcal{\bar F}_k^{ij} \cdot \ n_{ij}\right) \, dS, \qquad i=1,\ldots,N_p, \ k=1,\ldots,N, \ m=1,\ldots,\Ng,
	\ee  
	with $n_{ij}$ representing the unit outward normal pointing from element $P_i$ to the neighbor element $P_j$, $\partial P_{ij}$ denoting the face shared between element $P_i$ and $P_j$ and $\mathcal{\bar F}_k^{ij} $ being the numerical flux along the segment shared by the two control volumes $P_i$ and $P_j$. Here, we rely on a simple and very robust  Rusanov flux \cite{Rusanov:1961a} to evaluate the term $\mathcal{\bar F}_k^{ij} \cdot \ n_{ij}$. High order boundary-extrapolated data coming from the left and right piecewise polynomials defined in \eqref{DG} are used to fed the Riemann solver, which is projected in the normal direction to the boundary exploiting the Galilean invariant properties exhibited by the governing equations. The numerical flux is therefore evaluated as  
	\begin{equation}
		\mathcal{\bar F}_k^{ij} \cdot \ n_{ij} =  
		\frac{1}{2} \left( v_k f_{h,k}^+(x,t) + v_k f_{h,k}^-(x,t)  \right) \cdot n_{ij}  - 
		\frac{1}{2} s_{\max} \left( f_{h,k}^+(x,t) - f_{h,k}^-(x,t) \right), \qquad i=1,\ldots,N_p, \ k=1,\ldots,N, 
		\label{eq.rusanov} 
	\end{equation} 
	where $f_{h,k}^-(x,t)^{-}$ and $f_{h,k}^+(x,t)$ are the inner and outer high order boundary extrapolated values, respectively, relative to element  $P_i$ and $P_j$. The numerical dissipation of the Riemann solver is proportional to $s_{\max}$, that denotes the maximum eigenvalue of the system which in this case simply reads $\max_\mathcal{K}(|v_k|)$. Let us now introduce the definition of the element mass matrices
	\be\label{mass_m}
	\mathcal{M}_{ml}=\int_{P_i}\varphi_l(x)\varphi_m(x) \, dx, \qquad l=1,\ldots,\Ng, \ m=1,\ldots,\Ng,
	\ee
	and of the stiffness matrices
	\be\label{stiff_m}
	\nabla\mathcal{S}_{ml}=\left(\int_{P_i} \varphi_l(x) \partial_{x_1}\varphi_m(x) \, dx,\int_{P_i} \varphi_l(x) \partial_{x_2}\varphi_m(x) \, dx\right), \qquad l=1,\ldots,\Ng, \ m=1,\ldots,\Ng. 
	\ee
	The DG scheme \eqref{DGscheme} can be compactly written for $i=1,\ldots,N_p$, $k=1,\ldots,N$ and $l=1,\ldots,\mathcal{M}$ as
	\begin{equation}
		\sum_{l=1}^{\Ng}\mathcal{M}_{ml}\partial_t \hat f_{k,l} + \sum_{j=1}^{N_{V_i}} \int_{\partial P_{ij}} \varphi_m\left( 
		\mathcal{\bar F}_k^{ij} \cdot \ n_{ij}\right) \, dS -\sum_{l=1}^{\Ng}\hat f_{k,l}v_k\cdot\nabla\mathcal{S}_{ml} =\int_{P_i} Q_k\left(f_{h,\mathcal{K}}\right)\varphi_m \, dx.\label{DGscheme_compact}
	\end{equation}
Let us notice that, in practice, the source term on the right hand-side of \eqref{DGscheme_compact} is computed relying on the same polynomial expansion used to approximate the distribution function $f_k(x,t)$. This gives
\be\label{Bolt_DG}
Q_k\left(f_{h,\mathcal{K}}\right)\approx \sum_{l=1}^{\Ng}\varphi_l(x)\hat Q_{k,l}(t), \ x\in P_i. 
\ee

\paragraph{WENO limiter for DG on polygonal meshes} We discuss now the limiter procedure adopted to control spurious oscillations of the numerical solution in the presence of discontinuities. In fact, if the order of the polynomial reconstruction is larger then one, the numerical solution $f_{h,\mathcal{K}}$ has to undergo a process of limiting, in order to produce a stable and essentially non-oscillatory solution. The approach we employ follows the lines of the method discussed in \cite{WENODGlimiter,BosDG3}. The idea consists in using the same polynomials computed for the DG solutions in a WENO-type fashion with the aim to obtain a nonlinear convex combination based on oscillation indicators and nonlinear weights, like in high order finite volume schemes \cite{JiangShu1996,BosDim,BosDim2}. 

The first objective is then to detect the cells where a possible non physical oscillation is produced as an effect of the DG polynomial approximation. These cells are referred to as troubled cells and, in principle, they should be located only across strong discontinuities or around regions characterized by rapid variation of the macroscopic quantities, so that the excellent resolution properties typically exhibited by DG methods are not spoiled by the reduction of the polynomial reconstruction where not strictly needed. 
To detect these troubled elements, we rely on the KXRCF shock detection technique developed in \cite{KXRCF}. The element boundary is then split into two portions, that is $\partial P_i^{+}$ and $\partial P_i^{-}$, which correspond to the flow going out ($u \cdot n \geq 0$) and coming into ($u\cdot n < 0$) the cell, respectively. Here, $u$ is the macroscopic velocity of the element under consideration and $n(x)$ is the outward pointing unit normal vector referred to each side of the polygonal cell. Both the density $\rho$ and the total energy $E$ are regarded as \textit{indicator variables}. Then, a cell $P_i$ is marked as troubled if at least one of the two following conditions holds true:
\be
\frac{\left| \int_{\partial P_i^{-}} \left( \rho_i - \rho_{j} \right) \, dS \right|}{h_i^{\frac{M+1}{2}}|\partial P_i^{-}| \cdot |\rho_i|} \geq C_k, \qquad \frac{\left| \int_{\partial P_i^{-}} \left( E_i - E_{j} \right) \, dS \right|}{h_i^{\frac{M+1}{2}}|\partial P_i^{-}| \cdot |E_i|} \geq C_k.
\ee
In the above formulae, the subscript $j$ is used to denote the quantities (i.e. $\rho_j$ and $E_j$) referred to the neighbor element $P_j$ of $P_i$ across the cell boundary $\partial P_i^{-}$, while $C_k=1$ is a constant that is adopted in this work. If the cell is marked as troubled, then the corresponding troubled cell indicator $g_i$ is set to one, thus $g_i=1$. 
This indicator function permits to detect the space elements which possibly present spurious oscillations. Moreover, for the sake of safety, even those elements that are about to be crossed by a shock or a rapid variation of the macroscopic quantities, but have still to enter the wave, should be marked by the indicator function as an element which may exhibit spurious oscillations. Therefore, we modify $g_i$ by introducing a new function $\tilde g_i$ which takes into account also the Neumann neighborhood $\mathcal{D}(P_i)$ of cell $i$: 
\begin{equation}
	\tilde g_i =\max(g_i,\max \limits_{j \in \mathcal{D}(P_i)} {g_j}).
	\label{eqn.flattenerNode}
\end{equation}
Therefore the cell is flagged as troubled if $\tilde g_i>0$ and a remedy has to be found. Other detectors can be employed as well for avoiding spurious oscillations, see for instance the flattener indicator \cite{BalsaraFlattener} or the recently introduced MOOD paradigm \cite{CDL1,ALEDG,DGCWENO}. 

We are now ready to discuss the WENO-type limiting strategy that is applied only to those elements that have been flagged as troubled by the indicator previously introduced. Let us suppose that the troubled element is the one identified by $P_i$, then the WENO stencil $\mathcal{S}_i$ for cell $P_i$ is composed of the element itself and its Neumann neighbors, that is 
\begin{equation}
	\mathcal{S}_i = P_i \,\, \cup \left( P_j \in \mathcal{D}(P_i) \right). 
	\label{eqn.limstencil}
\end{equation} 
Let also $f_{h,k,i}$ and $f_{h,k,j}$ be the high order numerical solution in $P_i$ and $P_j$, respectively, given in terms of the DG approximation which in absence of oscillations are employed for the determination of the numerical fluxes in \eqref{eq.rusanov}. To ensure conservation, the integral average of all polynomials belonging to $\mathcal{S}_i$ in the WENO framework must be equal to the integral average of cell $P_i$, as usually required by finite volume reconstruction techniques \cite{BosDim,BosDim2}. Because of the chosen modal basis functions \eqref{eq.Dubiner_phi_spatial}-\eqref{eqn.phi_cons}, only the first mode assures such conservation \eqref{eqn.phi_cons} since
all the remaining higher order modes vanish as effect of the integration over the cell element. As a consequence, one can construct polynomials centered over the element $P_i$ by using the neighbor polynomials $f_{h,k,j}$ simply by \textit{shifting} their cell averages with the one of the troubled cell, i.e. cell $P_i$, and leaving untouched the higher order modes. This strategy is an extension of the parallel axis theorem in rational mechanics and has successfully been applied to finite volume and DG schemes on unstructured meshes for the solution of hyperbolic conservation laws \cite{WENOAOunstr,BosDG3}. Such modified polynomials are denoted by $\tilde{f}_{h,k,j}$, and they exhibit the same mean value of the one defined in element $P_i$. Specifically, the shifted polynomials are expressed in terms of the expansion coefficients as
\begin{equation}
	\left. \begin{array}{ll} \tilde{f}_{k,l,j} = \hat{f}_{k,l,i} & \textnormal{if} \quad l=1 \\ \tilde{f}_{k,l,j} = \hat{f}_{k,l,j} & \textnormal{if} \quad l \in [2,\Ng] \end{array} \right\} \qquad \forall P_j \in \mathcal{D}(P_i) .
	\label{eqn.shiftPoly}
\end{equation}
The above formula exactly holds true if the scaling factor between elements $P_i$ and $P_j$ is the same, i.e. $h_i=h_j$. In principle, this is never the case on unstructured meshes, therefore the shifting procedure must take into account a \textit{rescaling} of the basis functions. Here, we limit us to give an example for the second order case. For $M=1$ the basis functions \eqref{eq.Dubiner_phi_spatial} for cells $P_i$ and $P_j$ explicitly write
\begin{equation}
	\varphi_i = \left( 1, \frac{x_1 - x_{m_i,1}}{h_i}, \frac{x_2 - x_{m_i,2}}{h_i} \right), \qquad \varphi_j = \left( 1, \frac{x_1 - x_{m_j,1}}{h_j}, \frac{x_2 - x_{m_j,2}}{h_j} \right).
\end{equation}
Equating high order modes of $\varphi_i$ and $\varphi_j$ and imposing integral conservation on the cell $P_i$, the shifted \textit{dof} $\tilde{f}_{k,l,j}$ are then computed as 
\begin{equation}
	\tilde{f}_{k,l,j} = \hat{f}_{k,l,j} \frac{h_j}{h_i}, \qquad l=1,2.
\end{equation}
All the details of the aforementioned procedure can be found in \cite{WENOAOunstr}, including explicit examples of shifting polynomials up to fourth order of accuracy. Notice that, if the basis functions are defined in a \textit{reference element}, no scaling factor is needed and thus formula \eqref{eqn.shiftPoly} already provides the shifted degrees of freedom, see \cite{WAOFVun} for further details.
The new limited DG solution $f_{h,k,i}$ for cell $P_i$ (we use the same notation for the DG and the limited DG solution) is then computed by a \textit{nonlinear} weighting among the above-defined stencil polynomials $\tilde{f}_{h,k,j}$, i.e. the ones using the coefficients denoted with the tilde symbol in \eqref{eqn.shiftPoly}. The final nonlinear DG-limited polynomial $f_{h,k,i}$ is then given by 
\begin{equation}
	f_{h,k,i} = \sum \limits_{s=1}^{ \#\mathcal{D}(P_i)+1} \omega_{k,i,s} \, \tilde{f}_{h,k,s},
	\label{eqn.wenoDG} 
\end{equation}
with $\#\mathcal{D}(P_i)$ representing the number of neighbor elements of cell $P_i$. Finally, the nonlinearity is introduced as usual in the nonlinear WENO weights $\omega_{k,i,s}$
\begin{equation}
	\omega_{k,i,s} = \frac{\tilde{\omega}_{k,i,s}}{\sum_{s=1}^{\#\mathcal{D}(P_i)+1} \tilde{\omega}_{k,i,s}}, \qquad 
	\tilde{\omega}_{k,i,s} = \frac{\lambda_{i,s}}{\left(\sigma_{k,i,s} + \epsilon \right)^r}, 
	\label{eqn.omega} 
\end{equation}
through the oscillation indicators $\sigma_{k,i,s}$ with $s\in[1, \#\mathcal{D}(P_i)+1]$, hence including the central element $P_i$ as well. The oscillation indicators are simply given by
\be
\sigma_{k,i,s} = \sum \limits_{l=2}^\mathcal{M} \left(\hat f_{k,l,s}\right)^2,
\ee
with $\hat f_{k,l,s}$ denoting the expansion coefficients of the shifted polynomials \eqref{eqn.shiftPoly} for each stencil element $s$. Finally, we use $\epsilon=10^{-12}$, $r=4$, and linear weights expressed by
\begin{equation}
	\begin{array}{ll} \lambda_{i,s}=1- \left(10^{-2} \cdot  \#\mathcal{D}(P_i) \right) & \textnormal{if} \quad s=1, \\ \lambda_{i,s}=\left(1-\lambda_1\right) \cdot  \#\mathcal{D}(P_i)^{-1} & \textnormal{if} \quad s \in [2, \#\mathcal{D}(P_i)+1], \end{array}
	\label{eqn.linweights}
\end{equation}
thus ensuring $\sum_s \lambda_{i,s} = 1$. This concludes the presentation of the DG discretization. 


\subsection{Time Discretization}
\label{time_discr}
We only focus on the challenging case of the Boltzmann operator, thus we consider $Q(f)=Q_{B}(f)$ and we introduce implicit-explicit time integration techniques. These methods are needed to handle the different scales induced by the collision and the transport operators of the underlying kinetic equation. In particular, in kinetic theory one is often interested in having a method which is able to handle different regimes, dense or rarefied, with the same time step $\Delta t$ and the same accuracy (in time and space). These regimes are identified thanks to the adimensionalization proposed in Section 2 by the different values of the Knudsen number $\varepsilon$. 
In particular, the numerical scheme should be capable of describing the fluid limit, i.e. we would the discretization of the Boltzmann equation to be equivalent to a discretization of the compressible Euler equations \eqref{eq:Euler} in the limit $\varepsilon\to 0$ without however being limited by too small time steps, $\Delta t\approx \varepsilon$, caused by the fast scale collision dynamics. This request is equivalent to the notion of \textit{asymptotic-preserving (AP)} schemes \cite{ACTA,degondrev,DegondAP,Jin_review} and the schemes here introduced, both IMEX-RK as well as IMEX-LMM, belong to this special class. Specifically, we consider the extension to the DG framework of the schemes discussed from the theoretical point of view in \cite{Dimarco_stiff1} and \cite{Dimarco_stiff4}. To that aim, let us also observe that the usage of IMEX-BDF methods, i.e. of a particular subclass of the LMM schemes, permits to describe the Navier-Stokes regime \eqref{eq:NavierStokes} with the Boltzmann model under the same constraints on the time step, i.e. with CFL conditions only driven by the hyperbolic transport term.  This is a remarkable property that up to our knowledge has been only theoretically proved in \cite{Dimarco_stiff4} but never employed in applications.

In order to introduce the time discretization, we first rewrite equation \eqref{DGscheme_compact} by adding and subtracting the BGK operator $Q_{BGK}(f)$, thus relying on a \textit{ penalization method}. This is done in order to treat implicitly only the BGK operator instead of the Boltzmann collision operator. In fact, the inversion of the Boltzmann operator would require the solution of a nonlinear system that in principle has to be avoided due to the already very high computational cost needed for the solution of the multidimensional kinetic equations. The penalization technique here proposed has been used for the first time in \cite{Filbet_Jin}. We then have 
for $i=1,\ldots,N_p$, $k=1,\ldots,N$ and $l=1,\ldots,\mathcal{M}$
\begin{eqnarray}
	&&\sum_{l=1}^{\Ng}\mathcal{M}_{ml}\partial_t \hat f_{k,l} + \sum_{j=1}^{N_{V_i}} \int_{\partial P_{ij}} \varphi_m\left( 
	\mathcal{\bar F}_k^{ij} \cdot \ n_{ij}\right) \, dS -\sum_{l=1}^{\Ng}\hat f_{k,l}v_k\cdot\nabla\mathcal{S}_{ml} =\int_{P_i} Q_{B,k}\left(f_{h,\mathcal{K}}\right)\varphi_m \, dx-\int_{P_i} Q_{BGK,k}\left(f_{h,\mathcal{K}}\right)\varphi_m \, dx\nonumber\\&&+\int_{P_i} Q_{BGK,k}\left(f_{h,\mathcal{K}}\right)\varphi_m \, dx.\label{DGscheme_compact_pen}
\end{eqnarray}
The above expression can be written more compactly by inverting the mass matrices, hence obtaining
\be\label{DG_TI}
\partial_t \hat f_{h,\mathcal{K}}=\mathcal{L}(\hat f_{h,\mathcal{K}})+\mathcal{Q}_\mathcal{K}(\hat f_{h,\mathcal{K}}), \ee 
where $\mathcal{L}(\hat f_{h,\mathcal{K}})$ takes into account the DG discretization of the transport term $(v_k\cdot\nabla_x f_k), \ k=1,\ldots N$ (with the minus sign) plus the DG discretization of the difference between the two collision operators $Q_{B,k}-Q_{BGK,k}$. Finally, $\mathcal{Q}_\mathcal{K}(\hat f_{h,\mathcal{K}})$ contains the DG discretization of the collision term $Q_{BGK,k}(f_{\mathcal{K}})$ for $k=1,\ldots, N$ and all modes $l=1,\ldots,\mathcal{M}$.
We introduce now the IMEX Runge-Kutta methods \cite{Dimarco_stiff2} 
\begin{eqnarray}
	\hat F_{h,\mathcal{K}}^{(\iota)} &=& \displaystyle \hat f_{h,\mathcal{K}}^{n}+\Delta t \sum_{j=1}^{\iota-1} \ta_{\iota j} \mathcal{L}(\hat f^{(j)}_{h,\mathcal{K}})+\Delta t\sum_{j=1}^{\nu} a_{\iota j}\mathcal{Q}_\mathcal{K}(\hat f^{(j)}_{h,\mathcal{K}}), \label{eq:GIMEX} \\
	\hat f_{h,\mathcal{K}}^{n+1} &=& \displaystyle \hat f_{h,\mathcal{K}}^{n}+\Delta t
	\sum_{\iota=1}^{\nu}\tw_{\iota}\mathcal{L}(\hat F^{(\iota)}_{h,\mathcal{K}})+\Delta
	t\sum_{\iota=1}^{\nu}w_{\iota}\mathcal{Q}_\mathcal{K}(\hat F_{h,\mathcal{K}}^{(\iota)}).
	\label{eq:GIMEX1}
\end{eqnarray}
The matrices $ \tA=(\ta_{\iota j} )$,
$\ta_{\iota j} = 0$ for $j \geq \iota$ and $A = (a_{\iota j})$ are
$\nu\times \nu$ matrices such that the resulting scheme is explicit in the operator $\mathcal{L}(\hat f_{h,\mathcal{K}})$ and implicit in $\mathcal{Q}_\mathcal{K}(\hat f_{h,\mathcal{K}})$. The coefficient vectors $\tw =( \tw_{1},..,\tw_{\nu})^{T}$, $w =(w_{1},..,w_{\nu})^{T}$ complete the definition of the scheme. The schemes considered in this work are all diagonally implicit
($a_{\iota j} = 0,$ for $j> \iota$) thus the linear transport term is always evaluated explicitly. They are  conveniently represented with a compact
notation by a double Butcher tableau, as reported in Table \ref{butcher},
\begin{table}[h!]
	\begin{center}
		\begin{tabular}{l|r c}
			$\tc$ & $\tA$   \\
			\hline\\
			& $\tw^{T}$
		\end{tabular} \ \ \ \ \ \ \qquad
		\begin{tabular}{l|r c}
			$c$  & $A$   \\
			\hline\\
			& $w^{T}$
		\end{tabular}
		\caption{Butcher tableau for IMEX-RK schemes.}
	    \label{butcher}
	\end{center}
\end{table}
where the coefficients $\tc$ and $c$ are given by the
usual relation \be
\tc_{\iota}=\sum_{j=1}^{\iota-1}\ta_{\iota j}, \qquad
c_{\iota}=\sum_{j=1}^{\iota}a_{\iota j}.\ee
In practice, in the numerical section, we consider three time integration IMEX-RK schemes. The first one is the standard first order implicit-explicit Euler scheme for which we do not report the Butcher tableau. The second is the second order ARS(2,2,2) \cite{Ascher} scheme 
{\small
	\[
	\begin{array}{c|ccc}
		0 & 0 & 0 & 0   \\
		\gamma   & \gamma & 0 & 0 \\
		1   & \delta & 1-\delta & 0\\
		\hline
		& \delta &  1-\delta & 0
	\end{array}\qquad
	\begin{array}{c|ccc}
		0 & 0 & 0 & 0 \\
		\gamma   & 0 & \gamma & 0 \\
		1   & 0 & 1-\gamma & \gamma\\
		\hline
		& 0 & 1-\gamma & \gamma
	\end{array}
	\]
}
with $\gamma=1-1/\sqrt{2}$ and $\delta=1-1/(2\gamma)$, while the third is the third order BPR(3,4,3) \cite{BPR}
{\small
	\[
	\begin{array}{c|cccccc}
		0 & 0 & 0 & 0 & 0 & 0  \\
		1   & 1 & 0 & 0 & 0 & 0 \\
		2/3   & 4/9 & 2/9 & 0 & 0 & 0\\
		1   & 1/4 & 0 & 3/4 & 0 & 0\\
		1   & 1/4 & 0 & 3/4 & 0 & 0\\
		\hline
		& 1/4 & 0 & 3/4 & 0 & 0\\
	\end{array}\qquad
	\begin{array}{c|cccccc}
		0 & 0 & 0 & 0 & 0 & 0   \\
		1   & 1/2 & 1/2 & 0 & 0& 0 \\
		2/3   & 5/18 & -1/9 & 1/2 & 0& 0 \\
		1   & 1/2 & 0 & 0 & 1/2& 0 \\
		1   & 1/4 & 0 & 3/4 & -1/2& 1/2 \\
		\hline
		& 1/4 & 0 & 3/4 & -1/2& 1/2 \\
	\end{array}
	\]
}
The methods deserve some remarks.
\begin{remark}\ \\
	\vspace{-0.4cm}
\begin{itemize}
	\item The stability properties of the IMEX-RK schemes, the capability of describing different collisional regimes and the consistency with the underlying compressible Euler model have been discussed in \cite{Dimarco_stiff2} in the case of finite difference methods for the Boltzmann and the BGK equations. The case of finite volume methods has been discussed for the BGK model in \cite{BosDim} and for the Boltzmann model in \cite{BosDim2}.
	\item Despite the fact that the schemes  \eqref{eq:GIMEX}-\eqref{eq:GIMEX1} are in part implicit, the schemes can be explicitly solved. In fact, a direct evaluation of the implicit BGK collision operator can be simply done by observing that the Maxwellian distribution depends only on the first three moments of the distribution function which can be explicitly obtained by integrating \eqref{eq:GIMEX} and \eqref{eq:GIMEX1} in velocity space. By integration, all collision terms in \eqref{eq:GIMEX} and \eqref{eq:GIMEX1} disappear due to the conservation properties of $Q(f)$ recalled in \eqref{eq:QC}. Thus, the schemes \eqref{eq:GIMEX}-\eqref{eq:GIMEX1} result in a high order explicit in time approximation of the moments equations which define the conserved quantities $U$ at time $t^{n+1}$ and permit the evaluation of the Maxwellian $M(x,v,t^{n+1})$.
	\item The penalization  \eqref{DGscheme_compact_pen} can be performed by using different simplified collision operator instead of the BGK operator $Q_{BGK}(f)$. The only condition that it has to be fulfilled is that the penalization and the Boltzmann operators share the same equilibrium state. In particular, the use of the so-called ES-BGK relaxation model \cite{Holway} which is known to represent a closer approximation to the Boltzmann equation is an interesting direction of research for applications.
	\end{itemize}
\end{remark}

We discuss now the second class of IMEX method employed in this work. The IMEX Linear Multistep Methods \cite{Ascher2, Dimarco_stiff4} applied to equation \eqref{DG_TI} read
\begin{equation}
	\hat f^{n+1}_{h,\mathcal{K}} = - \sum_{j=0}^{\nuu-1} a^{LM}_j \hat f_{h,\mathcal{K}}^{n-j}+ \Delta t \sum_{j=0}^{\nuu-1} b^{LM}_j\mathcal{L}(\hat f_{h,\mathcal{K}}^{n-j})+ \Delta t \sum_{j=-1}^{\nuu-1} c^{LM}_j \mathcal{Q}_{\mathcal{K}}(\hat f^{n-j}_{h,\mathcal{K}}),
	\label{eq:IMEX_ms}
\end{equation}
where it holds $c^{LM}_{-1} \neq 0$ since these methods are implicit in $\mathcal{Q}_{\mathcal{K}}$. Those schemes for which $c^{LM}_j=0$, $j=0,\ldots,\nuu-1$ are referred to as implicit-explicit backward differentiation formula, IMEX-BDF, and they are the sole multistep methods considered here. The coefficients $a_j^{LM},b_j^{LM},c_j^{LM}, \ j=0,\ldots,\nuu$ are determined in order the full scheme to have a given order of accuracy $p$. To satisfy this condition on accuracy, the coefficients must obey the following relations
\begin{equation}
	\begin{aligned}
		&1+\sum_{j=0}^{\nuu-1} a^{LM}_j =0,\\
		1-\sum_{j=1}^{\nuu-1}& j a^{LM}_j =\sum_{j=0}^{\nuu-1} b^{LM}_j =\sum_{j=-1}^{\nuu-1} c^{LM}_j,\\
		\frac12+\sum_{j=1}^{\nuu-1} \frac{j^2}{2}a^{LM}_j&=-\sum_{j=1}^{\nuu-1} jb^{LM}_j=c^{LM}_{-1}-\sum_{j=1}^{\nuu-1}j c^{LM}_j,\\
		\vdots\\
		\frac1{p!}+\sum_{j=1}^{\nuu-1} \frac{(-j)^p}{p!} a^{LM}_j=-\sum_{j=1}^{\nuu-1}&\frac{(-j)^{p-1}}{(p-1)!} jb^{LM}_j=\frac{c^{LM}_{-1}}{(p-1)!}+\sum_{j=1}^{\nuu-1}\frac{(-j)^{p-1}}{(p-1)!} c^{LM}_j.
	\end{aligned}
	\label{eq:GIMEXcond}
\end{equation}
It can be shown (see \cite{Ascher2}) that a scheme with $\nuu$ stages is at maximum of order $\nuu$ and that, if $p\leq \nuu$, it always exists a multistep method with $s$ stages of order $s$. In Table \ref{tb:examples}, we report the three IMEX-BDF multistep schemes in vector form used in the numerical part. 
\begin{table}[t]
	\caption{The IMEX-BDF methods up to order $3$}
	\begin{center}
		{\small
			\begin{tabular}{l|c c c}
				\hline\\[-.25cm]
				Scheme & $a^{LM}$ & $b^{LM}$ & $(c^{LM}_{-1},c^{LM})$\\
				\hline\\[-.25cm]
				IMEX-BDF1 & $-1$ & $1$ & $\left(1,0\right)$\\[+.25cm]
				IMEX-BDF2 & $\left(-\frac43,\frac13\right)$ & $\left(\frac43,-\frac23\right)$ & $\left(\frac23,0,0\right)$\\[+.25cm]
				IMEX-BDF3 & $\left(-\frac{18}{11},\frac{9}{11},-\frac{2}{11}\right)$ & $\left(\frac{18}{11},-\frac{18}{11},\frac{6}{11}\right)$ & $\left(\frac{6}{11},0,0,0\right)$\\[+.25cm]
								\hline
			\end{tabular}
		}
	\end{center}
	\label{tb:examples}
\end{table}
The methods deserve some remarks.
\begin{remark}\ \\
	\vspace{-0.4cm}
	\begin{itemize}
		\item If the vector of initial steps $\hat f_{h,\mathcal{K}}^{n-j}, \ j=0,\ldots, \nuu -1$ is well-prepared,  i.e. the moments $\langle \phi \hat f_{h,\mathcal{K}}^{n-j}\rangle$ satisfy the compressible Euler equations, then in the limit $\varepsilon\to 0$, the schemes \eqref{eq:IMEX_ms}  become a class of explicit multistep schemes characterized by the coefficients $a^{LM}_j$ and $b^{LM}_j$ for the compressible Euler system \eqref{eq:Euler}. We refer to \cite{Dimarco_stiff4} for details. 
		\item Under the hypothesis in which the initial steps are well-prepared with respect to the Navier-Stokes limit, i.e. the moments $\langle \phi \hat f_{h,\mathcal{K}}^{n-j}\rangle$ satisfy the Navier-Stokes equation \eqref{eq:NavierStokes} for $j=0,\ldots, \nuu -1$, it can be shown that for small values of $\varepsilon$, the IMEX-BDF schemes \eqref{eq:IMEX_ms} become a class of explicit multistep schemes for the Navier-Stokes system \eqref{eq:NavierStokes}. We refer again to \cite{Dimarco_stiff4} for details. 
		\item As for the case of IMEX-RK schemes, despite the fact that the schemes \eqref{eq:IMEX_ms} are in part implicit, they are explicitly solvable. In fact, the implicit Maxwellian state defining the BGK operator at time $t^{n+1}$ depends on the moments of the distribution function at the same time $t^{n+1}$. These can be obtained by integrating equation \eqref{eq:IMEX_ms} in velocity space and observing that $\mathcal{Q}_{\mathcal{K}}$ as well as the DG discretization of the Boltzmann operator embedded in $\mathcal{L}(\hat f_{h,\mathcal{K}})$ disappear as the operators are mass, momentum and energy preserving. Thus, the schemes \eqref{eq:IMEX_ms} result in a high order explicit in time approximation of the moments equations which define the conserved quantities $U$ at time $t^{n+1}$ and permit the evaluation of the Maxwellian $M(x,v,t^{n+1})$. 
		\item Again as for the IMEX-RK case, the penalization  \eqref{DGscheme_compact_pen} can be performed by using different simplified collision operator instead of the BGK operator $Q_{BGK}(f)$ giving rise to more precise penalization techniques. 
		\item Compared to the IMEX-RK schemes, the IMEX-BDF need only one evaluation of the collision and of the transport operators since high order in time is reached by using the previous time levels $t^{n-j}, \ j=0,\ldots, \nuu-1$. This results in a more efficient class of numerical time integration methods. This is particularly interesting in the DG framework due their higher computational cost compared to FV methods.   
	\end{itemize}
\end{remark}

To conclude this section and the presentation of the numerical schemes, we need to define the time step restrictions. In all the numerical simulation discussed in the next section, the time step is fixed with a CFL-type stability condition
\be
\Delta t =\frac{\textnormal{CFL}}{2(M+1)}  \left(\frac{\min_{\Omega}h_i}{\max_\mathcal{K}(|v_k|)} \right),
\label{eq:Time2}
\ee
where we set $\textnormal{CFL}\leq 1/d_x$ as usual on unstructured meshes. Let observe that as opposite to the case of the numerical methods for conservation laws here the time step is fixed at the beginning of the simulation being proportional to the largest, in modulus, velocity admitted in the system that is obtained through the initial truncation of the velocity space discussed in \eqref{disc_space}. This aspect is extremely important for the application of IMEX-BDF schemes where the time step must be constant among the multistep levels of the time integrator scheme.


\section{Numerical results}
\label{sec.validation}
This section is devoted to the numerical validation of the proposed schemes. First, an extensive numerical convergence study is carried out which includes both homogeneous and non-homogeneous test problems for the Boltzmann model. Space-time convergence analysis is also carried out for different values of the Knudsen number in the non homogeneous case. Second, two benchmark compressible gas dynamics test cases are considered, namely the Lax shock tube problem and the multidimensional explosion problem. These tests are also employed for comparing the results against the simpler BGK model, as well as for measuring the efficiency of BDF linear multistep methods over classical RK time integrators. Third, a more realistic simulation is proposed, that involves the fluid flow around an airfoil profile with space-time-dependent inflow boundary conditions which induce a rapid change in the flow patterns that the scheme must be capable to capture.

We assume a Cartesian mesh for the discretization of the velocity space, thus the mesh spacing is given by $\Delta v_{x_1}=\Delta v_{x_2}:=\Delta v$. A total number of $32 \times 32=1024$ equal elements is used, unless otherwise stated. This also corresponds to the number of modes $N_M$ adopted in the spectral method for the solution of the Boltzmann collision operator. The ratio of specific heats of the gas is set to $\gamma=d_v/2=2$ for all the simulations. Boundary conditions are specified for each test case: wall boundaries are handled relying on the mass balance \eqref{eq:KE}, while Dirichlet conditions are used in part of the boundary which is not affected by any characteristic waves up to the final time of the simulation.

%
\subsection{Homogeneous test: BKW solution}
This first test aims at validating the correct implementation of the spectral method for the discretization of the Boltzmann operator. Therefore, the Boltzmann homogeneous equation is considered in the 2D velocity space and the analytical solution \cite{Bobylev1975}, that must be approximated by the spectral scheme, is given by
\begin{equation}\label{boby}
	f(v,t) = \frac{\exp{(-v^2/2S)}}{2\pi S^2} \Bigg[2S - 1 + \frac{1-S}{2S}v^2\Bigg]	\qquad	\text{with}	\qquad	S=1-\frac{\exp{(-t/8)}}{2}.	
\end{equation}
The velocity domain is defined as $\mathcal{V}=[-12;12]\times [-12;12]$, for which different discretizations are set, that also correspond to the number of modes $N_M$ of the Fourier expansion. Specifically, we use $N=N_M=\{16,24,32,48,64\}$. The errors in the $L_1$ and $L_\infty$ norms are obtained with the following formulae
\begin{equation}
	L_1 (\Delta v, t) = \frac{\sum_{k=1}^{N} |f_{k,\Delta v}(t) - f_{ref}(v_k,t)|}{\sum_{k=1}^{N} |f_{ref}(v_k,t)|}, \qquad L_\infty (\Delta v, t) = \max \limits_{k \in N} |f_{k,\Delta v}(t) - f_{ref}(v_k,t)|,
	\label{error}
\end{equation}
with $\Delta v$ representing the characteristic mesh size in the velocity space and $f_{k,\Delta v}(t)$ the approximated solution in the velocity points as defined in \eqref{disc_ord}. This is obtained after the back transformation into physical variables of the solution obtained by the Fourier method described in Section 3.2. The quantity $f_{ref}$ denote the reference solution computed as well on the same velocity points. The Table \ref{tab.conv_BKW} reports the convergence studies carried out at three output times, namely $t=3$, $t=6$ and $t=9$. As normally achieved by spectral methods, the convergence rate grows very fast, hence obtaining errors that become closer to machine precision as the number of modes is refined.

\begin{table}[!htbp]  
	\caption{Convergence study for the computation of the collision term $Q(f)$ for the BKW test. The errors are measured in normalized $L_1$ and $L_\infty$ norm at times $t=3$, $t=6$ and $t=9$ for different mesh discretizations in the velocity space $N=N_M=\{16,24,32,48,64\}$.}  
	\begin{center} 
		\begin{small}
			\begin{tabular}{c|cccc} 
				\multicolumn{5}{c}{\textbf{time = 3.0}} \\
				$N_M$ & $L_1$ & $\mathcal{O}(L_1)$ & $L_\infty$ & $\mathcal{O}(L_\infty)$ \\ 
				\hline
				16 & 2.480E+00 & -     & 3.358E-03 & - \\
				24 & 2.662E-01 & 5.50  & 4.118E-04 & 5.18 \\
				32 & 3.847E-02 & 6.72  & 6.739E-05 & 6.29 \\
				48 & 6.867E-05 & 15.61 & 6.029E-08 & 17.31 \\
				64 & 4.735E-09 & 33.31 & 3.167E-12 & 34.25 \\ 																
				\multicolumn{5}{c}{} \\
				\multicolumn{5}{c}{\textbf{time = 6.0}} \\
				$N_M$ & $L_1$ & $\mathcal{O}(L_1)$ & $L_\infty$ & $\mathcal{O}(L_\infty)$ \\ 
				\hline
				16 & 6.007E+00 & -     & 1.008E-03 & - \\
				24 & 1.241E-01 & 9.57  & 7.940E-05 & 6.27 \\
				32 & 2.439E-02 & 5.66  & 1.362E-05 & 6.13 \\
				48 & 1.577E-05 & 18.11 & 4.112E-09 & 19.99 \\
				64 & 4.500E-10 & 36.37 & 6.834E-13 & 30.25 \\ 																
				\multicolumn{5}{c}{} \\
				\multicolumn{5}{c}{\textbf{time = 9.0}} \\
				$N_M$ & $L_1$ & $\mathcal{O}(L_1)$ & $L_\infty$ & $\mathcal{O}(L_\infty)$ \\ 
				\hline
				16 & 1.763E+01 & -     & 9.685E-04 & - \\
				24 & 2.892E-01 & 10.14 & 2.517E-05 & 9.00 \\
				32 & 1.290E-02 & 10.81 & 3.003E-06 & 7.39 \\
				48 & 5.813E-06 & 19.00 & 5.658E-10 & 21.15 \\
				64 & 1.917E-09 & 27.87 & 3.355E-13 & 25.83 \\ 										
			\end{tabular}
		\end{small}
	\end{center}
	\label{tab.conv_BKW}
\end{table}

\begin{figure}[!htbp]
	\begin{center}
		\begin{tabular}{cc} 
			\includegraphics[width=0.47\textwidth]{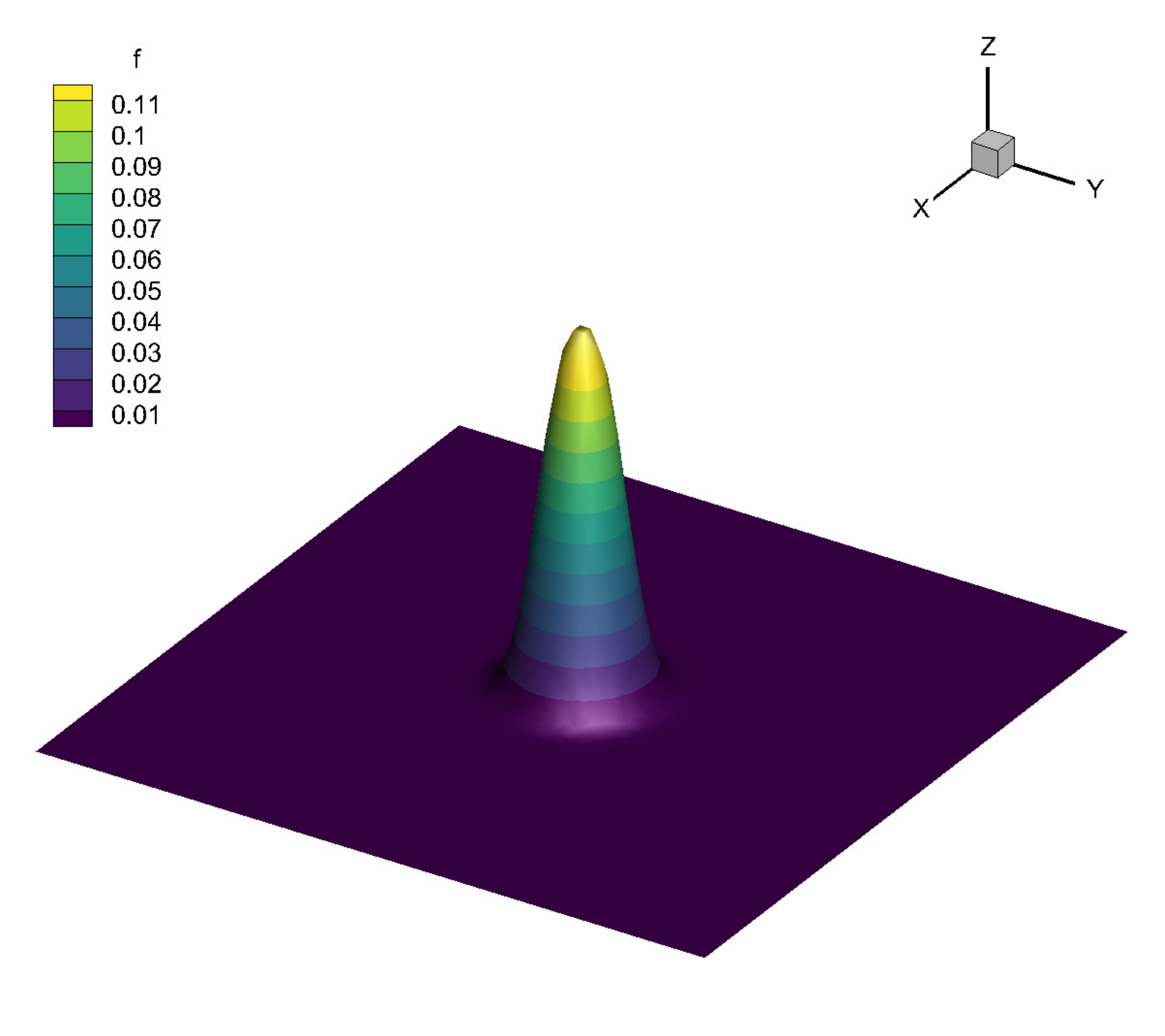}  &           
			\includegraphics[width=0.47\textwidth]{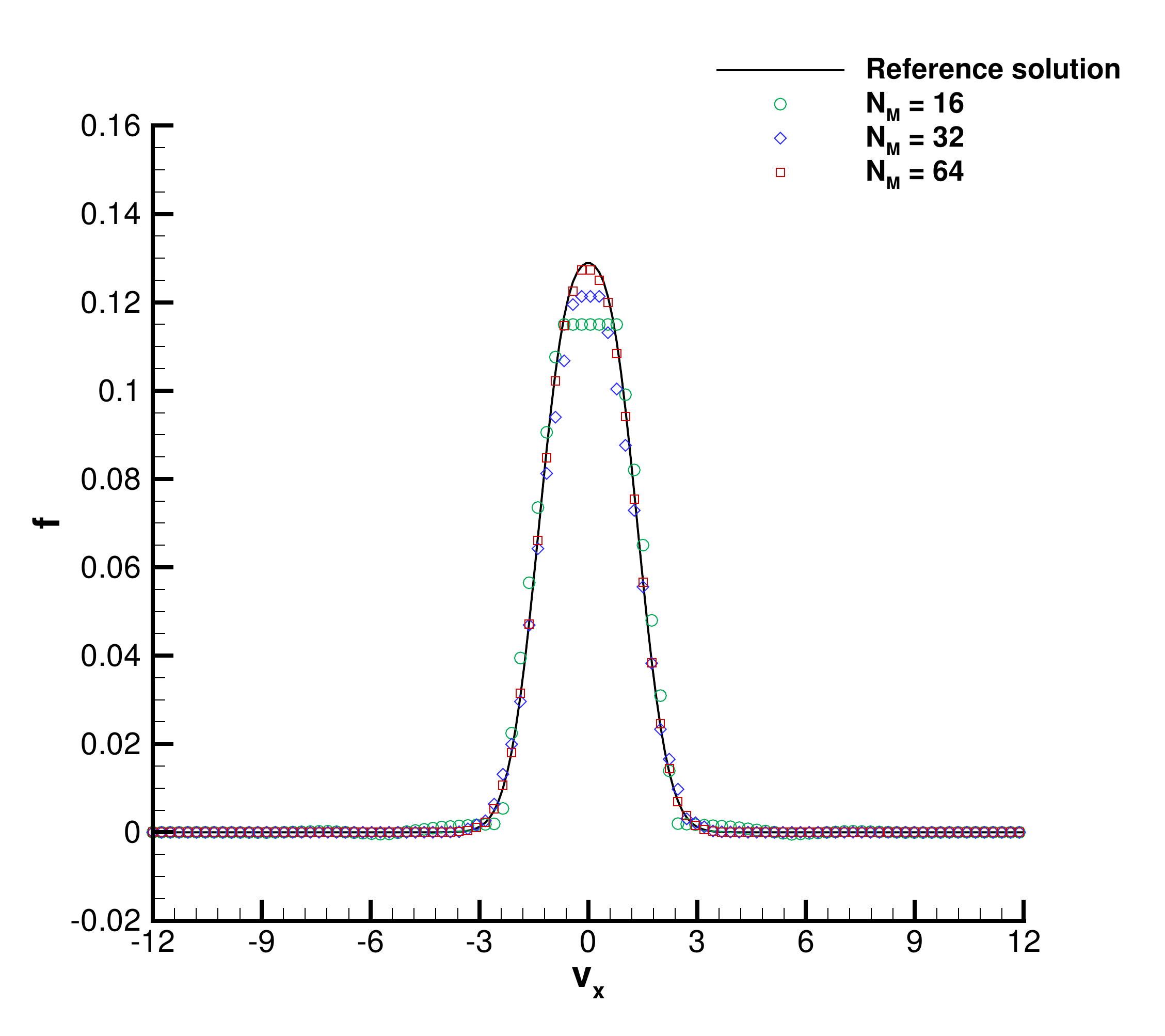} \\
			\includegraphics[width=0.47\textwidth]{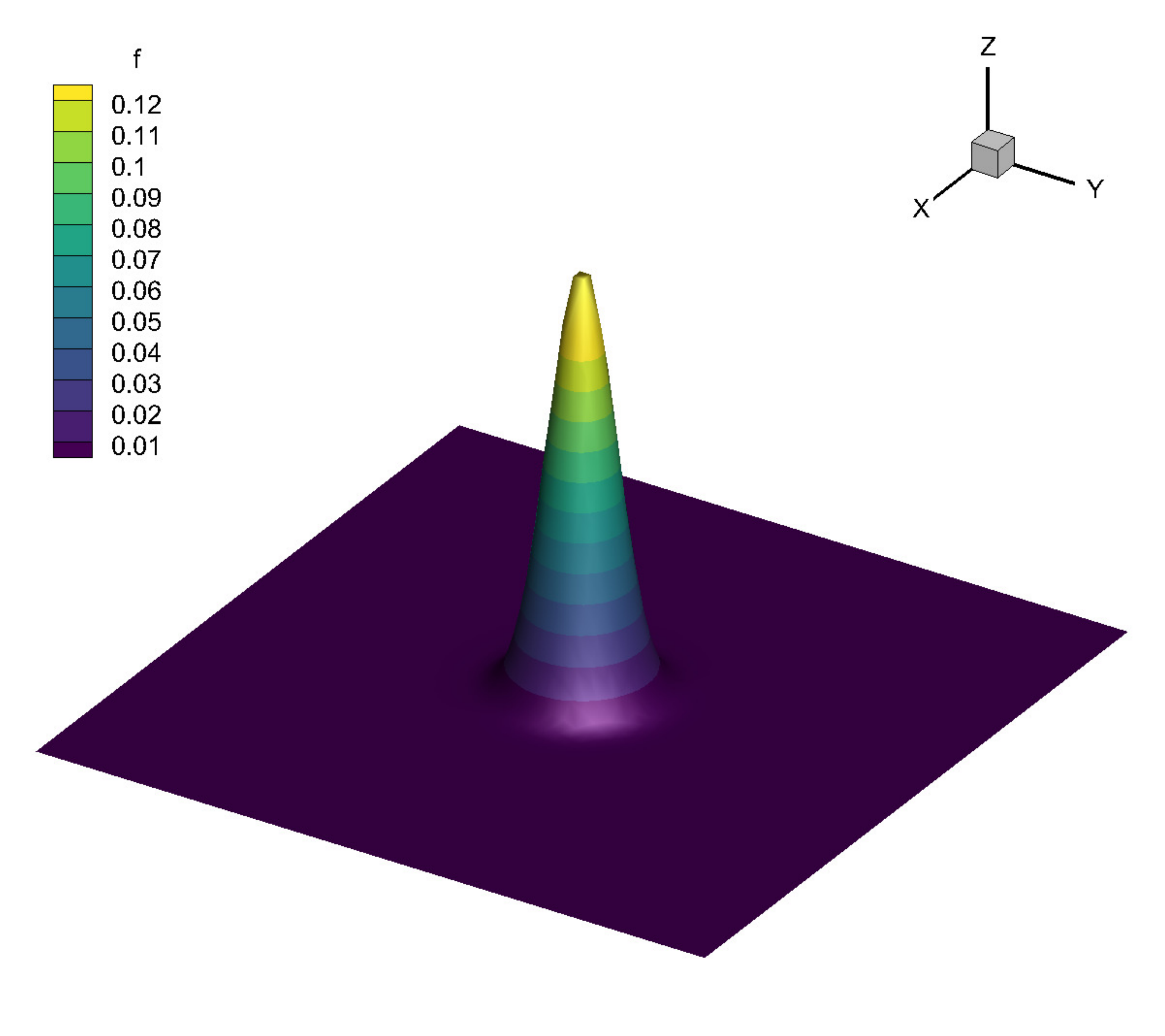}  &
			\includegraphics[width=0.47\textwidth]{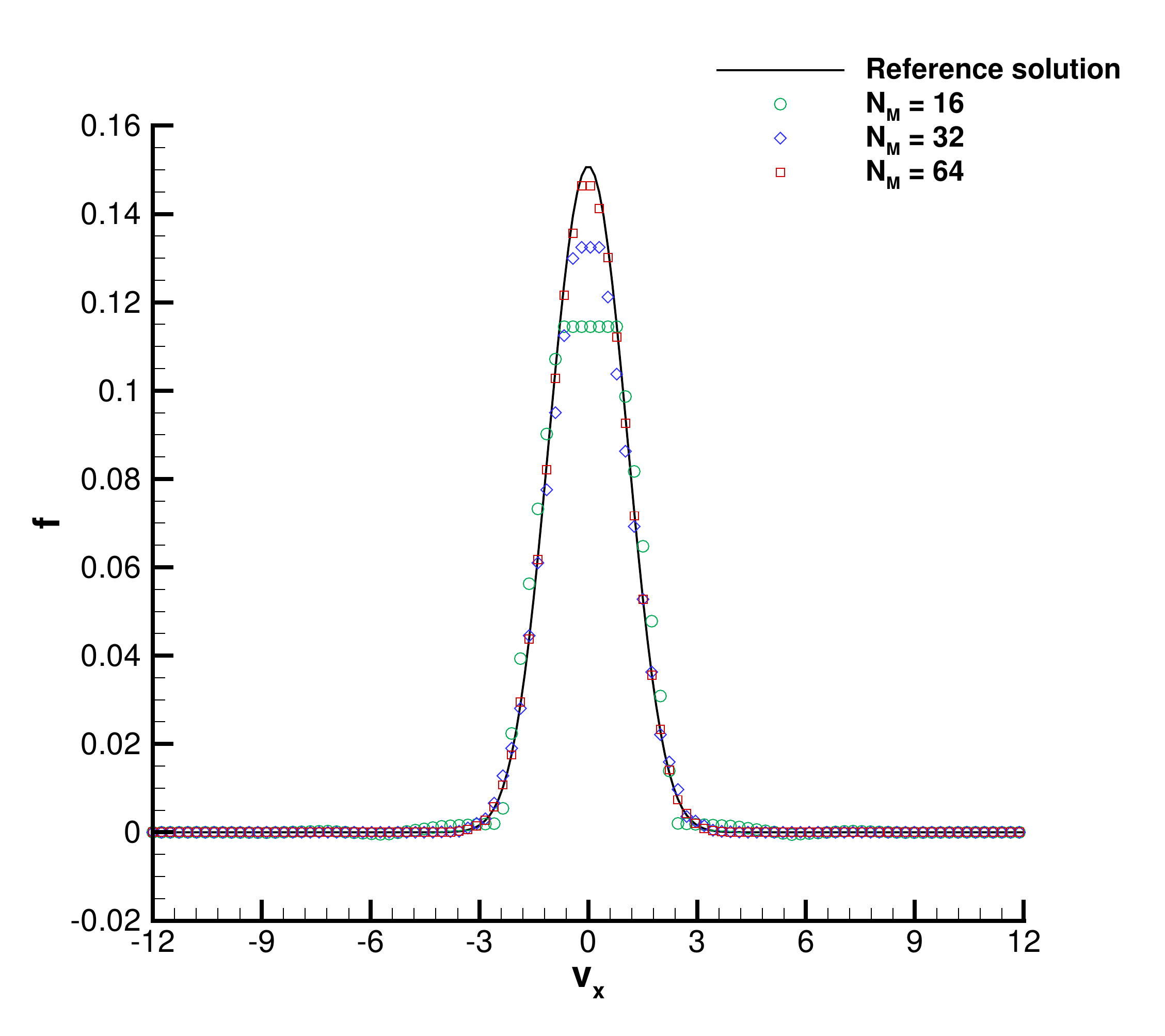} \\
		\end{tabular}
		\caption{BKW test. Distribution of $f$ (left) and comparison of the solution (right) for $N_M=16$, $N_M=32$ and $N_M=64$ at time $t=4$ (top row) and $t=8$ (bottom row). The reference solution is the exact solution \eqref{boby} represented using $N=256$ points. 
		The numerical solution is extracted as a 1D cut of 200 equidistant points along the $x$-axis at $y=N/2$.}
		\label{fig.BKW_f3D}
	\end{center}
\end{figure}

\begin{figure}[!htbp]
	\begin{center}
		\begin{tabular}{cc} 
			\includegraphics[width=0.47\textwidth]{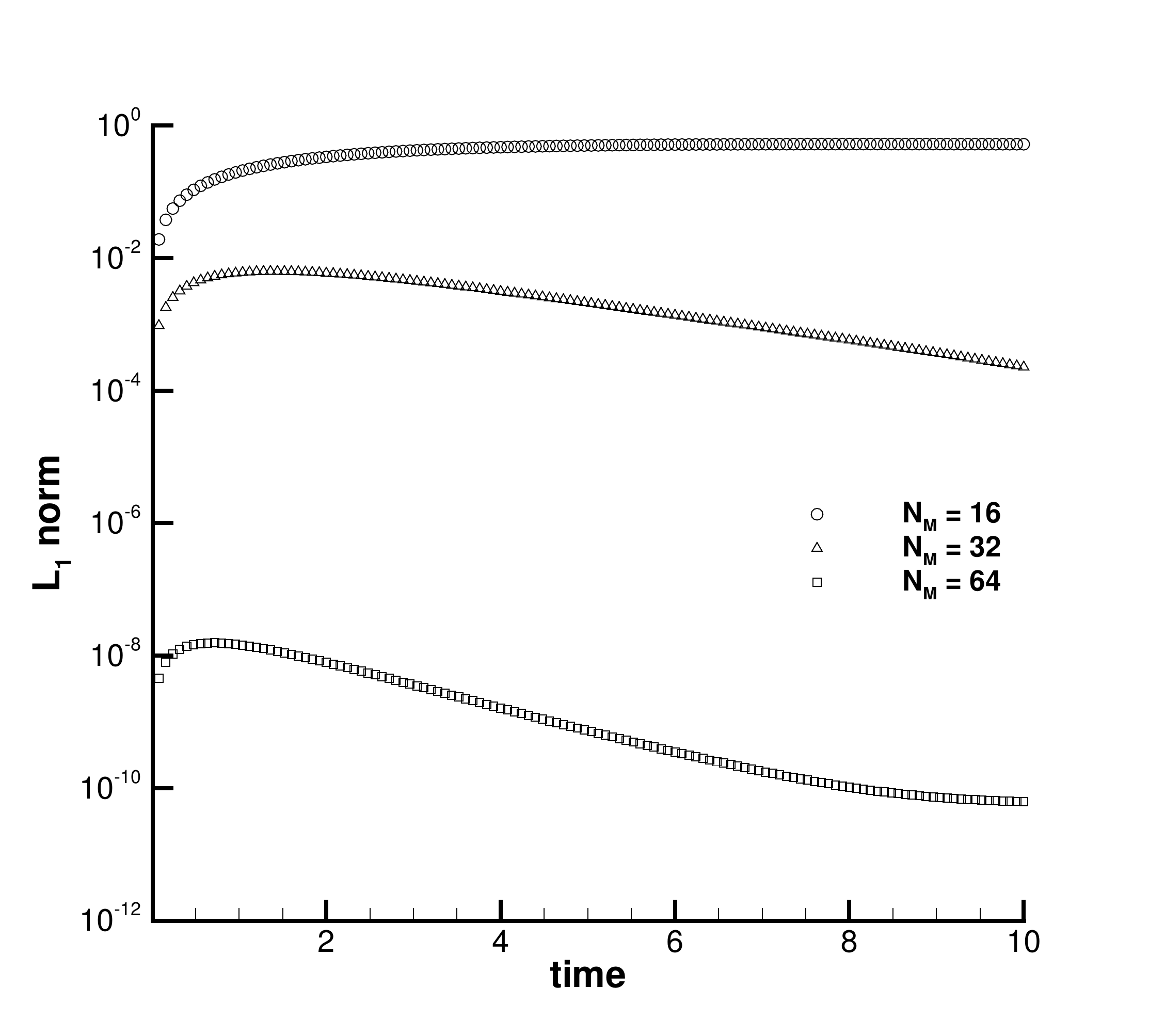}  &           
			\includegraphics[width=0.47\textwidth]{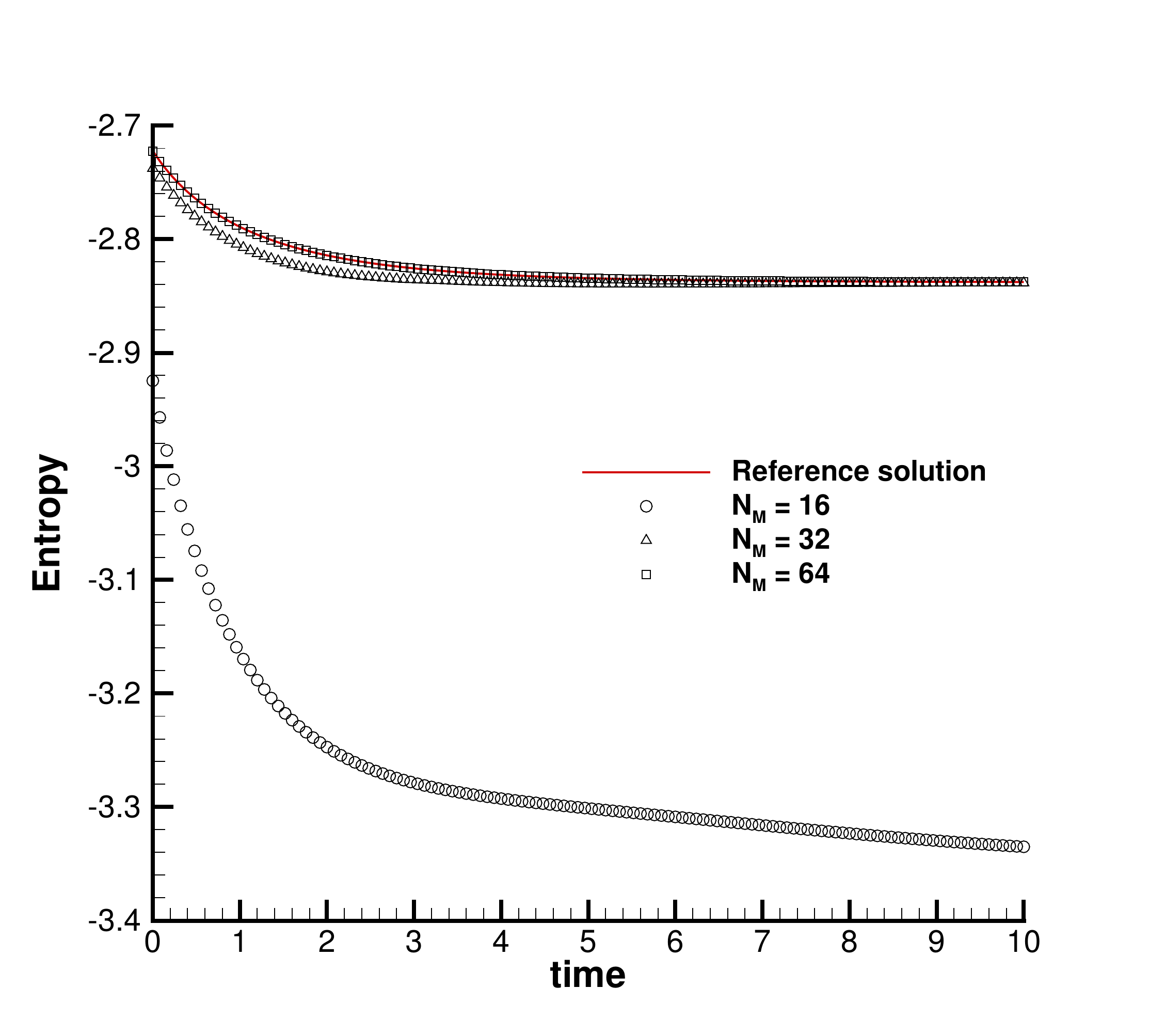} \\
		\end{tabular}
		\caption{BKW test. Time evolution of the $L_1$ errors (left) and entropy (right) for different discretizations of the velocity space with $N_M=16$, $N_M=32$ and $N_M=64$. The reference solution is corresponds to the exact solution \eqref{boby} project over $N=256$ points.}
		\label{fig.BKW_time}
	\end{center}
\end{figure}

Figure \ref{fig.BKW_f3D} depicts the numerical solution for the distribution of $f$ at two different output times ($t=4$ and $t=8$) as well as a comparison between numerical and analytical solution for different values of modes $N_M$ employed. Finally, the time evolution of the $L_1$ norm is displayed in Figure \ref{fig.BKW_time} for all three simulations. Entropy is also monitored, verifying that it is always a decreasing function. Relying on the mid-point rule, the entropy $\eta$ can be approximated as 
\begin{equation}
	\eta = \int \limits_{\mathbb{R}^2} f \, \ln f \, dv \approx \sum_{k \in N} f_k \, \ln f_k \,  \Delta v^2,
\end{equation}
and the entropy time evolution is shown in the right panel of Figure \ref{fig.BKW_time}. Let us notice that with $N_M=32$ a very precise solution is already obtained, thanks to the fast convergence of the spectral method. Consistently, there is almost no difference between the solution obtained with $N_M=64$ and the reference solution, 
meaning that for this test case the errors have already reached machine precision.

%
\subsection{Numerical convergence studies at different Knudsen numbers}
The numerical convergence of the novel DG schemes is studied by firstly considering the stiff limit of the Boltzmann equation, i.e. when it approaches the compressible Euler equations. We consider a problem as the one described in \cite{HuShuTri} where a two-dimensional smooth isentropic vortical flow is employed to measure the rate of convergence of the numerical method. Let $\Omega=[0;10] \times [0;10]$ be the computational domain, where periodic boundaries are set everywhere. A sequence of refined unstructured polygonal meshes is used to discretize the physical space, with characteristic mesh size $h(\Omega)=\big(\sum_{i=1}^{N_P} h_i\big)/N_P$ with $h_i=\sqrt{|P_i|}$ and $|P_i|$ denoting the surface of the cell $P_i$. The velocity space is bounded by $\mathcal{V}=[-10;10]\times [-10;10]$ and the initial condition is given by
\begin{equation}
	U = (\rho,u_x,u_y,T) = (1+\delta \rho, 1+\delta u_x, 1+\delta u_y, \delta T),
	\label{eq.ConvEul-IC}
\end{equation}
where the perturbations for temperature $\delta T$, density $\delta \rho$ and velocity $(\delta u_x, \delta u_y)$ are
\begin{eqnarray}
	\label{ShuVortDelta}
	\delta T    = -\frac{(\gamma-1)\beta^2}{8\gamma\pi^2}e^{1-R^2},\quad \delta \rho = (1+\delta T)^{\frac{1}{\gamma-1}}-1,\quad 
	\left(\begin{array}{c} \delta u_x \\ \delta u_y \end{array}\right) = \frac{\beta}{2\pi}e^{\frac{1-R^2}{2}} \left(\begin{array}{c} -(y-5) \\ \phantom{-}(x-5) \end{array}\right),  
\end{eqnarray}
with $\beta=5$ and $R=\sqrt{x^2+y^2}$ denoting the generic radial coordinate. The initial condition also corresponds to the analytical solution for this smooth stationary vortex flow. The final time of the simulation is chosen to be $t_f=0.1$ and the errors are reported in Tables \ref{tab.conv_0-RK} and \ref{tab.conv_0-BDF} for DG-IMEX-RK and DG-IMEX-BDF schemes, respectively. By setting $\varepsilon=0$ a fully space-time convergence analysis is carried out, showing that the formal second and third order of accuracy is properly achieved by the DG schemes. The errors are measured in the $L_1$, $L_2$ and $L_\infty$ norms for density and temperature as follows:
\begin{eqnarray}
	L_1 &=& \int \limits_{\Omega} \left| U_{ref}(x,y) - U_h(x,y) \right| dx \, dy, \\
	L_2 &=& \sqrt{\int \limits_{\Omega} \left( U_{ref}(x,y) - U_h(x,y) \right)^2 dx \, dy} \\
	L_\infty &=& \max \limits_{\Omega} \left| U_{ref}(x,y) - U_h(x,y) \right|.
\end{eqnarray}	
Here, $U_h(x,y)$ represents the high order DG solution for the macroscopic quantities obtained with our scheme, while $U_{ref}(x,y)$ is a prescribed \textit{reference solution}, that is given by \eqref{eq.ConvEul-IC}-\eqref{ShuVortDelta} for the limit case $\varepsilon=0$ being the solution studied stationary. We observe that, as expected \cite{Ascher2}, the IMEX-BDF schemes systematically give a lower error, thus exhibiting less numerical diffusion compared to IMEX-RK time stepping techniques.

\begin{table}[!htbp]  
	\caption{Numerical convergence results for the Boltzmann model using second and third order DG-IMEX-RK schemes at time $t_f=0.1$ with $\varepsilon=0$ on a sequence of refined polygonal meshes of size $h(\Omega)$. The errors are measured in $L_1$, $L_2$ and $L_\infty$ norm and refer to the variables $\rho$ (density) and $T$ (temperature).}  
	\begin{center} 
		\begin{small}
			\renewcommand{\arraystretch}{1.0}
			\begin{tabular}{c|cccccc} 
				\multicolumn{7}{c}{DG-IMEX-RK $\mathcal{O}2$} \\
				\hline
				$h(\Omega)$ & $\rho_{L_1}$ & $\mathcal{O}(\rho_{L_1})$ & $\rho_{L_2}$ & $\mathcal{O}(\rho_{L_2})$ &  $\rho_{L_\infty}$ & $\rho_{L_\infty}$ \\ 
				\hline
				3.83E-01 & 3.260E-02 & -    & 1.019E-02 & -    & 9.824E-03 & -    \\ 
			    1.96E-01 & 8.195E-03 & 2.06 & 2.609E-03 & 2.03 & 2.507E-03 & 2.04 \\ 
				1.32E-01 & 3.711E-03 & 1.99 & 1.175E-03 & 2.01 & 1.826E-03 & 0.80 \\ 
				9.90E-02 & 2.088E-03 & 2.01 & 6.677E-04 & 1.97 & 9.920E-04 & 2.13 \\ 
				\hline
				$h(\Omega)$ & $T_{L_1}$ & $\mathcal{O}(T_{L_1})$ & $T_{L_2}$ & $\mathcal{O}(T_{L_2})$ &  $T_{L_\infty}$ & $T_{L_\infty}$ \\ 
				\hline
				3.83E-01 & 8.503E-02 & -    & 2.788E-02 & -    & 2.523E-02 & -    \\ 
				1.96E-01 & 2.071E-02 & 2.11 & 6.958E-03 & 2.07 & 6.621E-03 & 2.00 \\ 
				1.32E-01 & 9.216E-03 & 2.04 & 3.105E-03 & 2.03 & 4.985E-03 & 0.71 \\ 
				9.90E-02 & 5.153E-03 & 2.03 & 1.747E-03 & 2.01 & 2.022E-03 & 3.15 \\
				\multicolumn{7}{c}{} \\
				\multicolumn{7}{c}{DG-IMEX-RK $\mathcal{O}3$} \\
				\hline
				$h(\Omega)$ & $\rho_{L_1}$ & $\mathcal{O}(\rho_{L_1})$ & $\rho_{L_2}$ & $\mathcal{O}(\rho_{L_2})$ &  $\rho_{L_\infty}$ & $\rho_{L_\infty}$ \\ 
				\hline
				3.83E-01 & 4.935E-03 & -    & 1.218E-03 & -    & 1.210E-03 & -    \\ 
				1.96E-01 & 7.147E-04 & 2.88 & 2.029E-04 & 2.68 & 2.940E-04 & 2.11 \\ 
				1.32E-01 & 2.322E-04 & 2.83 & 6.745E-05 & 2.77 & 1.417E-04 & 1.83 \\ 
				9.90E-02 & 1.022E-04 & 2.87 & 2.969E-05 & 2.87 & 4.909E-05 & 3.71 \\ 
				\hline
				$h(\Omega)$ & $T_{L_1}$ & $\mathcal{O}(T_{L_1})$ & $T_{L_2}$ & $\mathcal{O}(T_{L_2})$ &  $T_{L_\infty}$ & $T_{L_\infty}$ \\ 
				\hline
				3.83E-01 & 1.309E-02 & -    & 3.650E-03 & -    & 3.523E-03 & -    \\ 
				1.96E-01 & 1.550E-03 & 3.19 & 4.642E-04 & 3.08 & 5.624E-04 & 2.74 \\ 
				1.32E-01 & 4.764E-04 & 2.97 & 1.422E-04 & 2.97 & 3.150E-04 & 1.46 \\ 
				9.90E-02 & 2.046E-04 & 2.95 & 6.118E-05 & 2.95 & 8.494E-05 & 4.58 \\
			\end{tabular}
		\end{small}
	\end{center}
	\label{tab.conv_0-RK}
\end{table}

\begin{table}[!htbp]  
	\caption{Numerical convergence results for the Boltzmann model using second and third order DG-IMEX-BDF schemes at time $t_f=0.1$ with $\varepsilon=0$ on a sequence of refined polygonal meshes of size $h(\Omega)$. The errors are measured in $L_1$, $L_2$ and $L_\infty$ norm and refer to the variables $\rho$ (density) and $T$ (temperature).}  
	\begin{center} 
		\begin{small}
			\renewcommand{\arraystretch}{1.0}
			\begin{tabular}{c|cccccc} 
				\multicolumn{7}{c}{DG-IMEX-BDF $\mathcal{O}2$} \\
				\hline
				$h(\Omega)$ & $\rho_{L_1}$ & $\mathcal{O}(\rho_{L_1})$ & $\rho_{L_2}$ & $\mathcal{O}(\rho_{L_2})$ &  $\rho_{L_\infty}$ & $\rho_{L_\infty}$ \\ 
				\hline
				3.83E-01 & 3.249E-02 & -    & 1.019E-02 & -    & 9.824E-03 & -    \\ 
				1.96E-01 & 8.161E-03 & 2.06 & 2.604E-03 & 2.04 & 2.510E-03 & 2.04 \\ 
				1.32E-01 & 3.703E-03 & 1.99 & 1.170E-03 & 2.01 & 1.460E-03 & 1.41 \\ 
				9.90E-02 & 2.083E-03 & 2.01 & 6.651E-04 & 1.98 & 8.301E-04 & 1.90 \\ 
				\hline
				$h(\Omega)$ & $T_{L_1}$ & $\mathcal{O}(T_{L_1})$ & $T_{L_2}$ & $\mathcal{O}(T_{L_2})$ &  $T_{L_\infty}$ & $T_{L_\infty}$ \\ 
				\hline
				3.83E-01 & 8.478E-02 & -    & 2.787E-02 & -    & 2.522E-02 & -    \\ 
				1.96E-01 & 2.066E-03 & 2.11 & 6.958E-03 & 2.07 & 6.597E-03 & 2.00 \\ 
				1.32E-01 & 9.177E-03 & 2.04 & 3.094E-03 & 2.04 & 4.119E-03 & 1.18 \\ 
				9.90E-02 & 5.122E-03 & 2.04 & 1.742E-03 & 2.01 & 1.903E-03 & 2.70 \\
				\multicolumn{7}{c}{} \\
				\multicolumn{7}{c}{DG-IMEX-BDF $\mathcal{O}3$} \\
				\hline
				$h(\Omega)$ & $\rho_{L_1}$ & $\mathcal{O}(\rho_{L_1})$ & $\rho_{L_2}$ & $\mathcal{O}(\rho_{L_2})$ &  $\rho_{L_\infty}$ & $\rho_{L_\infty}$ \\ 
				\hline
				3.83E-01 & 4.897E-03 & -    & 1.215E-03 & -    & 1.217E-03 & -    \\ 
				1.96E-01 & 7.115E-04 & 2.88 & 2.023E-04 & 2.68 & 2.950E-04 & 2.12 \\ 
				1.32E-01 & 2.295E-04 & 2.84 & 6.666E-05 & 2.79 & 1.040E-04 & 2.62 \\ 
				9.90E-02 & 1.014E-04 & 2.86 & 2.941E-05 & 2.86 & 4.623E-05 & 2.84 \\ 
				\hline
				$h(\Omega)$ & $T_{L_1}$ & $\mathcal{O}(T_{L_1})$ & $T_{L_2}$ & $\mathcal{O}(T_{L_2})$ &  $T_{L_\infty}$ & $T_{L_\infty}$ \\ 
				\hline
				3.83E-01 & 1.301E-02 & -    & 3.644E-03 & -    & 3.517E-03 & -    \\ 
				1.96E-01 & 1.543E-03 & 3.18 & 4.626E-04 & 3.08 & 5.685E-04 & 2.72 \\ 
				1.32E-01 & 4.722E-04 & 2.98 & 1.411E-04 & 2.99 & 2.266E-04 & 2.31 \\ 
				9.90E-02 & 2.031E-04 & 2.95 & 6.065E-05 & 2.95 & 7.076E-05 & 4.07 \\
			\end{tabular}
		\end{small}
	\end{center}
	\label{tab.conv_0-BDF}
\end{table}

Next, a convergence study is performed for different values of the Knudsen number, namely $\varepsilon=10^{-3}$, $\varepsilon=10^{-2}$ and $\varepsilon=1$. Since no analytical solution is available in this case, the reference solution $U_{ref}(x,y)$, needed for computing the error norms, is obtained by running the test case on a very fine physical mesh while keeping the same discretization in the velocity space. Specifically, the reference triangular mesh is constructed by applying a refinement factor of $\xi=6$ to the finest grid used for the convergence analysis, thus subdividing each triangle into $36$ sub-triangles. Then, the numerical solution on each coarse grid is interpolated on the reference mesh and the errors can be computed. Tables \ref{tab.conv_Kn-RK} and \ref{tab.conv_Kn-BDF} collect the $L_1$ norms for density and temperature for the DG-IMEX-RK and DG-IMEX-BDF schemes, respectively. Convergence rates for second and third order of accuracy are observed for all three different Knudsen numbers, demonstrating that the space-time accuracy of the DG schemes is independent of the stiffness of the problem under consideration. 

\begin{table}[!htbp]  
	\caption{Numerical convergence results for the Boltzmann model using second and third order DG-IMEX-RK schemes at time $t_f=0.1$ with different Knudsen numbers ($\varepsilon=10^{-3}$, $\varepsilon=10^{-2}$, $\varepsilon=1$). A sequence of refined triangular meshes of size $h(\Omega)$ is used and the reference solution is computed on a fine mesh with refinement factor $\xi=6$. The errors are measured in $L_1$ norm and refer to the variables $\rho$ (density) and $T$ (temperature).}  
	\begin{center} 
		\begin{small}
			\renewcommand{\arraystretch}{1.0}
			\begin{tabular}{c|cccc|c} 
				\multicolumn{6}{c}{DG-IMEX-RK $\mathcal{O}2$} \\
				\hline
				$h(\Omega)$ & $\rho_{L_1}$ & $\mathcal{O}(\rho_{L_1})$ & $T_{L_1}$ & $\mathcal{O}(T_{L_1})$ &  $\varepsilon$ \\ 
				\hline
				5.75E-01 & 1.152E-01 & -    & 2.610E-01 & -    & 1E-03 \\ 
				2.87E-01 & 3.752E-02 & 1.62 & 7.891E-02 & 1.73 & 1E-03 \\ 
				1.92E-01 & 1.717E-02 & 1.93 & 3.268E-02 & 2.17 & 1E-03 \\ 
				1.44E-01 & 1.046E-02 & 1.72 & 1.985E-02 & 1.73 & 1E-03 \\ 
				\hline
				5.75E-01 & 1.116E-01 & -    & 2.525E-01 & -    & 1E-02 \\ 
				2.87E-01 & 3.474E-02 & 1.68 & 7.483E-02 & 1.75 & 1E-02 \\ 
				1.92E-01 & 1.531E-02 & 2.02 & 3.121E-02 & 2.16 & 1E-02 \\ 
				1.44E-01 & 9.136E-03 & 1.79 & 1.915E-02 & 1.70 & 1E-02 \\
				\hline
				5.75E-01 & 1.054E-01 & -    & 2.380E-01 & -    & 1E+00 \\ 
				2.87E-01 & 3.142E-02 & 1.75 & 6.809E-02 & 1.81 & 1E+00 \\ 
				1.92E-01 & 1.349E-02 & 2.09 & 2.841E-02 & 2.16 & 1E+00 \\ 
				1.44E-01 & 8.353E-03 & 1.67 & 1.806E-02 & 1.57 & 1E+00 \\
				\multicolumn{6}{c}{} \\
				\multicolumn{6}{c}{DG-IMEX-RK $\mathcal{O}3$} \\
				\hline
				$h(\Omega)$ & $\rho_{L_1}$ & $\mathcal{O}(\rho_{L_1})$ & $T_{L_1}$ & $\mathcal{O}(T_{L_1})$ &  $\varepsilon$ \\ 
				\hline
				5.75E-01 & 1.807E-02 & -    & 4.938E-02 & -    & 1E-03 \\ 
				2.87E-01 & 2.934E-03 & 2.62 & 6.675E-03 & 2.89 & 1E-03 \\ 
				1.92E-01 & 9.873E-04 & 2.69 & 1.941E-03 & 3.05 & 1E-03 \\ 
				1.44E-01 & 4.439E-04 & 2.78 & 8.192E-04 & 3.00 & 1E-03 \\ 
				\hline
				5.75E-01 & 1.740E-02 & -    & 4.831E-02 & -    & 1E-02 \\ 
				2.87E-01 & 2.600E-03 & 2.74 & 6.522E-03 & 2.89 & 1E-02 \\ 
				1.92E-01 & 8.016E-04 & 2.90 & 1.910E-04 & 3.03 & 1E-02 \\ 
				1.44E-01 & 3.427E-04 & 2.95 & 8.206E-04 & 2.94 & 1E-02 \\
				\hline
				5.75E-01 & 1.669E-02 & -    & 4.622E-02 & -    & 1E+00 \\ 
				2.87E-01 & 2.365E-03 & 2.82 & 6.118E-03 & 2.92 & 1E+00 \\ 
				1.92E-01 & 7.022E-04 & 3.00 & 1.778E-03 & 3.05 & 1E+00 \\ 
				1.44E-01 & 2.963E-04 & 3.00 & 7.581E-04 & 2.96 & 1E+00 \\
			\end{tabular}
		\end{small}
	\end{center}
	\label{tab.conv_Kn-RK}
\end{table}

\begin{table}[!htbp]  
	\caption{Numerical convergence results for the Boltzmann model using second and third order DG-IMEX-BDF schemes at time $t_f=0.1$ with different Knudsen numbers ($\varepsilon=10^{-3}$, $\varepsilon=10^{-2}$, $\varepsilon=1$). A sequence of refined triangular meshes of size $h(\Omega)$ is used and the reference solution is computed on a fine mesh with refinement factor $\xi=6$. The errors are measured in $L_1$ norm and refer to the variables $\rho$ (density) and $T$ (temperature).}  
	\begin{center} 
		\begin{small}
			\renewcommand{\arraystretch}{1.0}
			\begin{tabular}{c|cccc|c} 
				\multicolumn{6}{c}{DG-IMEX-BDF $\mathcal{O}2$} \\
				\hline
				$h(\Omega)$ & $\rho_{L_1}$ & $\mathcal{O}(\rho_{L_1})$ & $T_{L_1}$ & $\mathcal{O}(T_{L_1})$ &  $\varepsilon$ \\ 
				\hline
				5.75E-01 & 1.151E-01 & -    & 2.609E-01 & -    & 1E-03 \\ 
				2.87E-01 & 3.747E-02 & 1.62 & 7.885E-02 & 1.73 & 1E-03 \\ 
				1.92E-01 & 1.721E-02 & 1.92 & 3.263E-02 & 2.18 & 1E-03 \\ 
				1.44E-01 & 1.045E-02 & 1.73 & 1.985E-02 & 1.73 & 1E-03 \\ 
				\hline
				5.75E-01 & 1.114E-01 & -    & 2.522E-01 & -    & 1E-02 \\ 
				2.87E-01 & 3.469E-02 & 1.68 & 7.476E-02 & 1.75 & 1E-02 \\ 
				1.92E-01 & 1.533E-02 & 2.01 & 3.117E-02 & 2.16 & 1E-02 \\ 
				1.44E-01 & 9.135E-03 & 1.80 & 1.915E-02 & 1.69 & 1E-02 \\
				\hline
				5.75E-01 & 1.054E-01 & -    & 2.377E-01 & -    & 1E+00 \\ 
				2.87E-01 & 3.140E-02 & 1.75 & 6.802E-02 & 1.81 & 1E+00 \\ 
				1.92E-01 & 1.348E-02 & 2.08 & 2.834E-02 & 2.16 & 1E+00 \\ 
				1.44E-01 & 8.354E-03 & 1.66 & 1.807E-02 & 1.56 & 1E+00 \\
				\multicolumn{6}{c}{} \\
				\multicolumn{6}{c}{DG-IMEX-BDF $\mathcal{O}3$} \\
				\hline
				$h(\Omega)$ & $\rho_{L_1}$ & $\mathcal{O}(\rho_{L_1})$ & $T_{L_1}$ & $\mathcal{O}(T_{L_1})$ &  $\varepsilon$ \\ 
				\hline
				5.75E-01 & 1.826E-02 & -    & 4.955E-02 & -    & 1E-03 \\ 
				2.87E-01 & 2.928E-03 & 2.64 & 6.672E-03 & 2.89 & 1E-03 \\ 
				1.92E-01 & 9.870E-04 & 2.68 & 1.938E-03 & 3.05 & 1E-03 \\ 
				1.44E-01 & 4.427E-04 & 2.79 & 8.186E-04 & 3.00 & 1E-03 \\ 
				\hline
				5.75E-01 & 1.756E-02 & -    & 4.855E-02 & -    & 1E-02 \\ 
				2.87E-01 & 2.600E-03 & 2.76 & 6.522E-03 & 2.90 & 1E-02 \\ 
				1.92E-01 & 8.022E-04 & 2.90 & 1.914E-04 & 3.02 & 1E-02 \\ 
				1.44E-01 & 3.428E-04 & 2.96 & 8.237E-04 & 2.93 & 1E-02 \\
				\hline
				5.75E-01 & 1.680E-02 & -    & 4.622E-02 & -    & 1E+00 \\ 
				2.87E-01 & 2.365E-03 & 2.83 & 6.133E-03 & 2.93 & 1E+00 \\ 
				1.92E-01 & 7.016E-04 & 3.00 & 1.835E-03 & 2.98 & 1E+00 \\ 
				1.44E-01 & 2.961E-04 & 3.00 & 7.911E-04 & 2.93 & 1E+00 \\
			\end{tabular}
		\end{small}
	\end{center}
	\label{tab.conv_Kn-BDF}
\end{table}

%
\subsection{Lax problem} 
\label{sec.Lax}
\begin{figure}[!htbp]
	\begin{center}
		\begin{tabular}{ccc} 
			\includegraphics[width=0.33\textwidth]{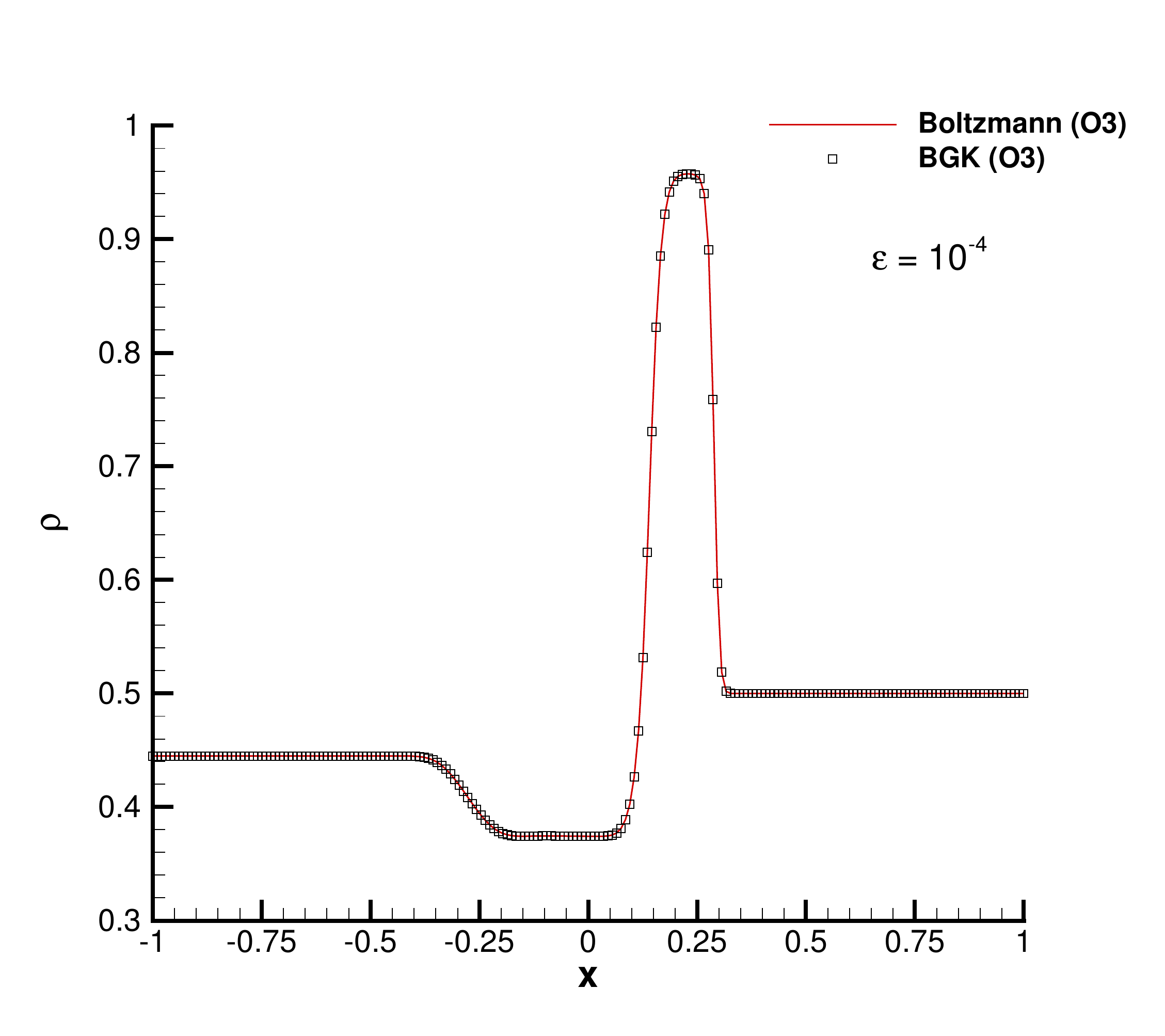}  &           
			\includegraphics[width=0.33\textwidth]{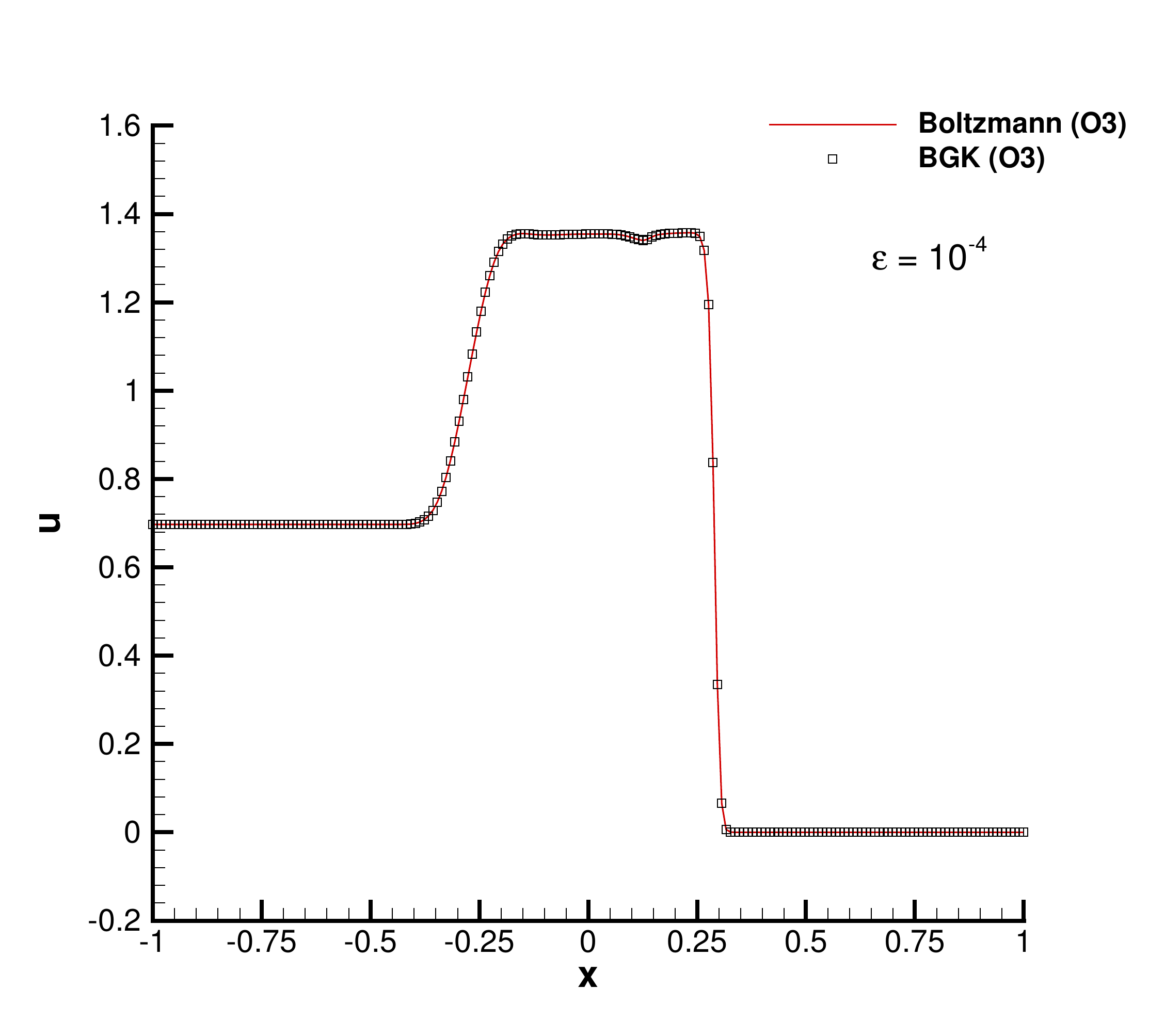}  &           
			\includegraphics[width=0.33\textwidth]{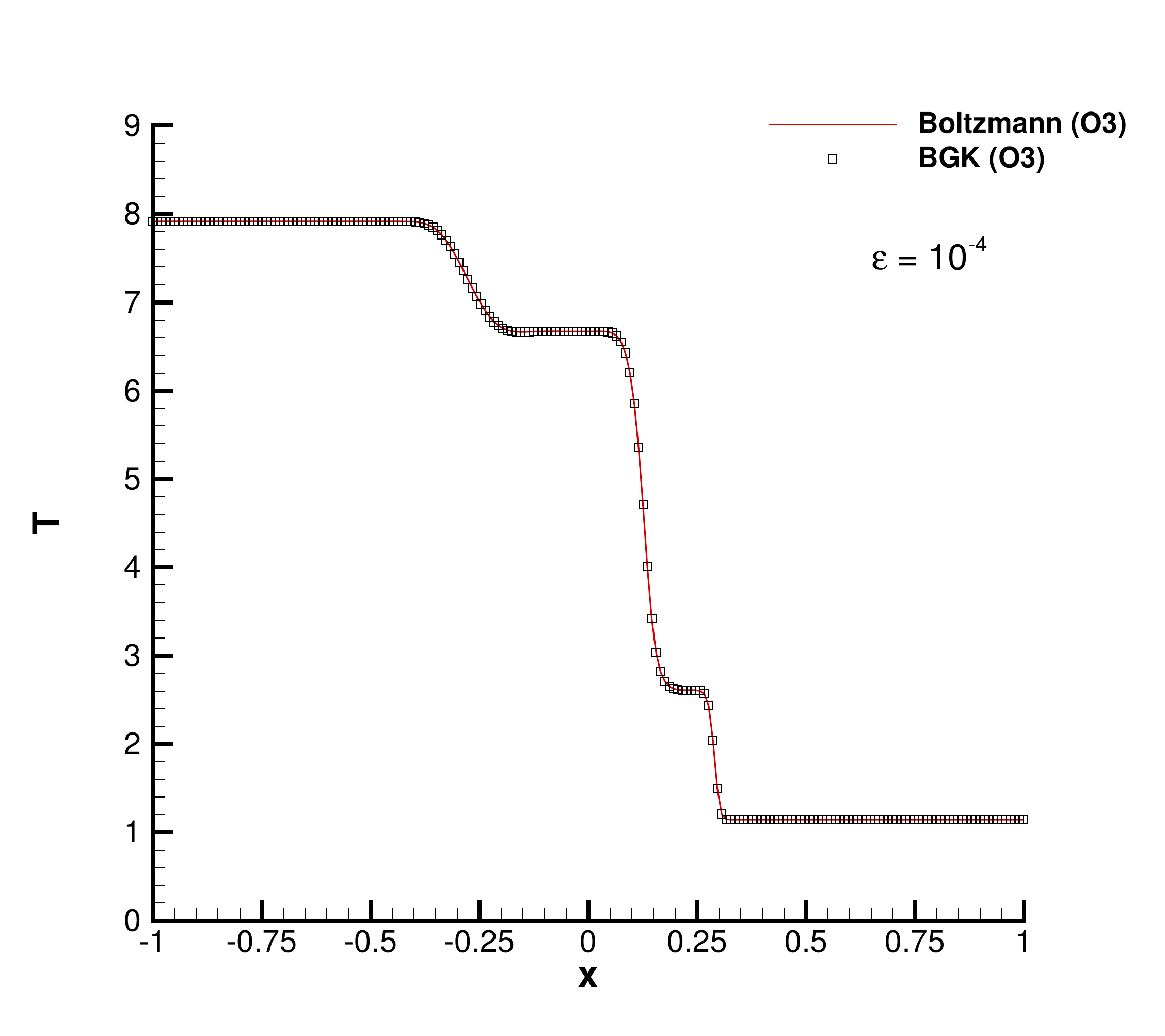} \\
			\includegraphics[width=0.33\textwidth]{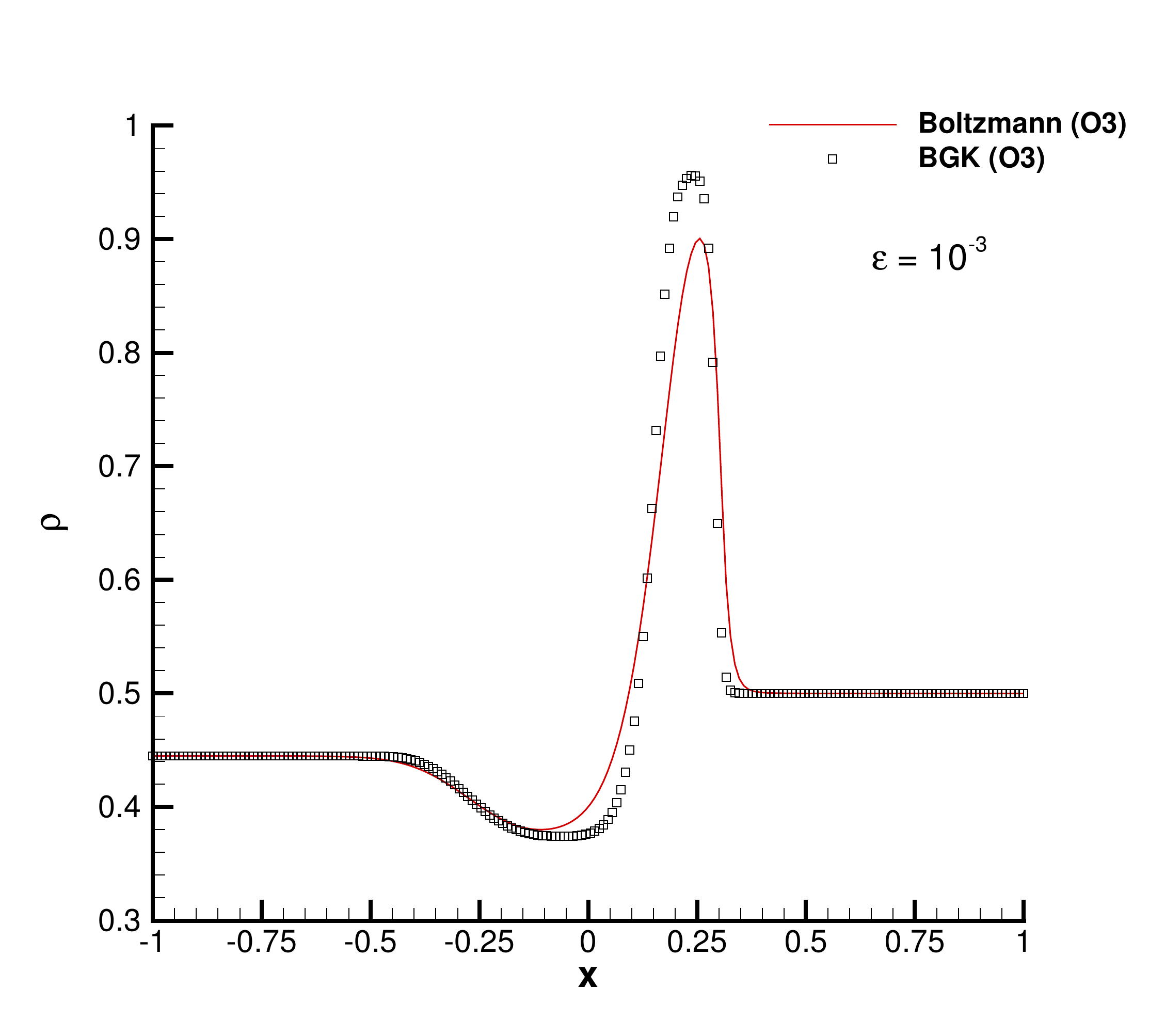}  &           
			\includegraphics[width=0.33\textwidth]{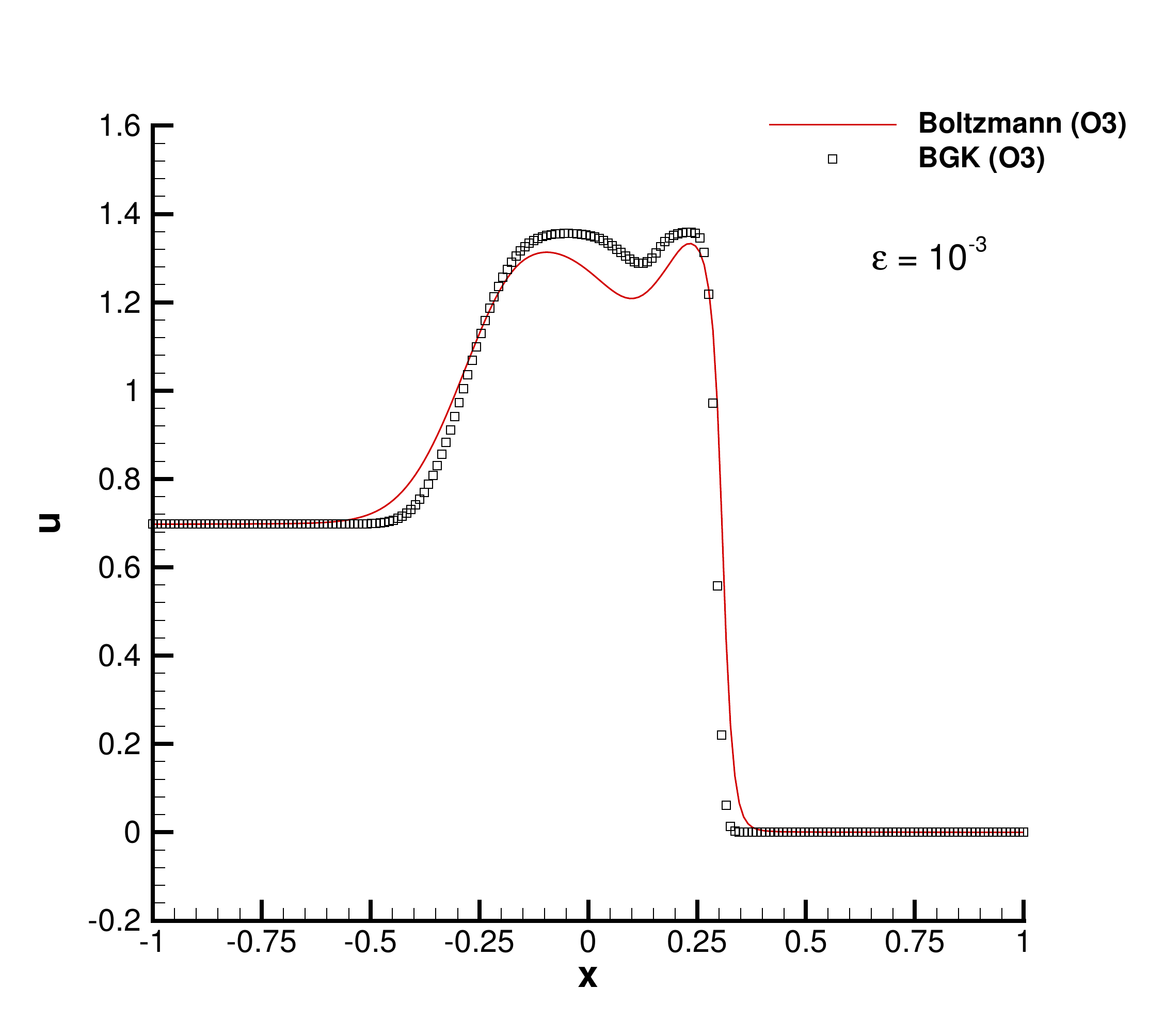}  &           
			\includegraphics[width=0.33\textwidth]{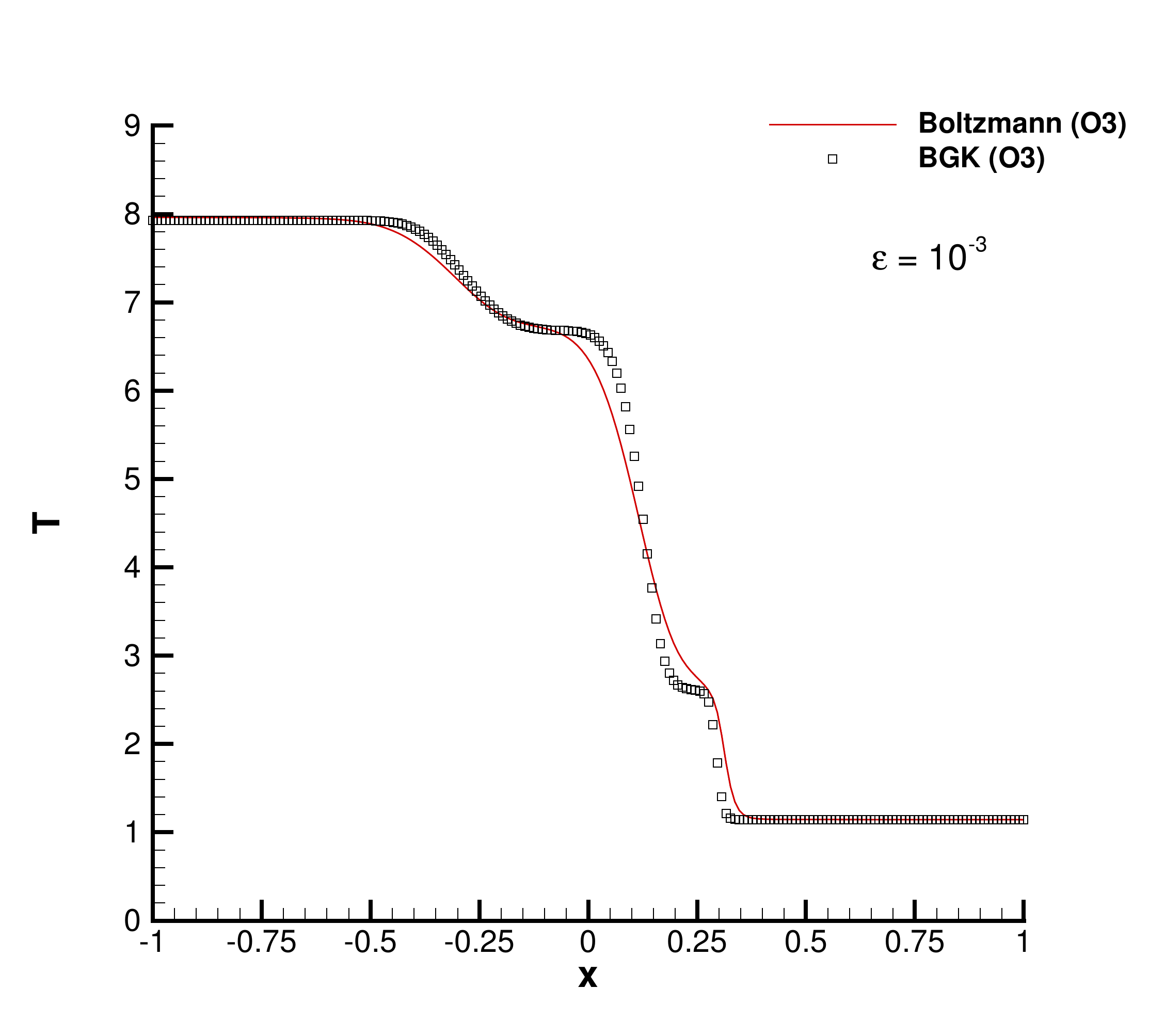} \\
			\includegraphics[width=0.33\textwidth]{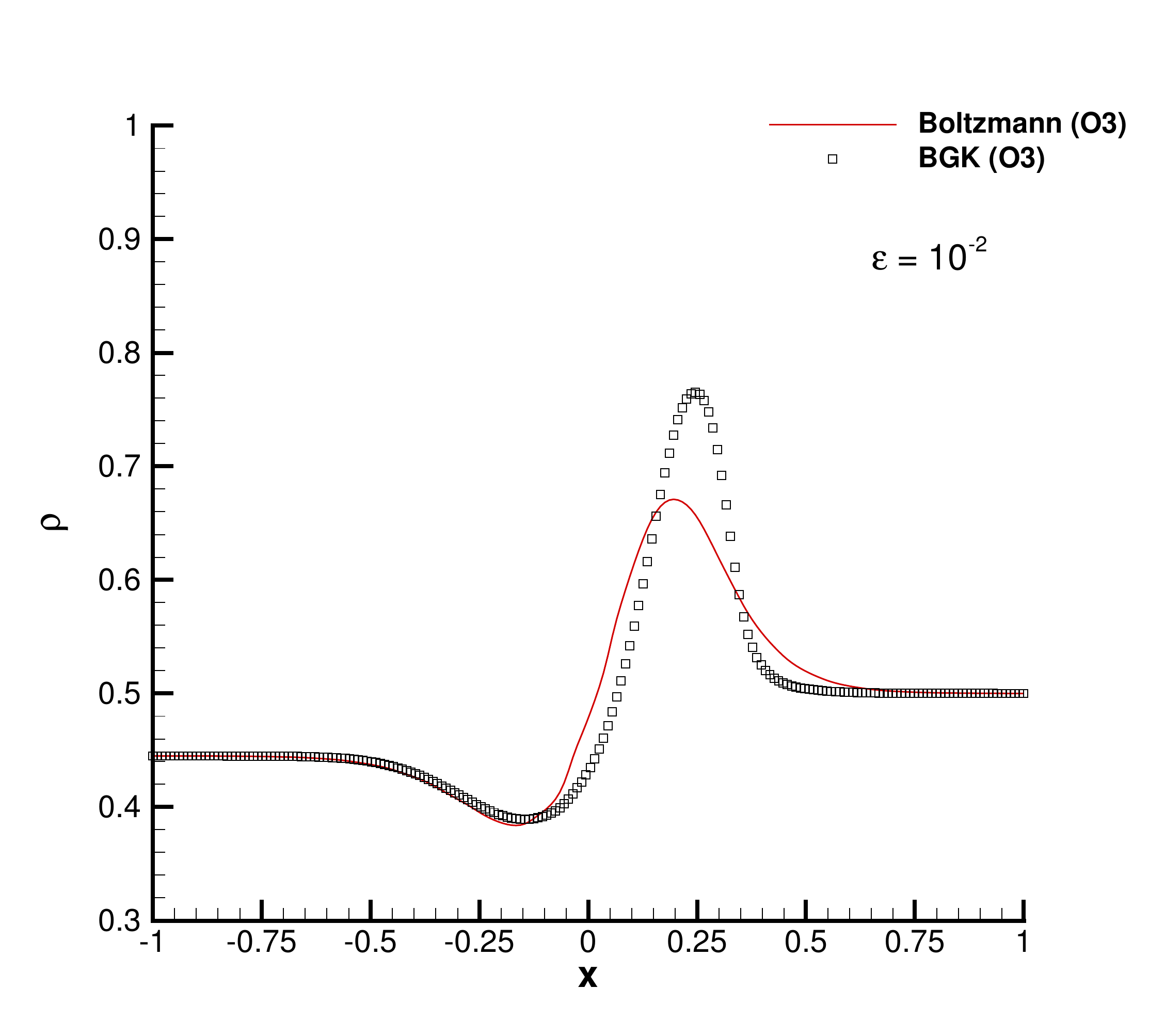}  &           
			\includegraphics[width=0.33\textwidth]{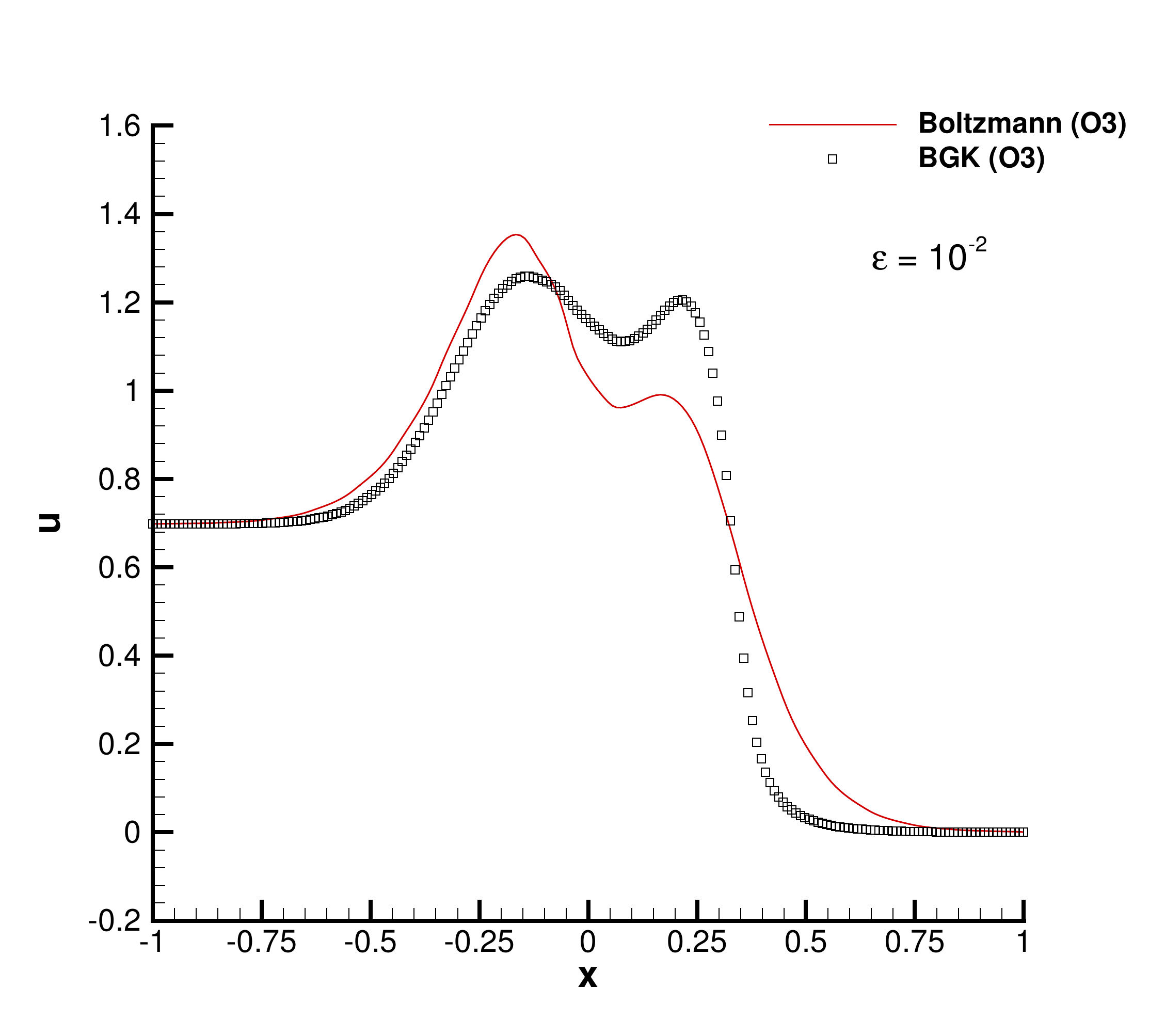}  &           
			\includegraphics[width=0.33\textwidth]{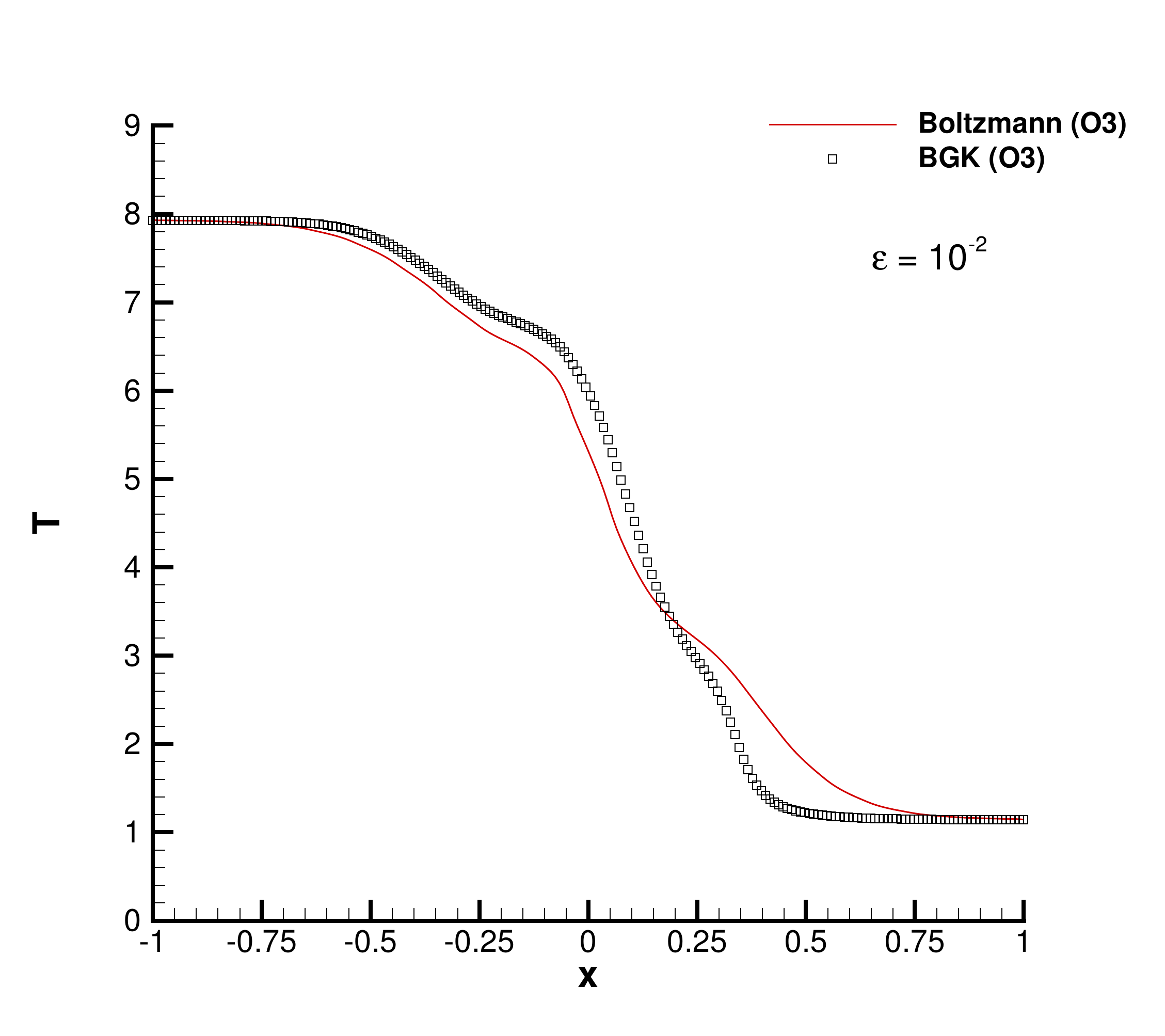} \\
		\end{tabular}
		\caption{Lax shock tube problem with $\varepsilon=10^{-4}$ (top row), $\varepsilon=10^{-3}$ (middle row) and $\varepsilon=10^{-2}$ (bottom row) for Boltzmann (red solid line) and BGK (black square scatters) model. 1D cut along the $x$-axis through the third order numerical results and comparison with exact solution for density, velocity in the $x$-direction and temperature.}
		\label{fig.Lax}
	\end{center}
\end{figure}

To verify the capability of the DG method to deal with classical fluid dynamic problems eventually involving shocks, the Lax shock tube Riemann problem is proposed hereafter. The computational domain is given by the rectangle $\Omega=[-1;1] \times [-0.05;0.05]$, with periodic boundary conditions imposed in the $y-$direction, while Dirichlet boundaries set along the $x-$direction. A total number of $[200 \times 10]$ triangular control volumes with characteristic mesh size of $h(\Omega)=0.01$ are used to discretize the physical space, whereas the velocity space $\mathcal{V}=[-15;15]\times [-15;15]$ counts a total number of $1024$ regular Cartesian control volumes giving rise to a problem with approximately $2\cdot 10^6$ mesh points. Although the set-up is typical of a one-dimensional problem, on unstructured meshes no element edges are aligned with the fluid motion, thus the capability of the method to maintain the one-dimensional structure of the solution can also be analyzed by running such problems. The initial condition is given in terms of two states $(U_L,U_R)$ separated at $x=0$, hence
\begin{equation}
	 U_L= \left(0.445,0.698,0,7.928 \right), \qquad U_R = \left(0.5,0,0,1.142\right).
\end{equation}
This test case is run with the third order DG-IMEX-BDF schemes considering three different values of the Knudsen number, that is $\varepsilon=10^{-4}$, $\varepsilon=10^{-3}$ and $\varepsilon=10^{-2}$. Furthermore, a comparison with the BGK model \eqref{BGK} is also proposed. Figure \ref{fig.Lax} shows a one-dimensional cut along the $x-$axis for the density, the horizontal velocity and the temperature profiles for the Boltzmann and BGK models for each value of the Knudsen number. The smaller is $\varepsilon$, the closer is the matching between the results for Boltzmann and BGK models, because the limit of the compressible Euler equations is approached. On the contrary, as the Knudsen number increases, differences between the two models arise since the source terms in the governing equations are responsible for a different modelling of the collisions among the gas particles.

\begin{figure}[!htbp]
	\begin{center}
		\begin{tabular}{ccc} 
			\includegraphics[width=0.33\textwidth]{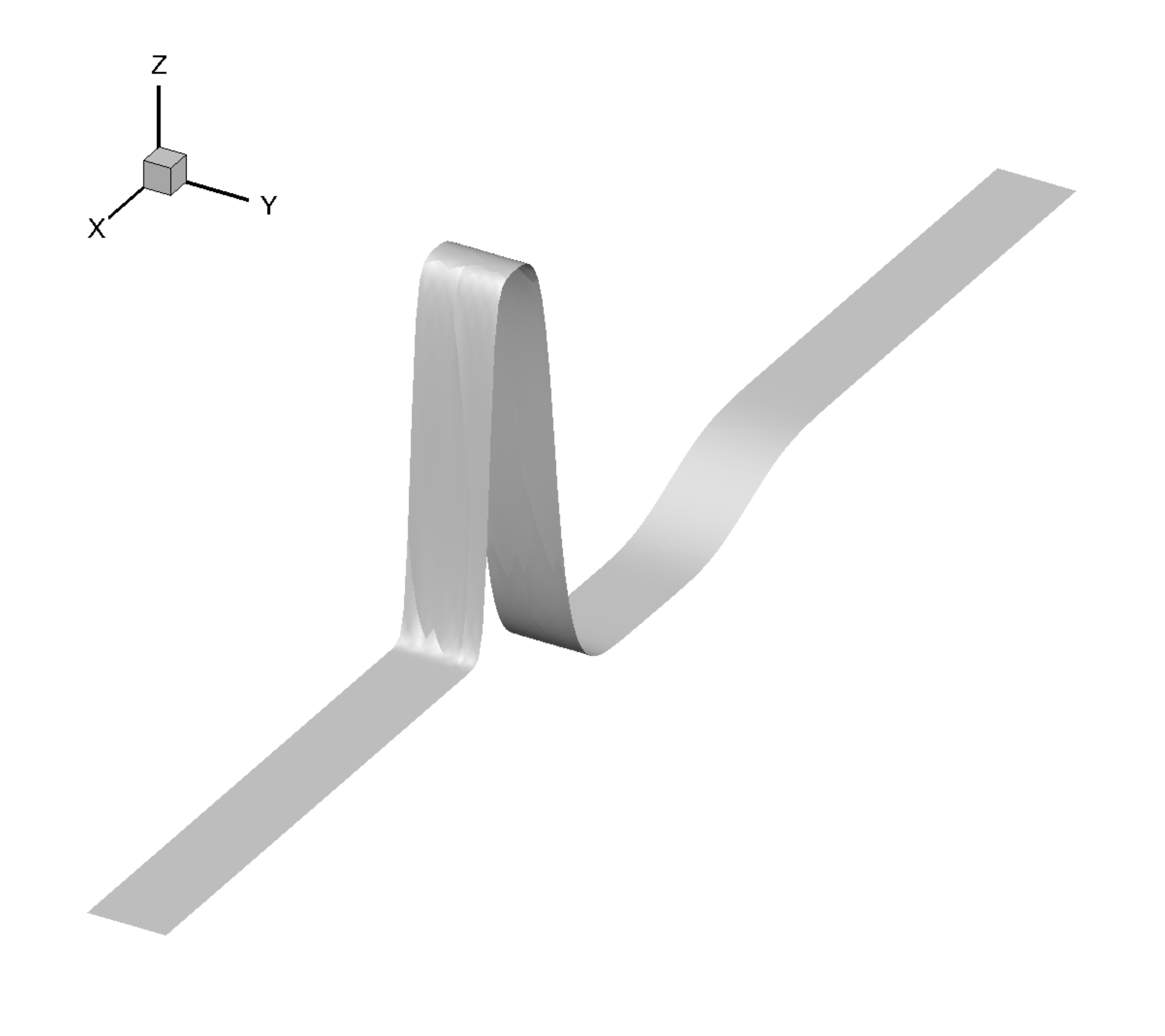}  &           
			\includegraphics[width=0.33\textwidth]{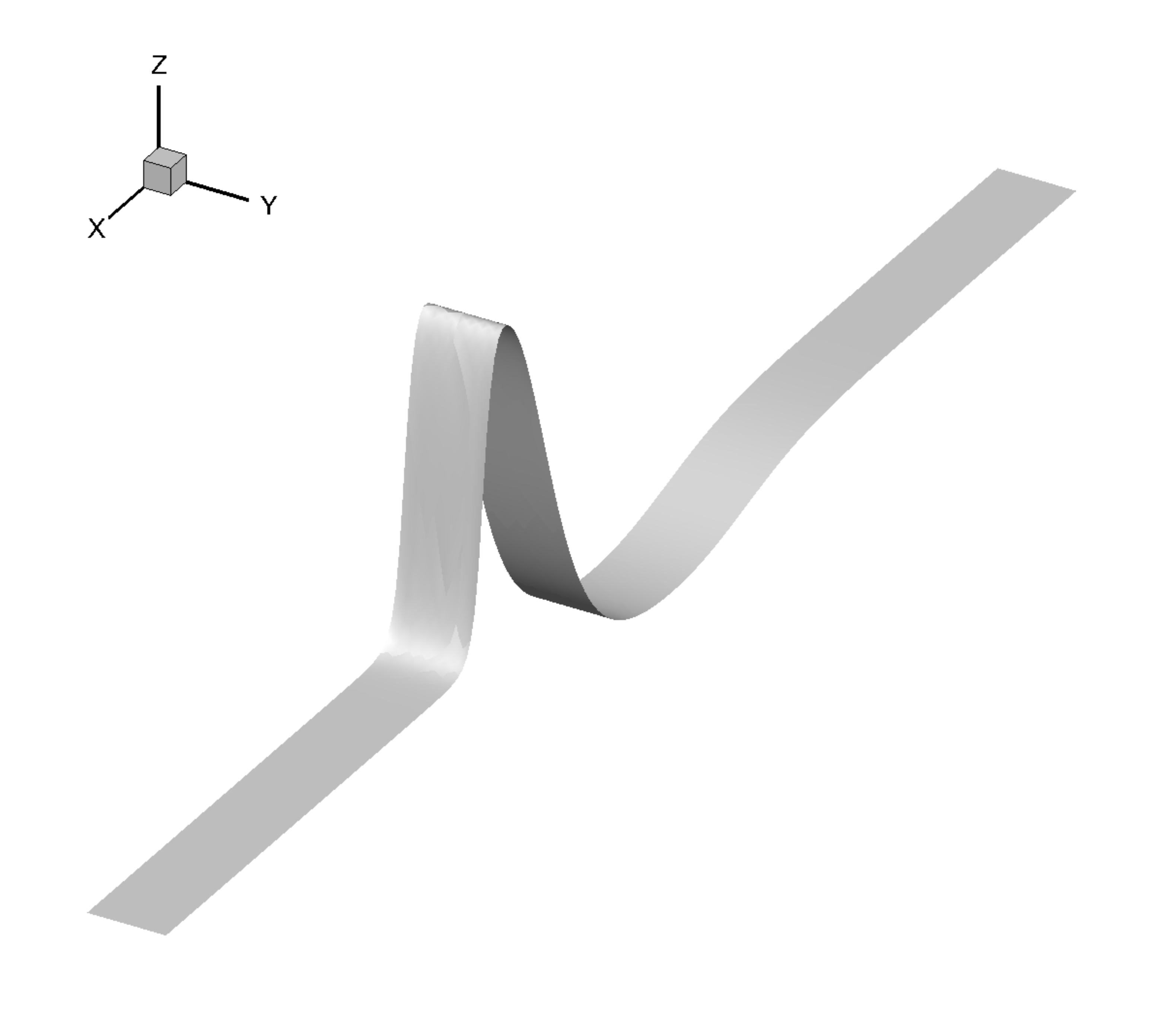}  &           
			\includegraphics[width=0.33\textwidth]{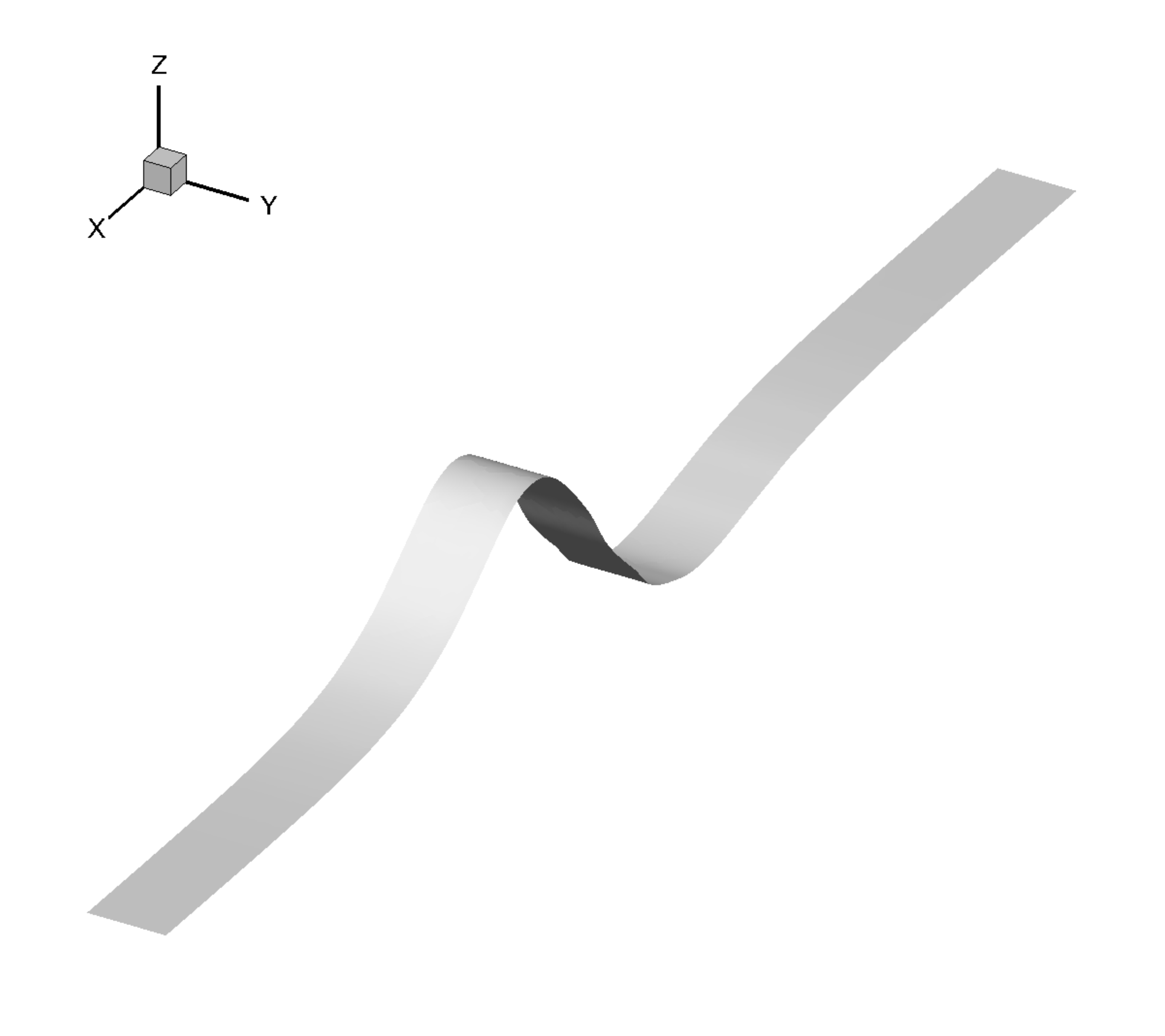} \\
		\end{tabular}
		\caption{Lax shock tube problem. Three dimensional view of density profile for the Boltzmann model with $\varepsilon=10^{-4}$ (left), $\varepsilon=10^{-3}$ (middle) and $\varepsilon=10^{-2}$ (right).}
		\label{fig.Lax-rho3D}
	\end{center}
\end{figure}

\begin{figure}[!htbp]
	\begin{center}
		\begin{tabular}{ccc} 
			\includegraphics[width=0.33\textwidth]{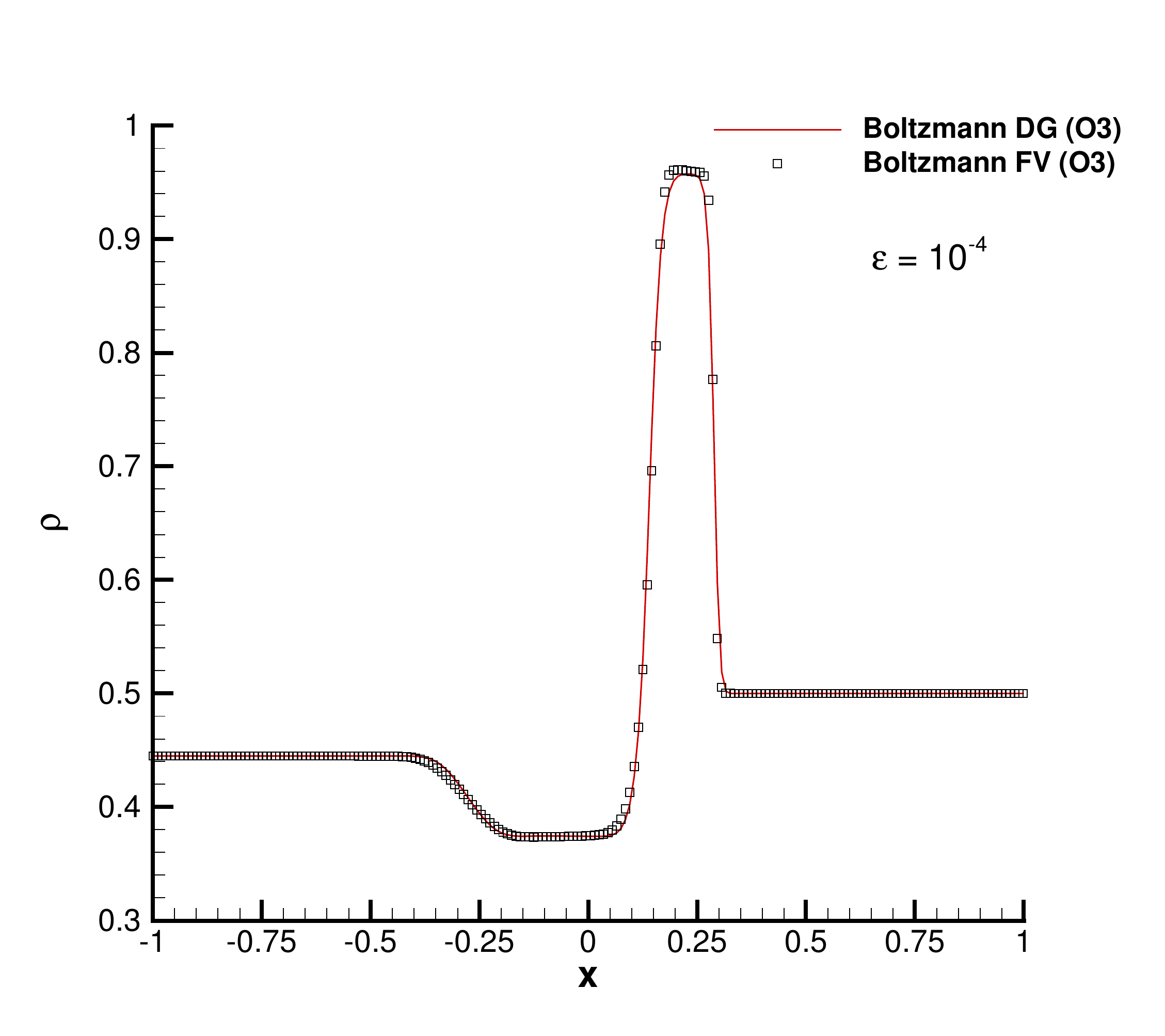}  &           
			\includegraphics[width=0.33\textwidth]{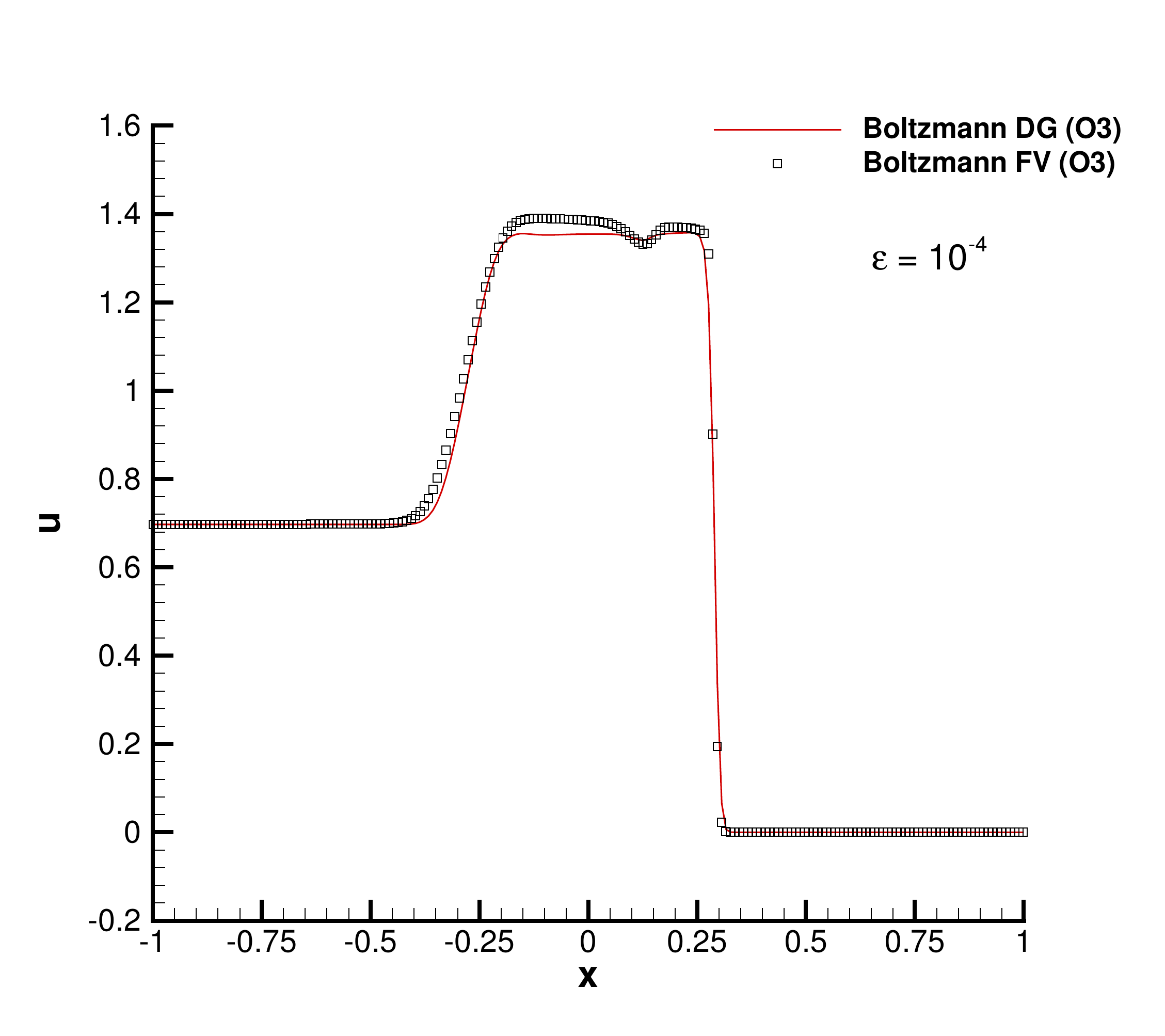}  &           
			\includegraphics[width=0.33\textwidth]{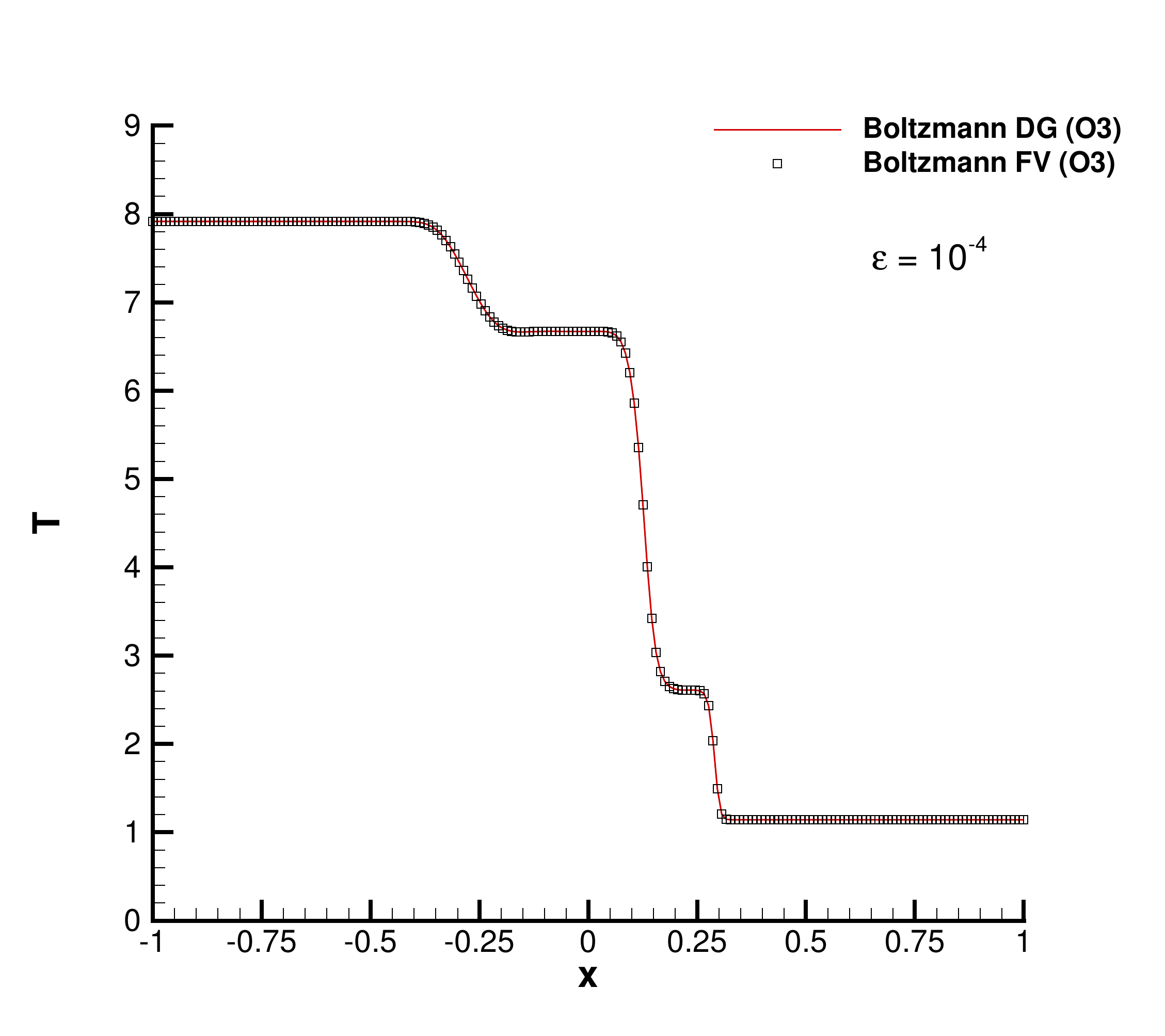} \\
				\includegraphics[width=0.33\textwidth]{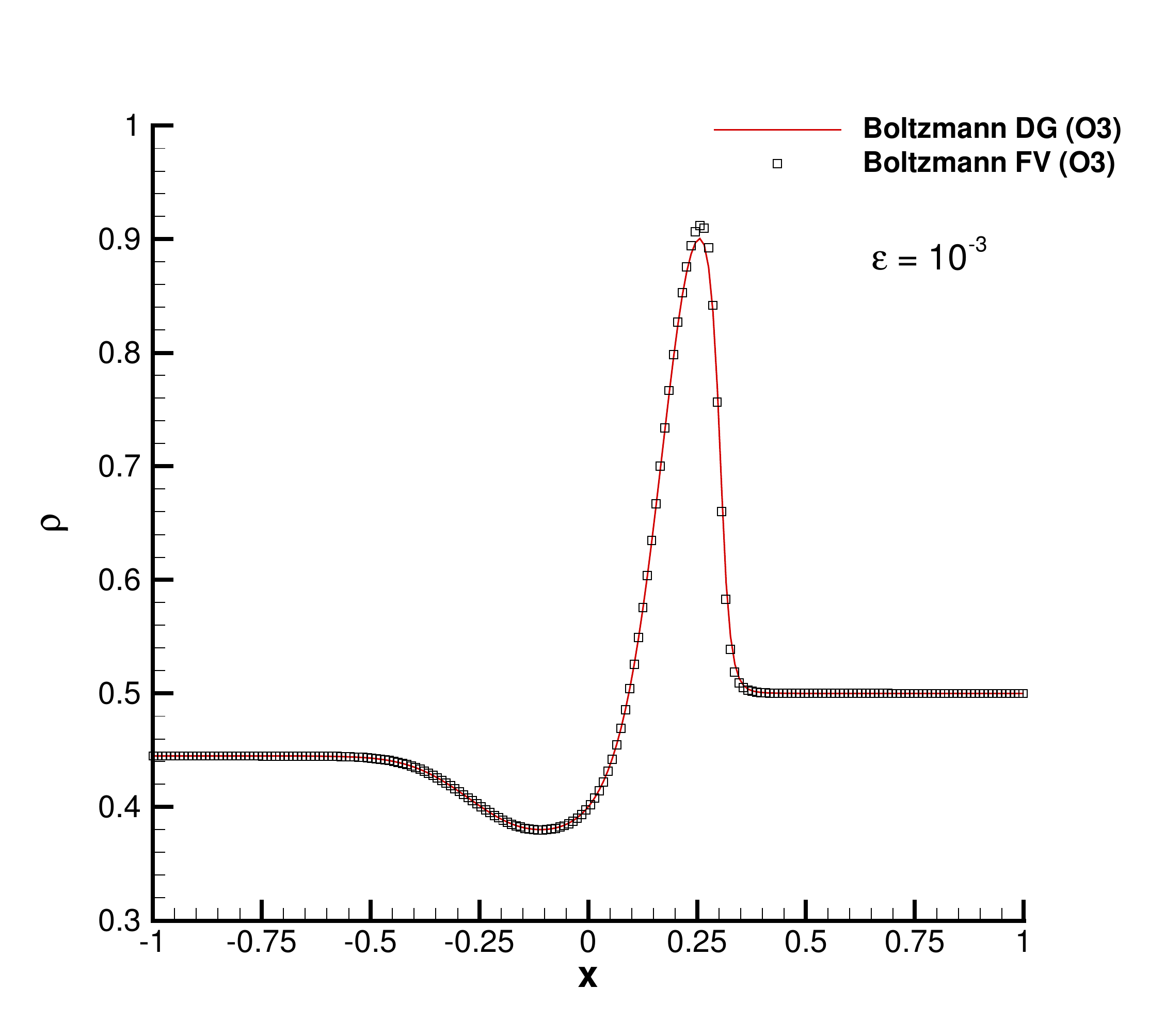}  &           
			\includegraphics[width=0.33\textwidth]{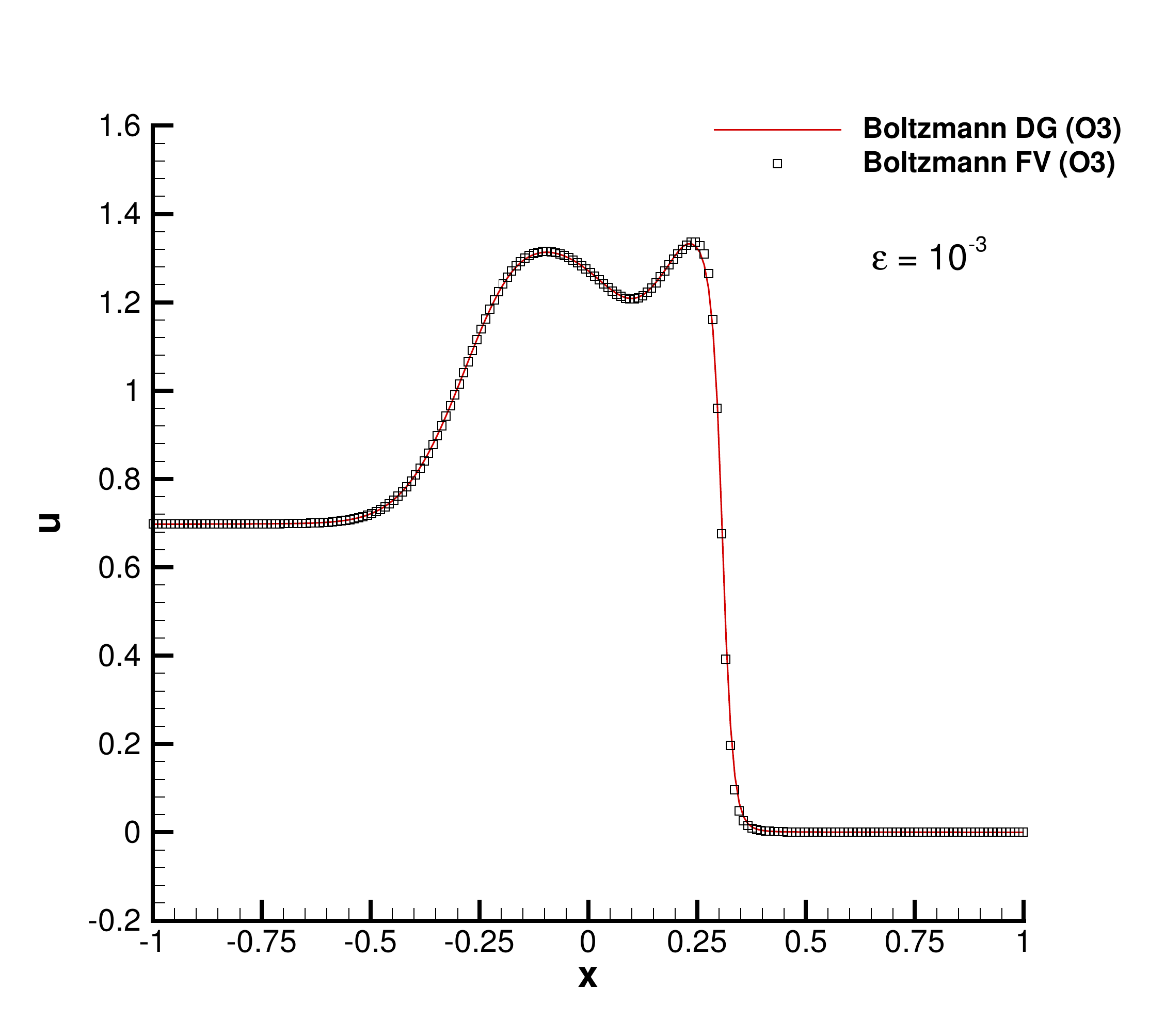}  &           
			\includegraphics[width=0.33\textwidth]{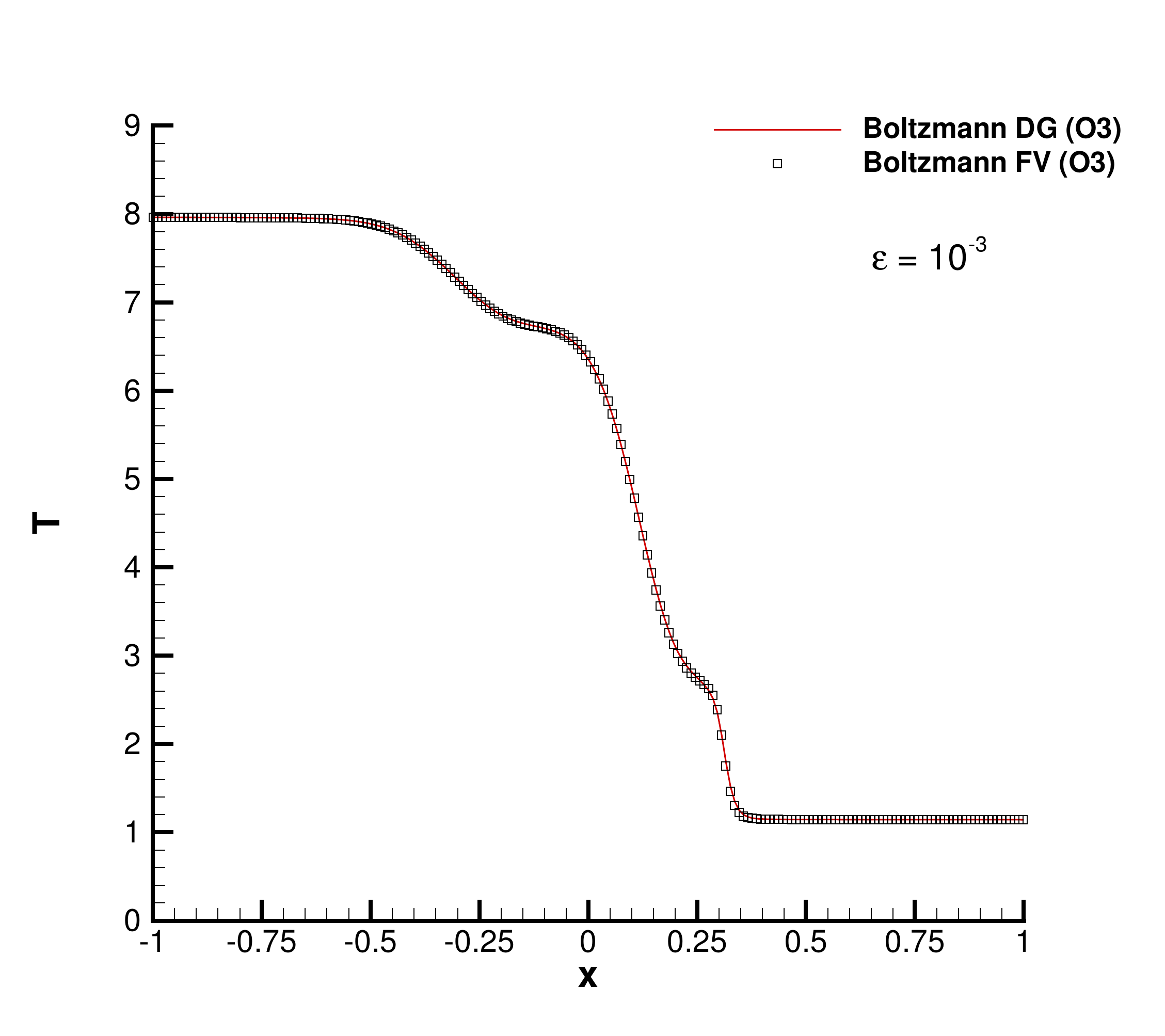} \\
				\includegraphics[width=0.33\textwidth]{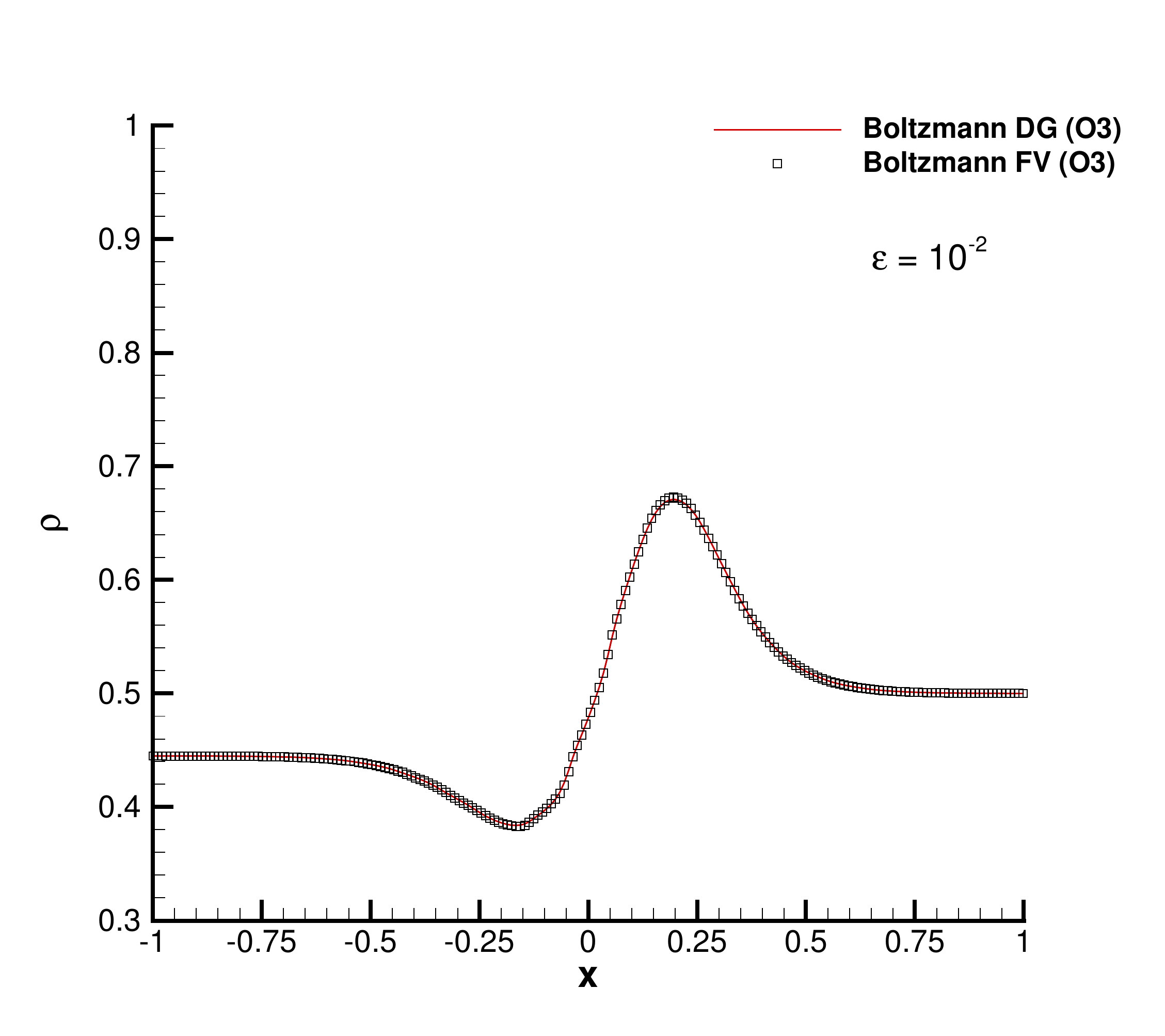}  &           
			\includegraphics[width=0.33\textwidth]{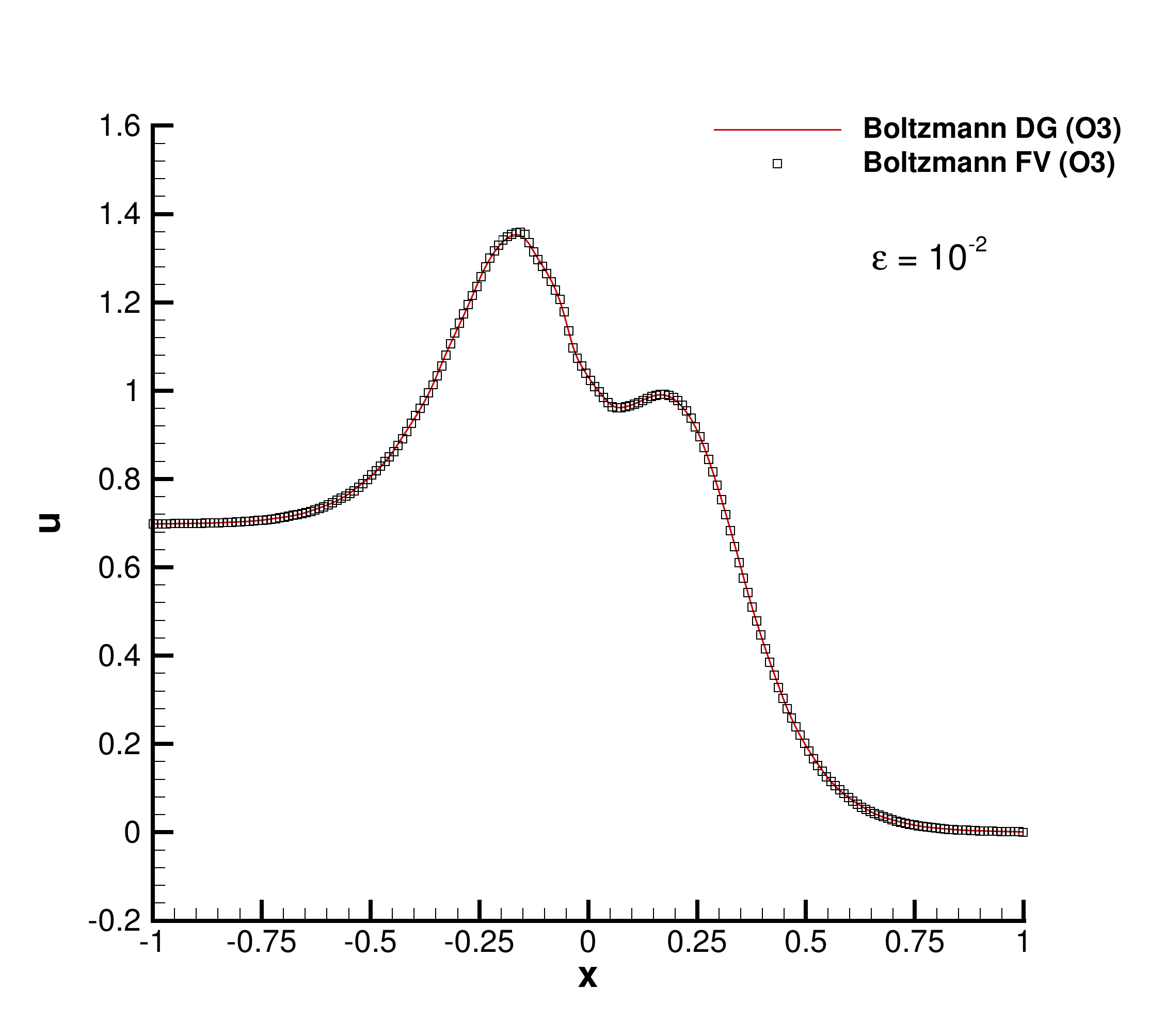}  &           
			\includegraphics[width=0.33\textwidth]{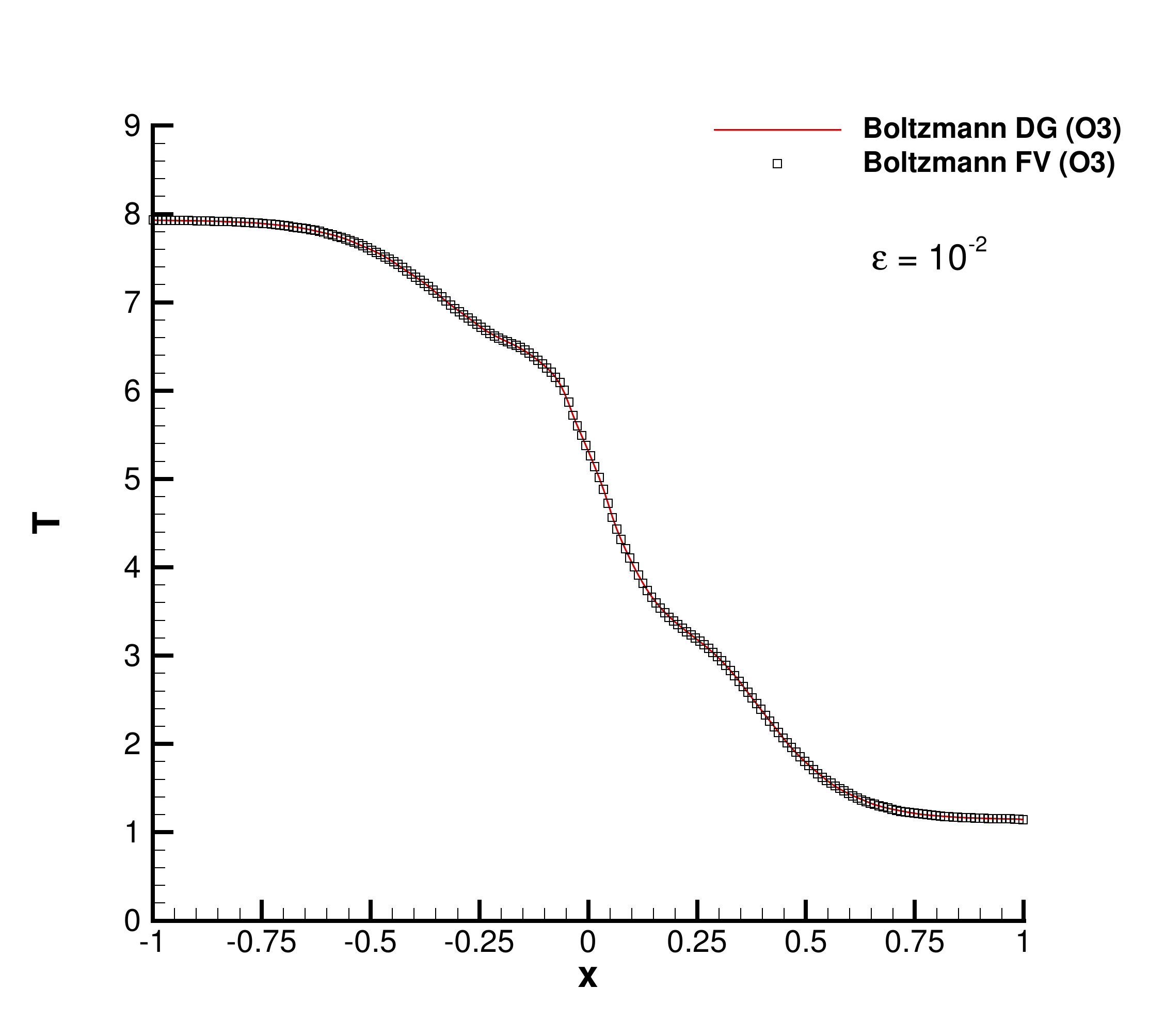} \\
		\end{tabular}
		\caption{Lax shock tube problem with $\varepsilon=10^{-4}$ (top row), $\varepsilon=10^{-3}$ (middle row) and $\varepsilon=10^{-2}$ (bottom row) for DG-IMEX-BDF (red solid line) and FV-IMEX-BDF (black square scatters) third order schemes. 1D cut along the $x$-axis through the third order numerical results and comparison with exact solution for density, velocity in the $x$-direction and temperature.}
		\label{fig.Lax_FV-DG}
	\end{center}
\end{figure}

Figure \ref{fig.Lax-rho3D} depicts instead a three-dimensional view of the density profile for all three simulations with different Knudsen number for the Boltzmann model, demonstrating that the solution maintains a one-dimensional structure, despite the unstructured meshes used to pave the computational domain.
Finally, the DG solution is compared against the results obtained with the third order IMEX-BDF finite volume method presented in \cite{BosDim2}, showing excellent agreement in Figure \ref{fig.Lax_FV-DG} for all macroscopic quantities and all the three values of the Knudsen number.

%
\subsection{Explosion problem}
\begin{figure}[!htbp]
	\begin{center}
		\begin{tabular}{cc} 
			\includegraphics[width=0.47\textwidth]{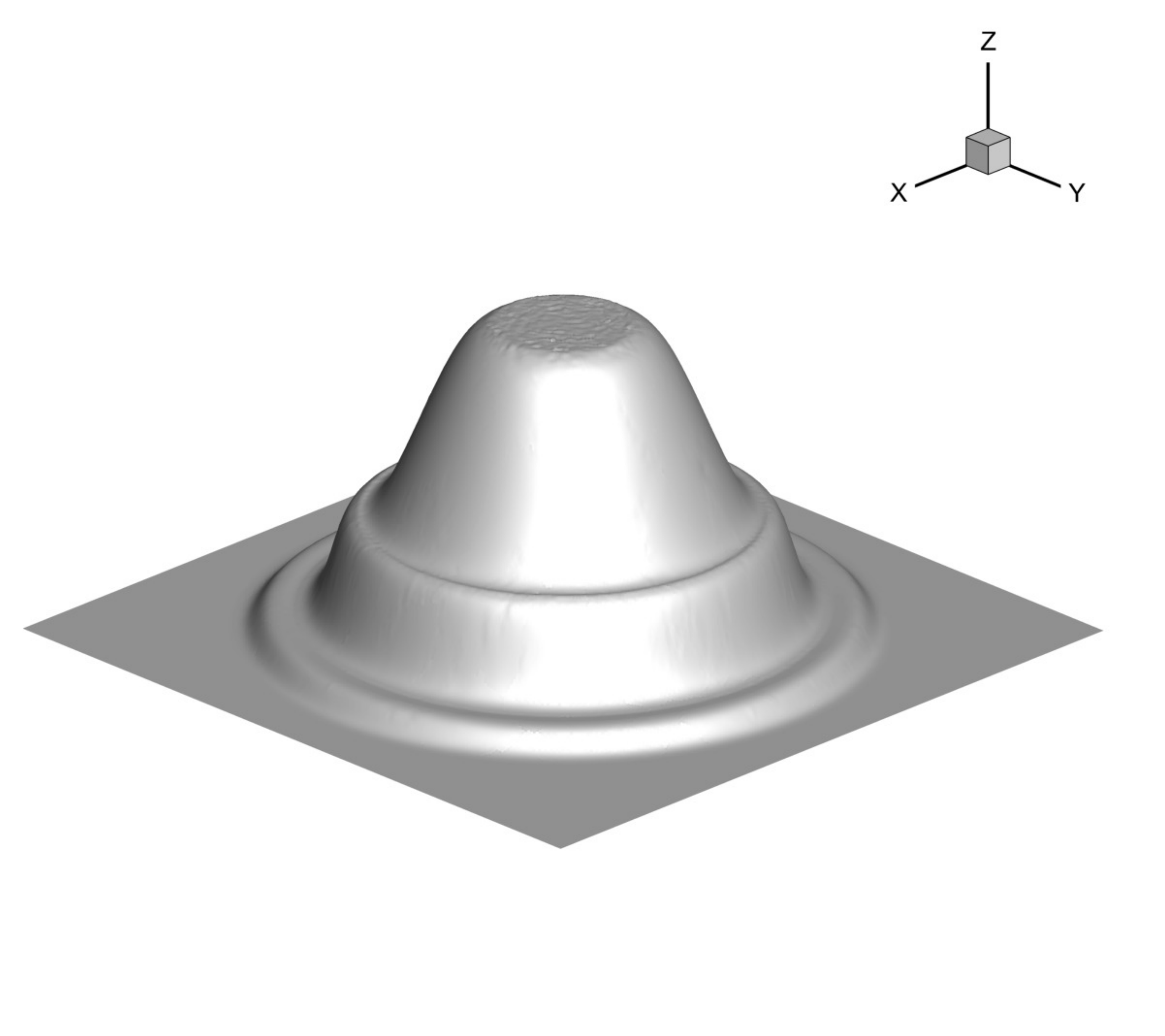}  &           
			\includegraphics[width=0.47\textwidth]{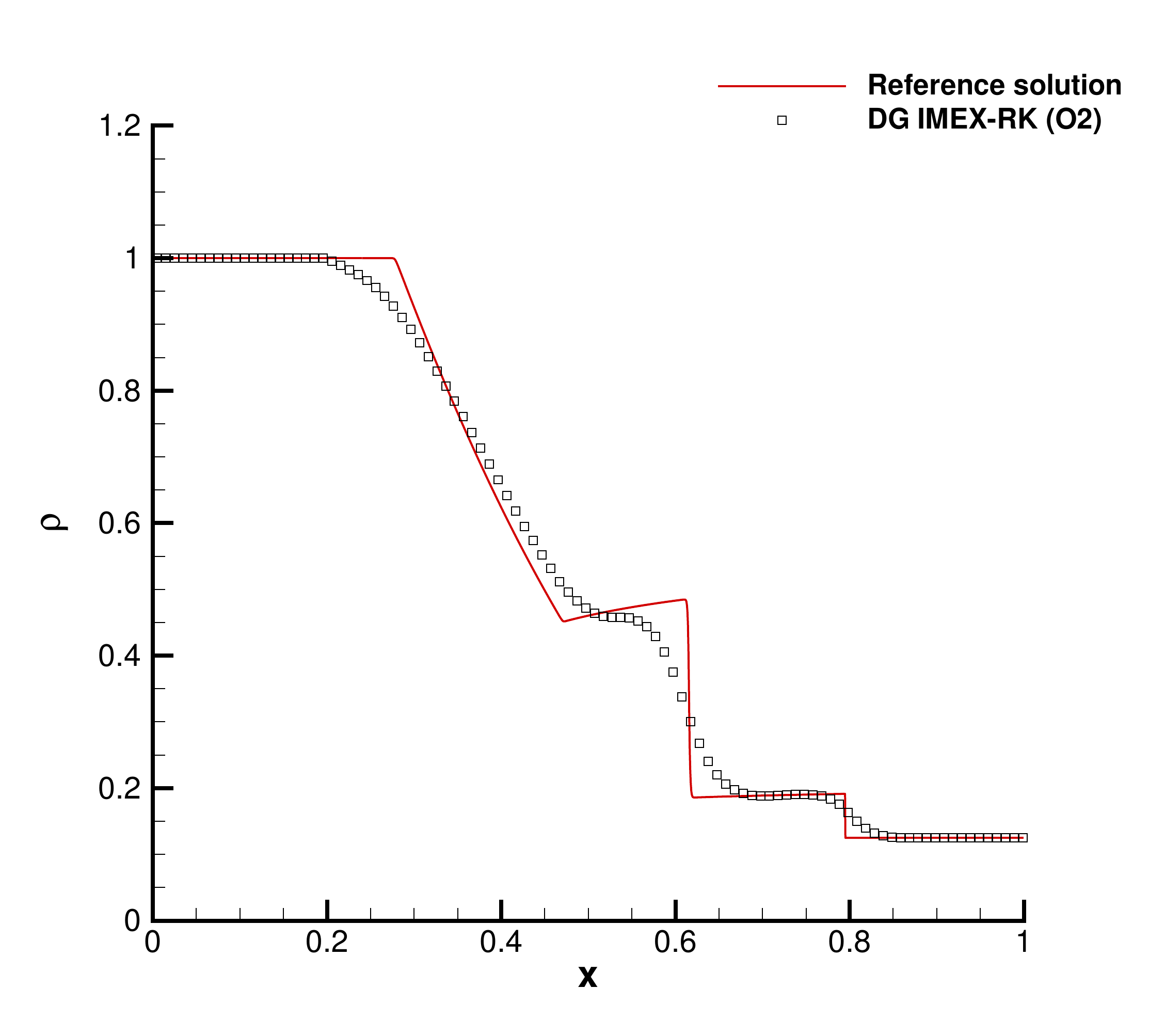}  \\           
			\includegraphics[width=0.47\textwidth]{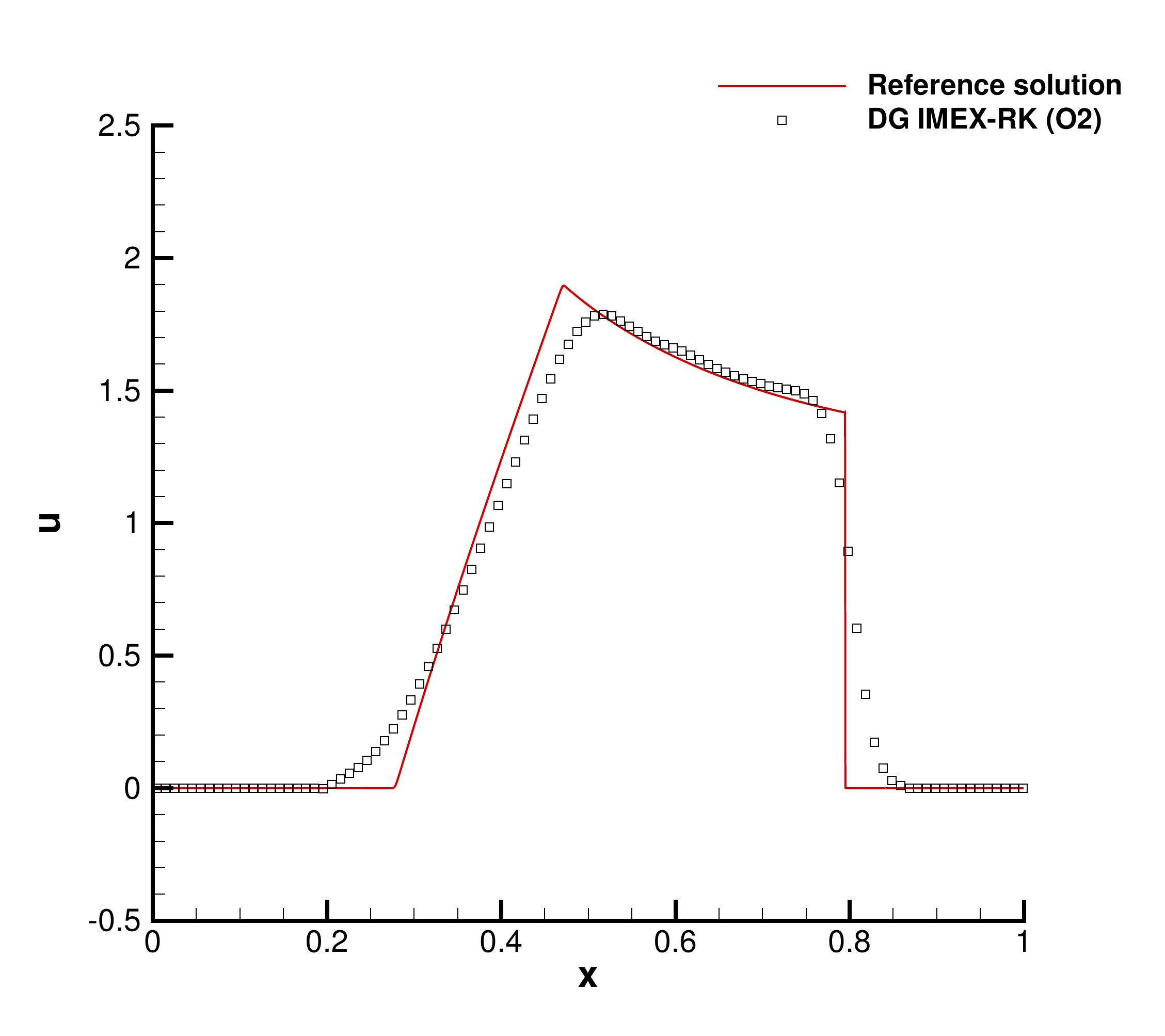}  &           
			\includegraphics[width=0.47\textwidth]{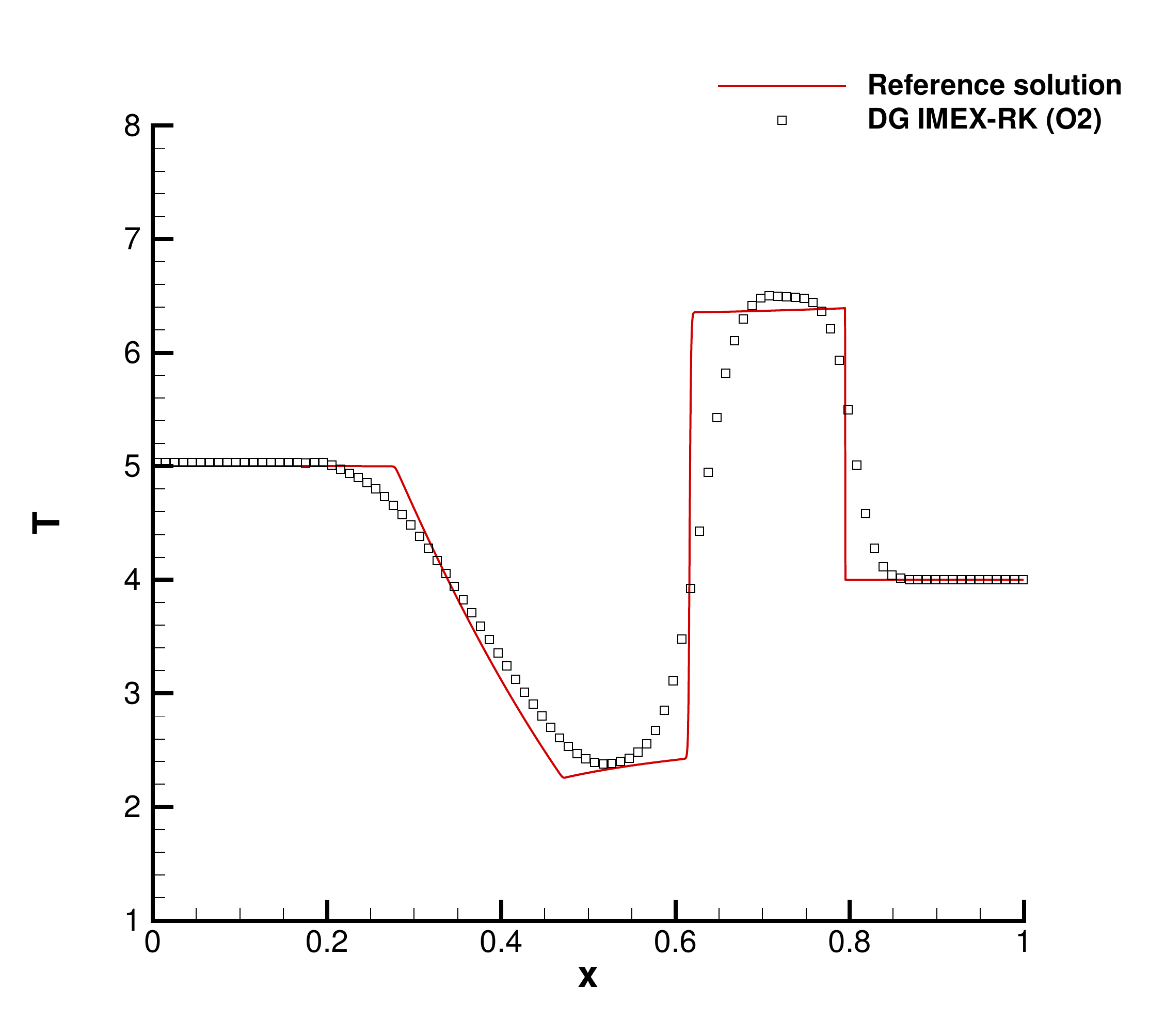}  \\ 
		\end{tabular}
		\caption{Explosion problem with $\varepsilon=10^{-5}$. Three-dimensional view of density profile for the Boltzmann model (top left). 1D cut along the $x$-axis through the second order numerical results and comparison with exact solution for density, velocity in the $x$-direction and temperature (from top to bottom right).}
		\label{fig.EP3D}
	\end{center}
\end{figure}

The explosion problem is often regarded as a multidimensional extension of the classical Sod shock tube test case. The physical space is given by the square $\Omega=[-1;1]\times[-1;1]$ with Dirichlet boundary conditions, which is discretized by $200 \times 200$ polygonal cells. The velocity space is defined in the interval $\mathcal{V}=[-15;15]\times [-15;15]$. The inner and outer states are initially prescribed, that are separated by a discontinuity at distance $R_d=0.5$ from the origin:
\begin{equation}
	\begin{cases} 
		U_{in} = \left(1,0,0,1 \right) & R \leq R_d, \\
		U_{out} = \left(0.125,0,0,0.8\right) & R > R_d,
	\end{cases}
\end{equation}
with $R=\sqrt{x^2+y^2}$ representing the radial position. The final time of the simulation is chosen to be $t_f=0.07$ and the second order DG-IMEX schemes with BDF and RK time stepping are used to run the simulation with Knudsen number $\varepsilon=10^{-5}$. Figure \ref{fig.EP3D} shows a three-dimensional view of the density profile for the DG-IMEX-RK at the final time as well as a comparison against the reference solution, given by the compressible Euler equations \cite{ToroBook}, for density, horizontal velocity and temperature. An overall good agreement can be noticed, and the DG scheme correctly captures the position of the shock wave and the plateau exhibited by the solution. 

\begin{figure}[!htbp]
	\begin{center}
		\begin{tabular}{ccc} 
			\includegraphics[width=0.33\textwidth]{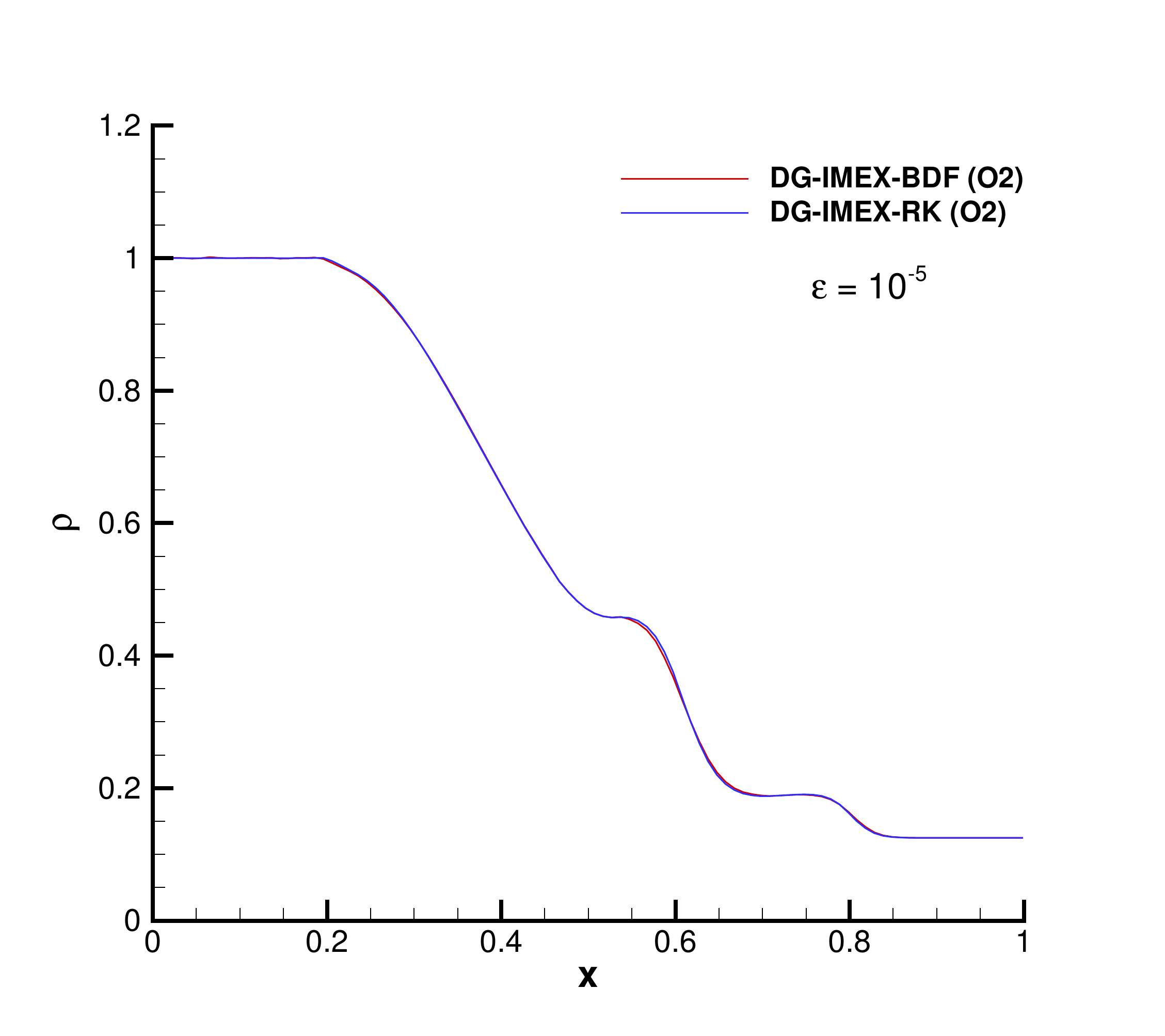}  &            
			\includegraphics[width=0.33\textwidth]{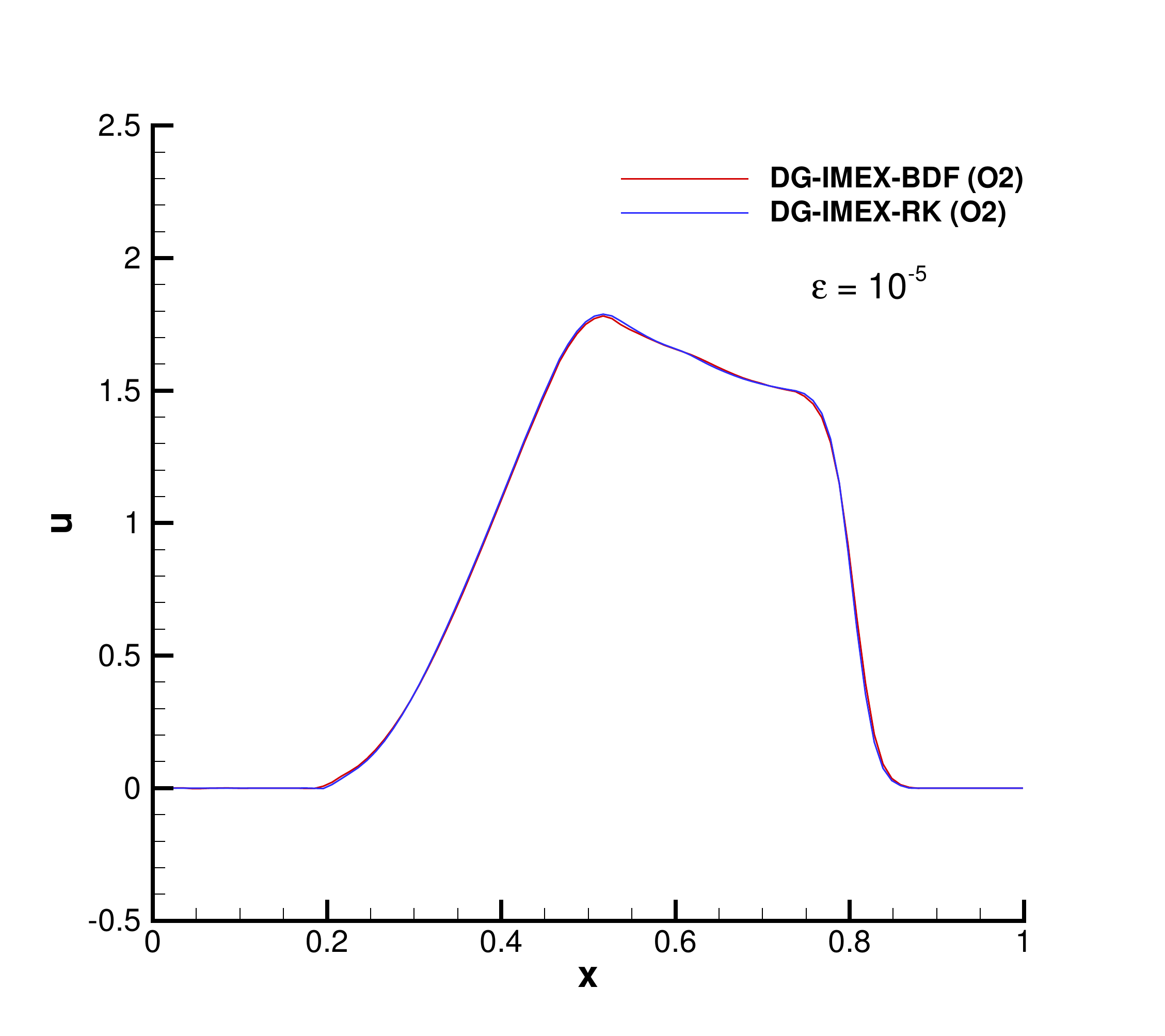}  &           
			\includegraphics[width=0.33\textwidth]{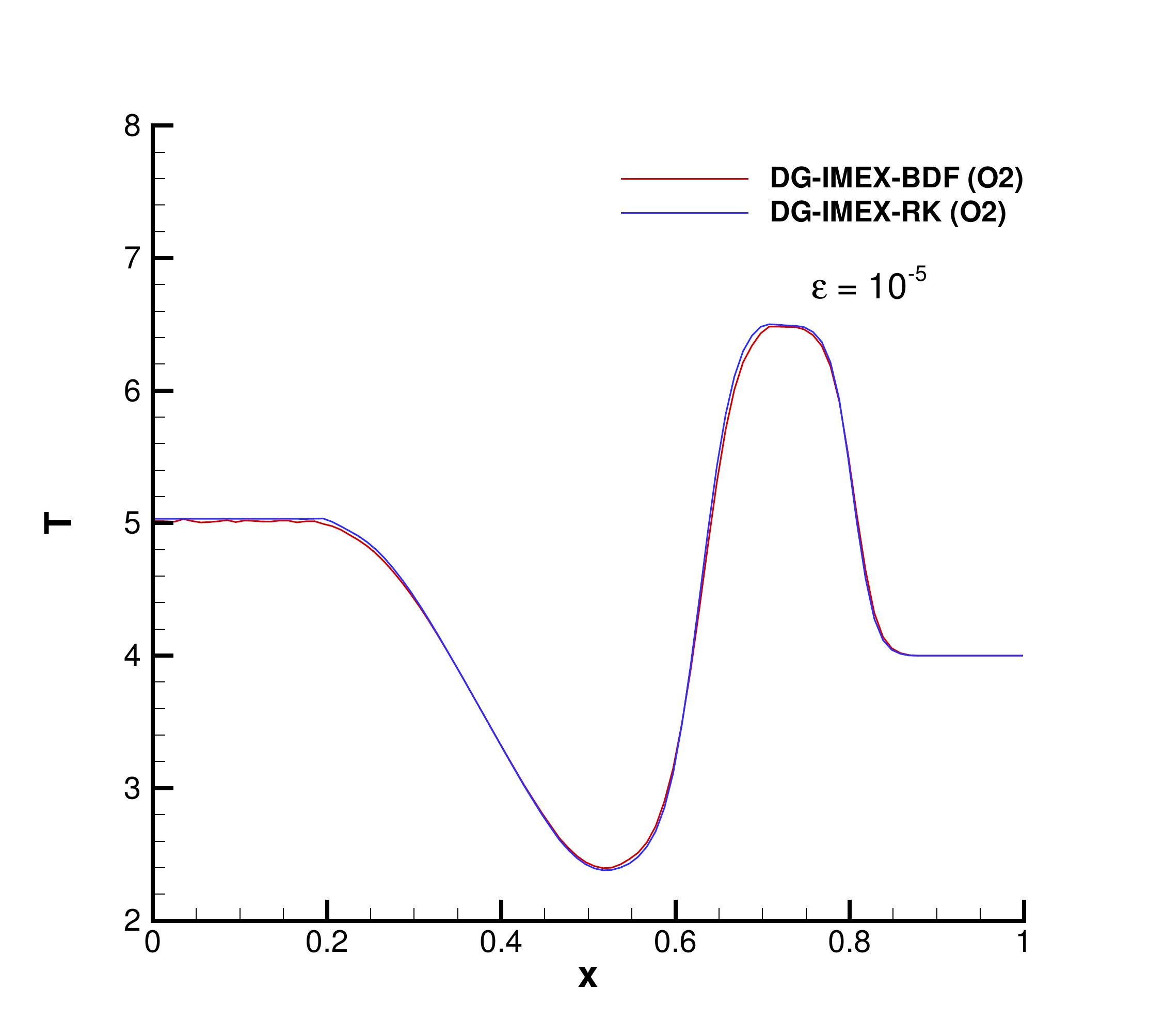}  \\ 
			\includegraphics[width=0.33\textwidth]{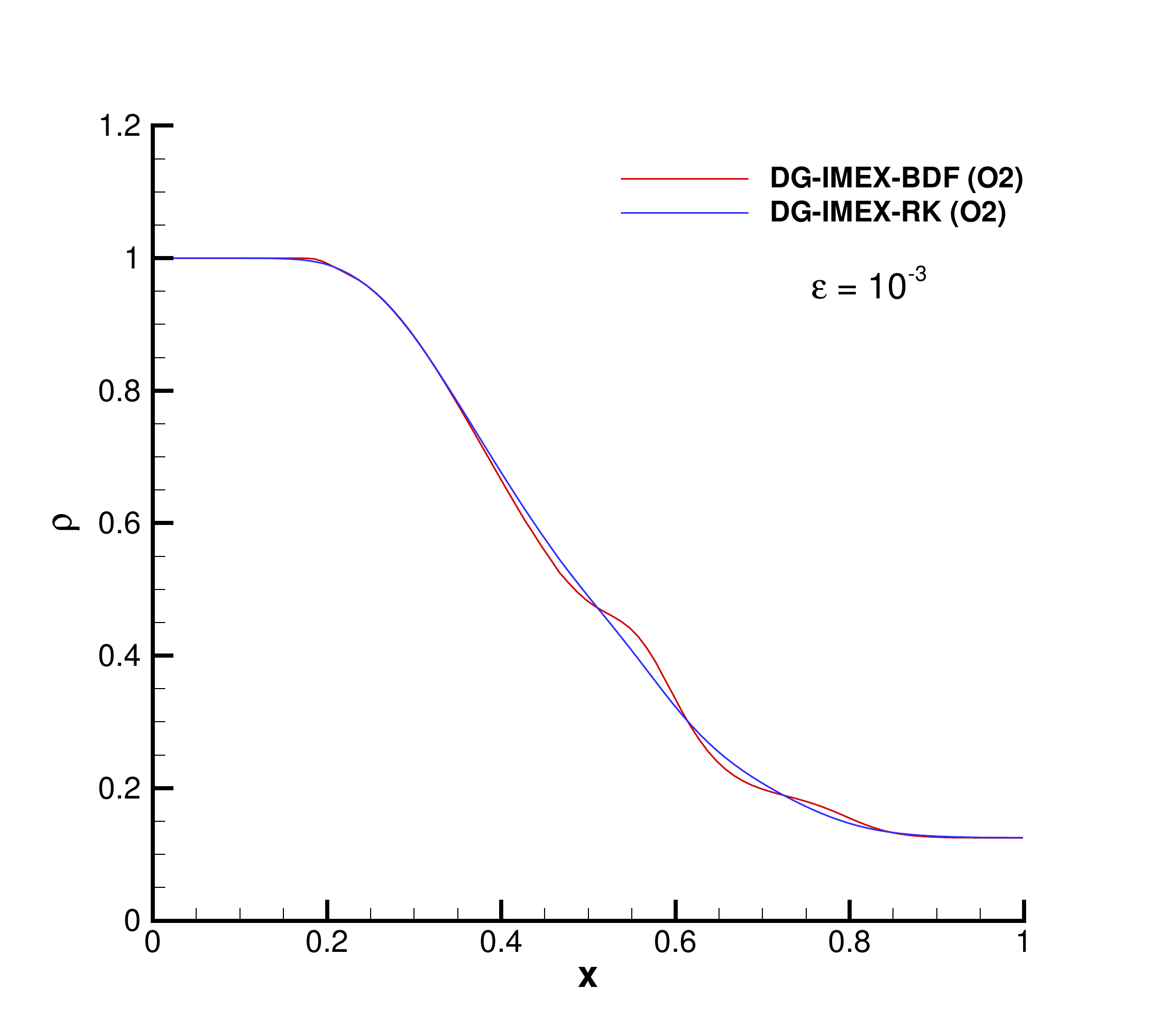}  &            
			\includegraphics[width=0.33\textwidth]{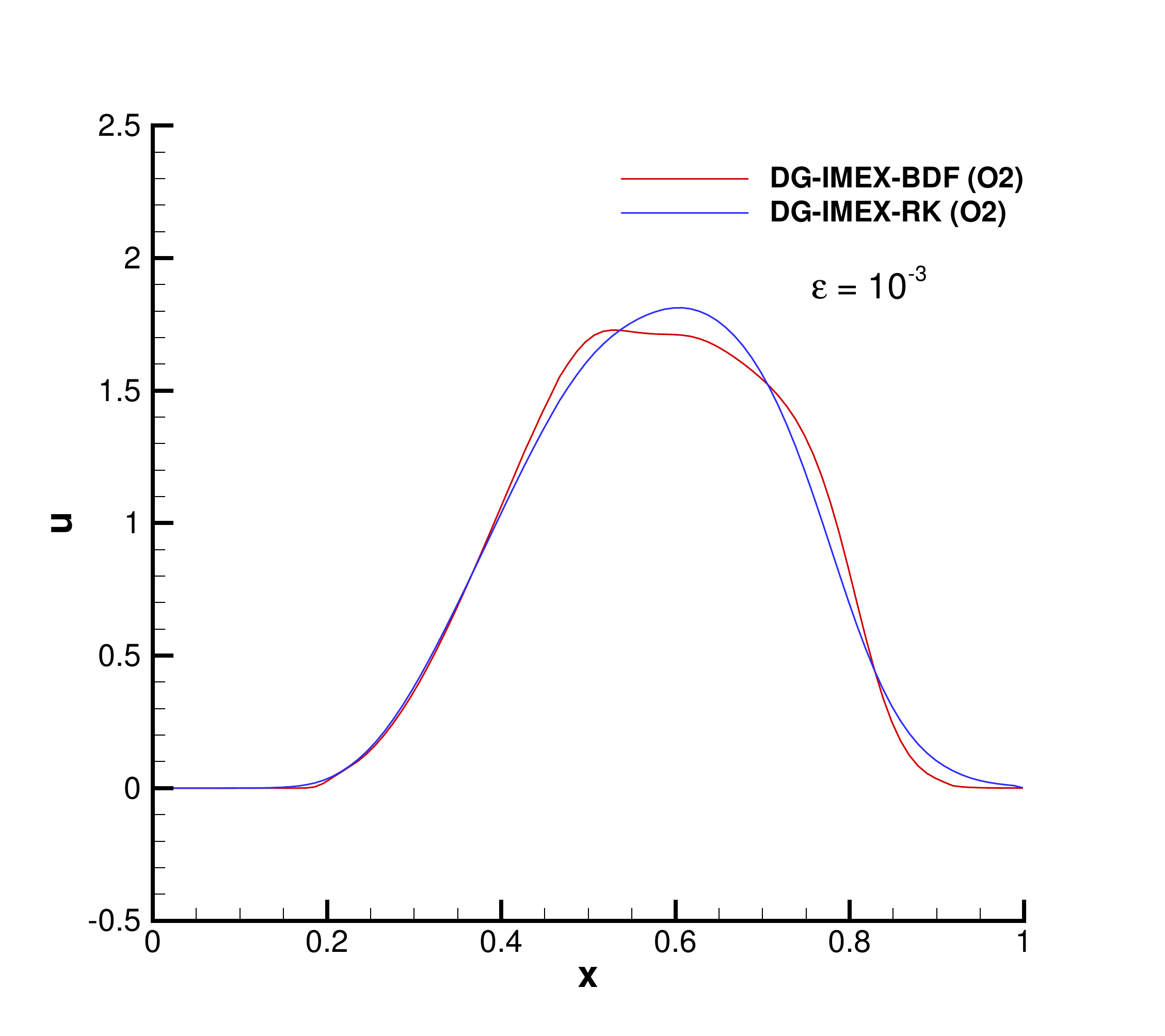}  &           
			\includegraphics[width=0.33\textwidth]{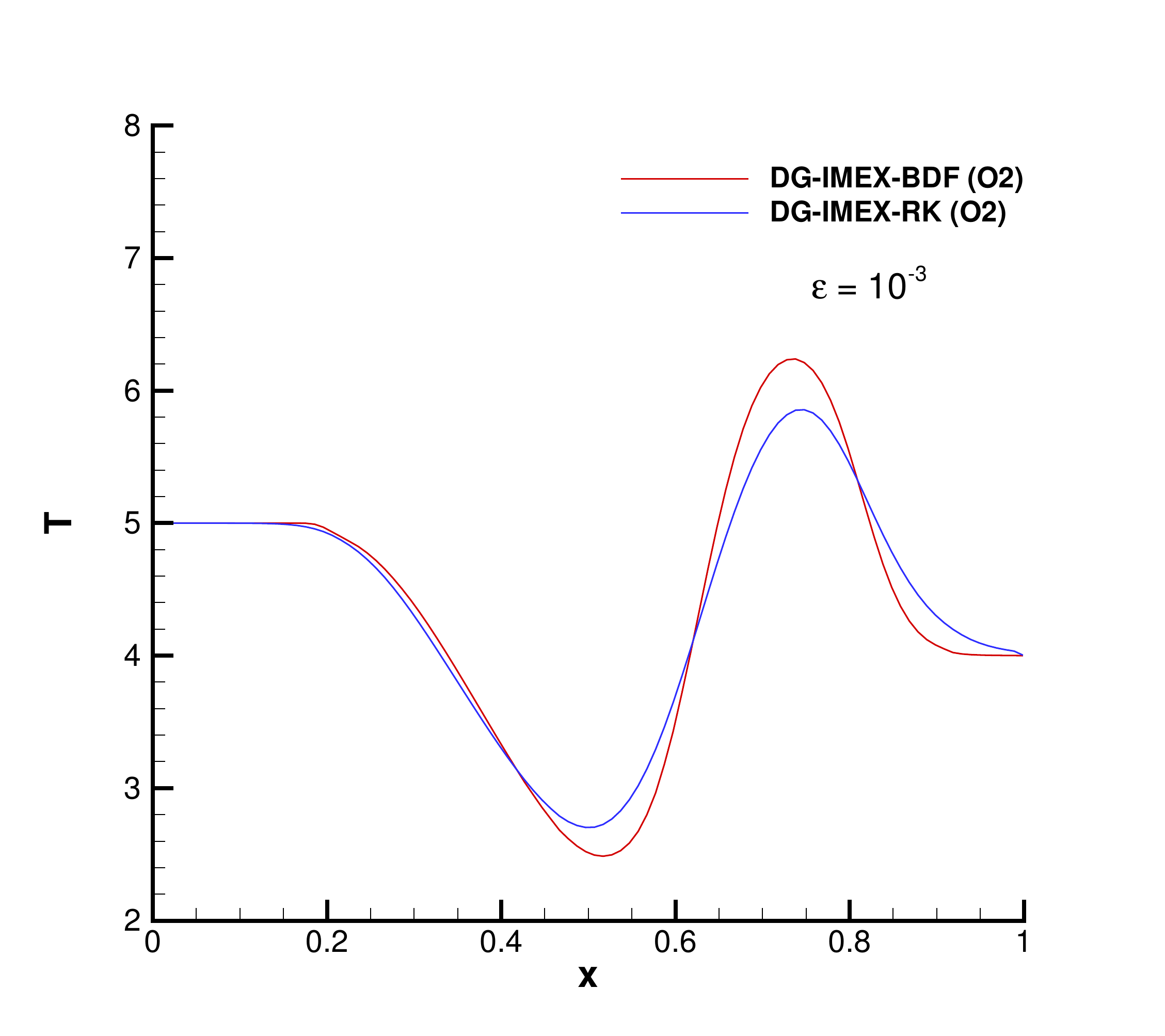}  \\ 
		\end{tabular}
		\caption{Explosion problem. Comparison between second order DG-IMEX-BDF (red solid line) and DG-IMEX-RK (blue solid line) schemes for $\varepsilon=10^{-5}$ (top row) and $\varepsilon=10^{-3}$ (bottom row). 1D cut along the $x$-axis through the third order numerical results and comparison with exact solution for density (left), velocity in the $x$-direction (middle) and temperature (right).}
		\label{fig.EP3D-comp}
	\end{center}
\end{figure}

\begin{figure}[!htbp]
	\begin{center}
		\begin{tabular}{cc} 
			\includegraphics[width=0.47\textwidth]{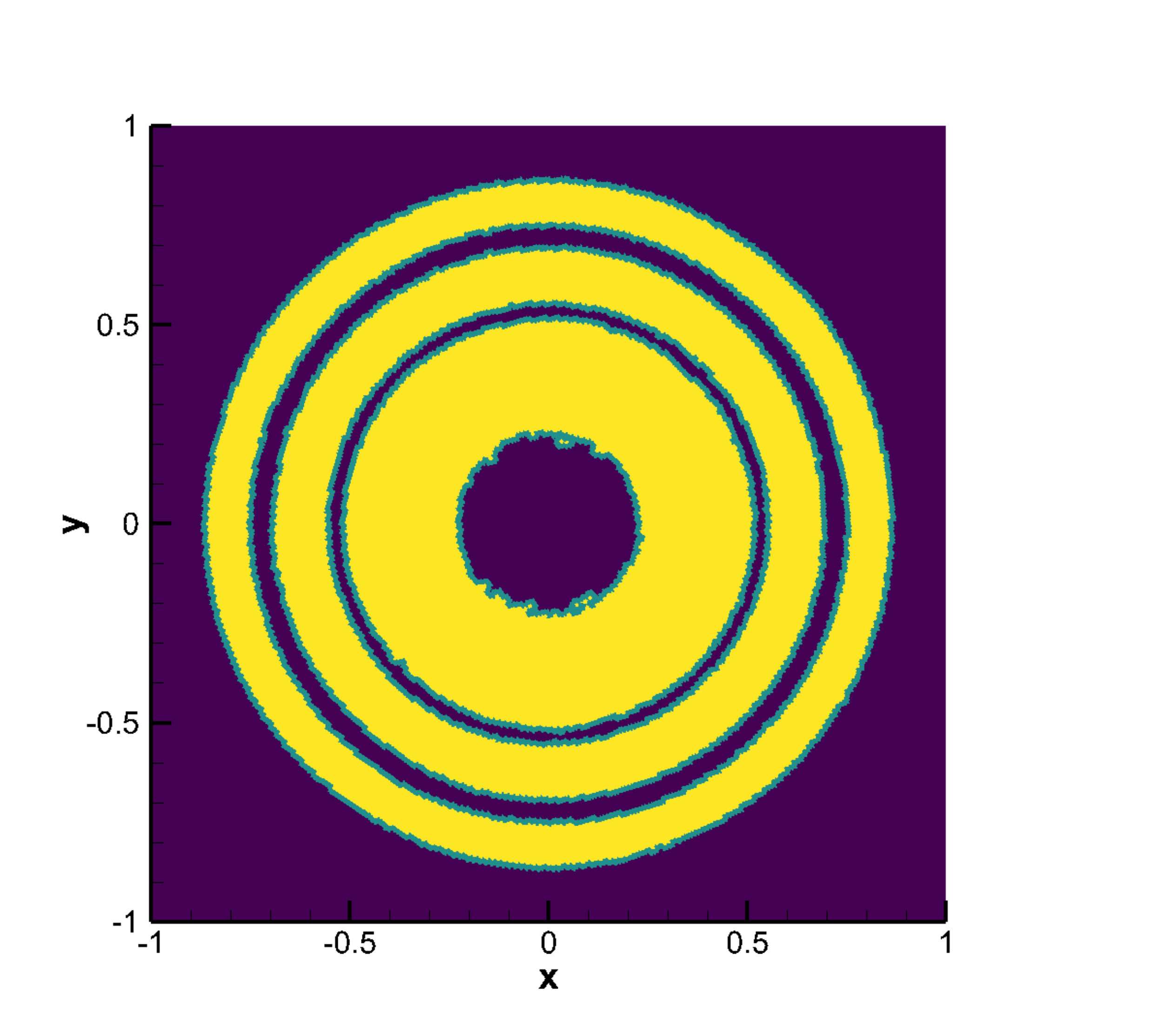}  &           
			\includegraphics[width=0.47\textwidth]{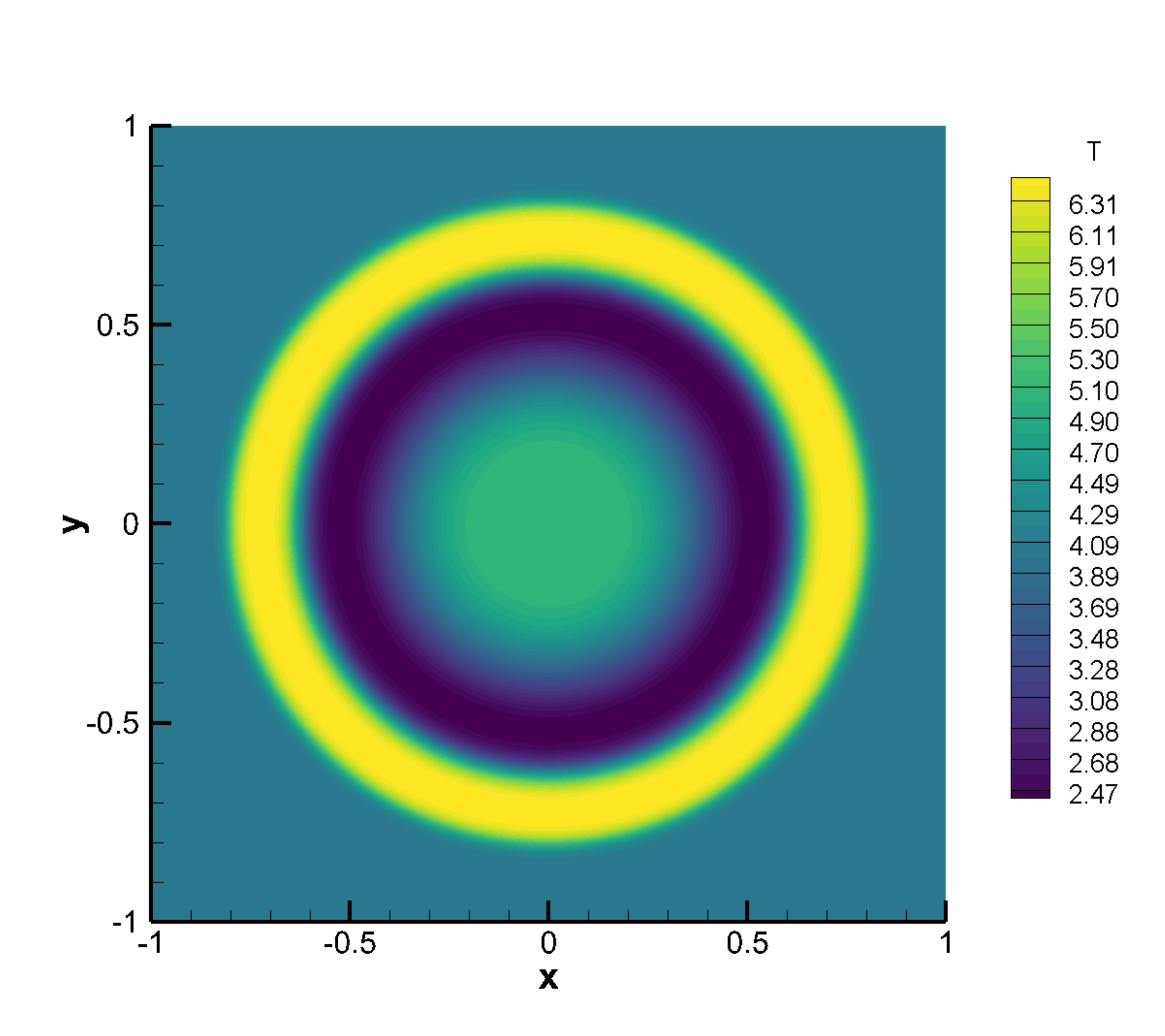}  \\           
		\end{tabular}
		\caption{Explosion problem with $\varepsilon=10^{-5}$. 
			Map of the troubled cells where the limiter is applied (left) and temperature distribution with 40 contours in the range $[2;6.5]$ (right).}
		\label{fig.EP3D-lim}
	\end{center}
\end{figure}

For comparison purposes, the same test is run again with both BDF and RK schemes in a rarefied gas regime, namely setting Knudsen number $\varepsilon=10^{-3}$. The results are depicted in Figure \ref{fig.EP3D-comp} for both values of the Knudsen number. In the fluid limit, no differences arise, whereas for $\varepsilon=10^{-3}$ the BDF schemes are less dissipative than the corresponding RK time stepping techniques as expected \cite{Ascher2}. 
Figure \ref{fig.EP3D-lim} depicts the map of the troubled cells for the DG limiter at the final time as well as the temperature contours on the two-dimensional domain, showing that the cylindrical symmetry of the reference solution is very well retrieved by the numerical results, despite the unstructured mesh composed of arbitrarily shaped polygons. Moreover, let us notice that the plateau between two discontinuities 
are not affected by the limiting procedure since those cells are not marked as troubled by the detector.
Finally, to measure the computational efficiency of the BDF linear multistep methods, the computational times for the Lax test and the explosion problem are monitored and reported in Table \ref{tab.efficiency}. The lower is the Knudsen number, the higher is the advantage of using the BDF time integrator over the RK time stepping. Indeed, the DG limiter is quite time consuming and in the limit $\varepsilon \to 0$ acts more frequently causing the increase of the computational effort for both schemes. Table \ref{tab.efficiency} contains the computational time $\tau_{dof}$ needed for updating one degree of freedom of the DG solution within one time step as well as the ratio $\mathcal{R}_{dof}$ between the values obtained for RK and BDF. The BDF schemes are up to $\approx 4$ times faster than the corresponding RK methods, while exhibiting less numerical diffusion as previously noticed in Figure \ref{fig.EP3D-comp}.

\begin{table}[!htbp]  
	\caption{Computational time for Lax test (Lax) and explosion problem (EP2D) with different Knudsen numbers run using DG-IMEX-RK and DG-IMEX-BDF time integrators. The number of degrees of freedom involved in the spatial discretization are addressed with ${dof}$, the number of time steps needed to reach the final time is given by $N_{\Delta t}$, while absolute time of each simulation measured in seconds [s] is provided for both schemes, that is $T_{RK}$ and $T_{BDF}$. The time needed for the update of one degree of freedom within one time step is indicated with $\tau_{RK}$ and $\tau_{BDF}$. Finally, the efficiency ratio is evaluated as $\mathcal{R}_{dof}=\tau_{RK}/\tau_{BDF}$.}  
	\begin{center} 
		\begin{small}
			\renewcommand{\arraystretch}{1.0}
			\begin{tabular}{l|rrr|cc|cc|c}
				\multicolumn{4}{c|}{} & \multicolumn{2}{c|}{DG-IMEX-RK} & \multicolumn{2}{c|}{DG-IMEX-BDF} & \\
				Test & Knudsen & ${dof}$ & $N_{\Delta t}$ & $T_{RK}$ [s] & $\tau_{RK}$ [s] & $T_{BDF}$ [s] & $\tau_{BDF}$ [s] & $\mathcal{R}_{dof}$\\ 
				\hline
				Lax  & $10^{-3}$ &  36468 & 5120 & 1.5127E+07 & 8.1015E-02 & 5.1338E+06 & 2.7495E-02 & 2.95 \\ 
				     & $10^{-5}$ &  36468 & 5120 & 1.9135E+07 & 1.0248E-01 & 5.1547E+06 & 2.7607E-02 & 3.71 \\
				EP2D & $10^{-3}$ & 186921 & 3236 & 4.0158E+07 & 6.6390E-02 & 2.1381E+07 & 3.5347E-02 & 1.88 \\ 
				     & $10^{-5}$ & 186921 & 3236 & 4.4380E+07 & 7.3370E-02 & 2.1433E+07 & 3.5434E-02 & 2.07 \\
			\end{tabular}
		\end{small}
	\end{center}
	\label{tab.efficiency}
\end{table}

%
\subsection{Fluid flow around NACA 0012 airfoil}
\begin{figure}[!htbp]
	\begin{center}
		\begin{tabular}{ccc} 
			\includegraphics[width=0.33\textwidth]{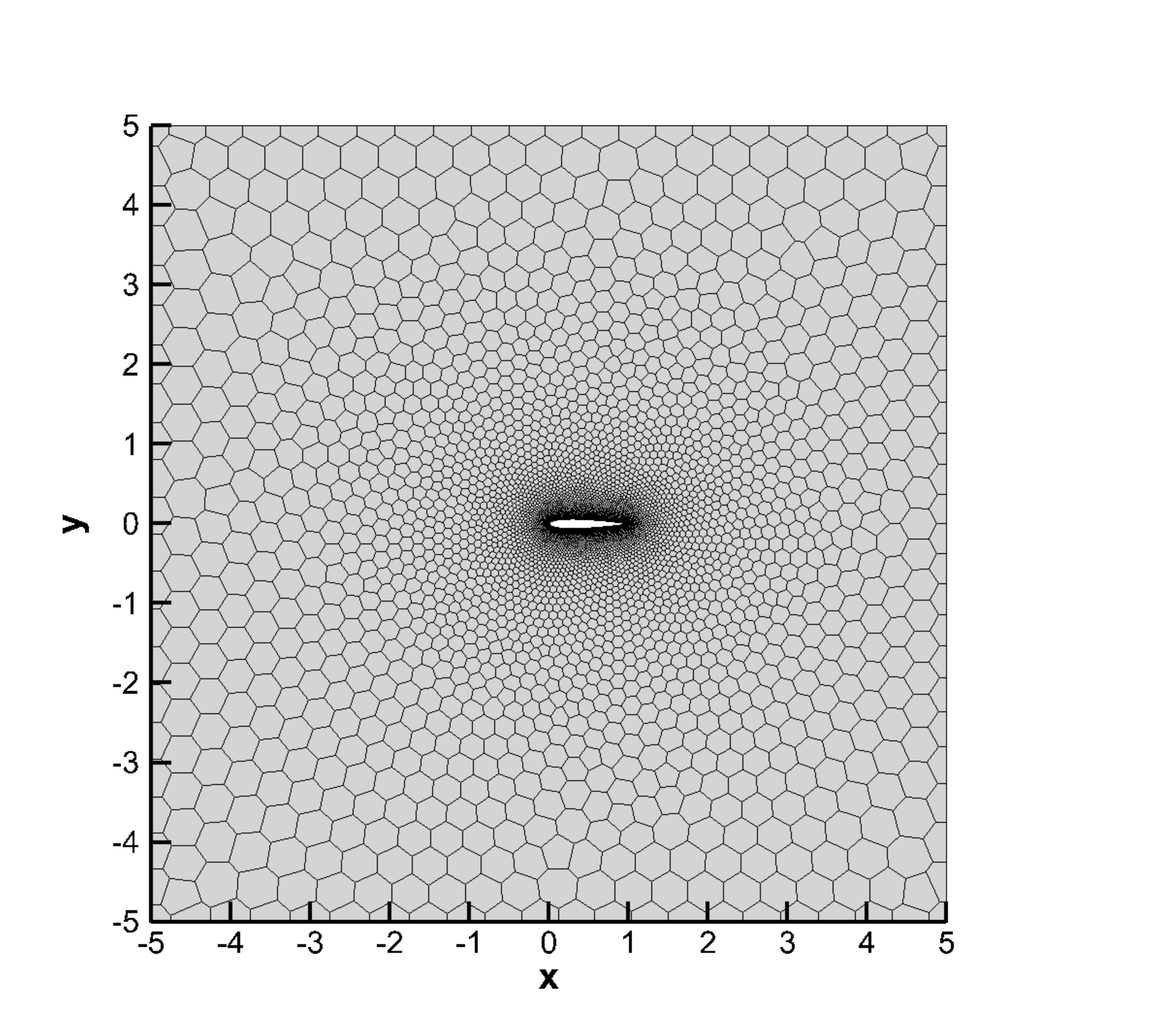}  &          
			\includegraphics[width=0.33\textwidth]{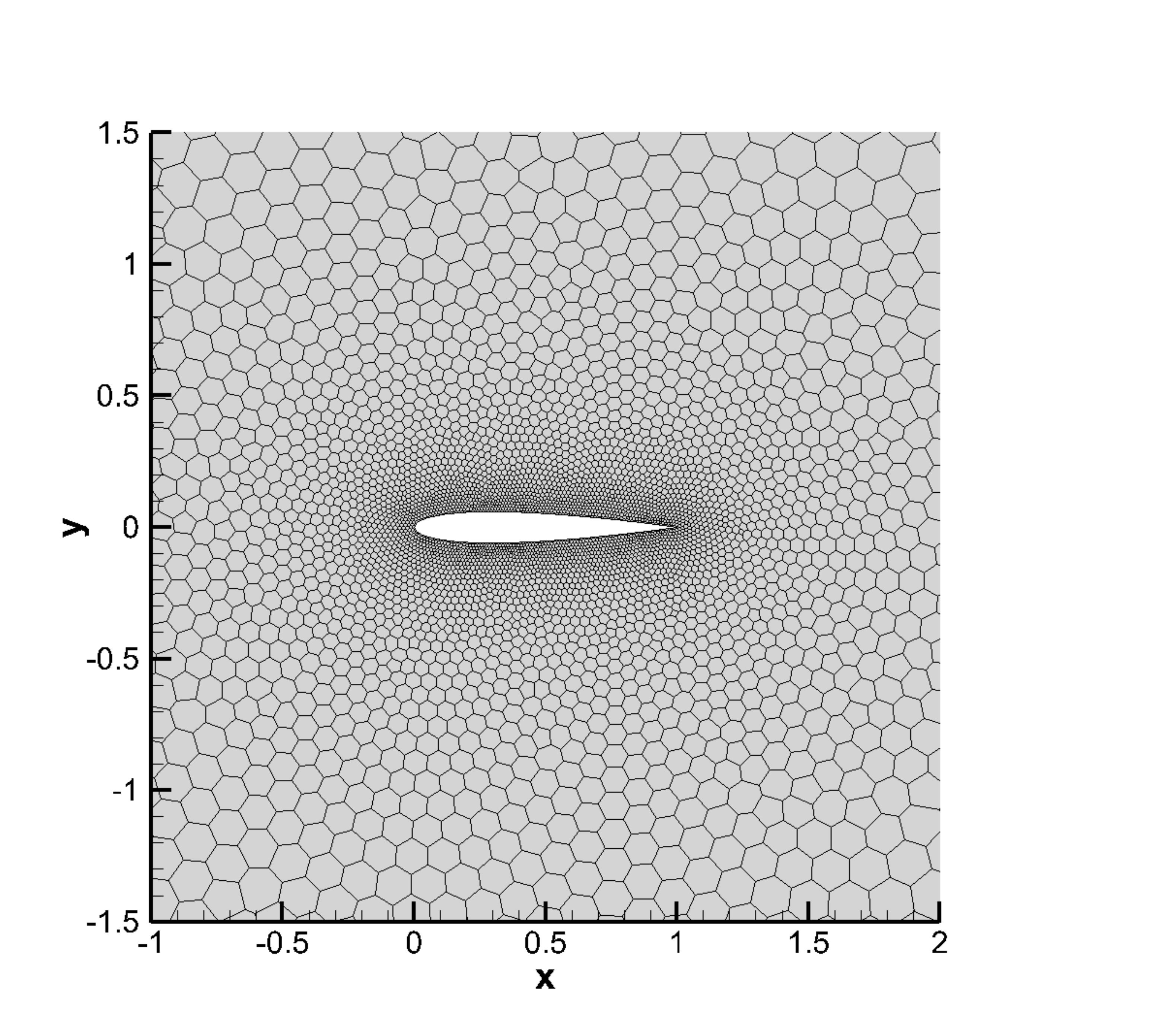}  &
			\includegraphics[width=0.33\textwidth]{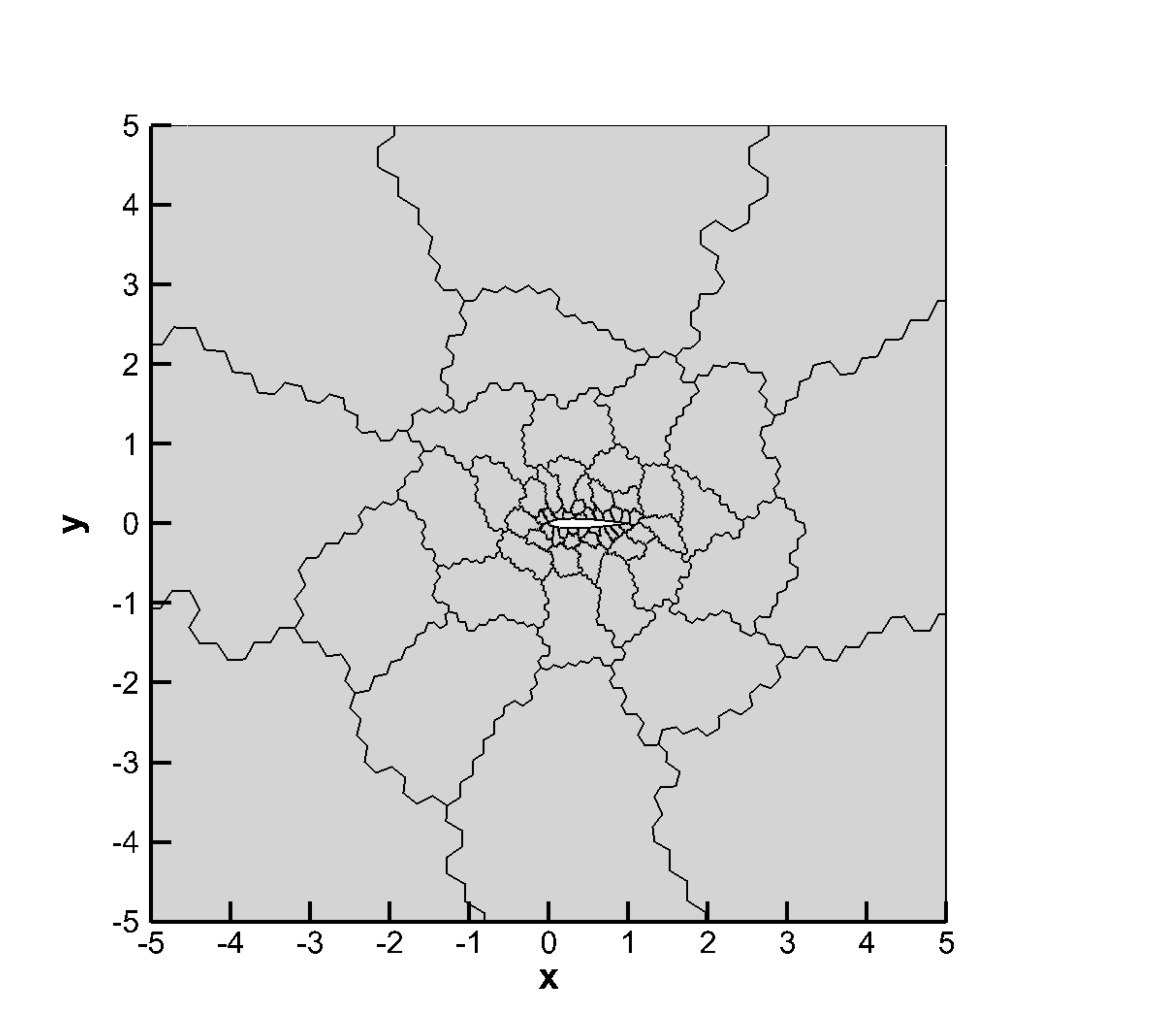}  \\           
		\end{tabular}
		\caption{Flow around NACA 0012 airfoil. Computational mesh with $N_P=10336$ polygonal cells (left), zoom around the airfoil profile discretized with 200 points (middle) and MPI domains which splits the computational domain into a total number of 64 sub-domains.}
		\label{fig.NACA-mesh}
	\end{center}
\end{figure}

As last test case we propose to simulate a more realistic scenario involving the flow around a NACA 0012 airfoil profile for different rarefied regimes. The original geometry of the NACA 0012 is slightly modified, so that the airfoil closes at chord $c=1$ with a sharp trailing edge. This altered two-dimensional geometrical configuration can be described by considering
\begin{equation}
	y = \pm 0.594689181 \, \left( 0.298222773\sqrt{x} - 0.127125232 x - 0.357907906 x^2 + 0.291984971 x^3 - 0.105174606 x^4 \right),
\end{equation}
which corresponds to a perfect scaled copy of the NACA 0012 profile with $x \in [0;1.008930411365]$. The airfoil is embedded in a square computational domain $\Omega=[-5;5]^2$ that is paved with a total number of $N_P=10336$ polygonal control volumes. Slip-wall conditions are prescribed on the airfoil boundary, that is discretized with 100 equidistant points in the $x-$direction both on the upper and the lower border. The characteristic mesh size is then progressively increased linearly with the distance from the airfoil, thus ranging from $h=1/100$ to $h=1/10$. On the outer boundary of the domain we set transmissive conditions everywhere, apart from the left side where a space-time-dependent inflow boundary condition will be imposed. The entire computational mesh and a zoom around the NACA profile is depicted in Figure \ref{fig.NACA-mesh}, together with the MPI partition of the domain within 64 CPUs. The velocity domain is discretized with $16^2=256$ Cartesian cells ranging in the interval $\mathcal{V}=[-5;5]\times [-5;5]$ giving rise to a problem involving approximately $2\cdot 10^6$ mesh points.

\begin{figure}[!htbp]
	\begin{center}
		\begin{tabular}{cc} 
			\includegraphics[width=0.47\textwidth]{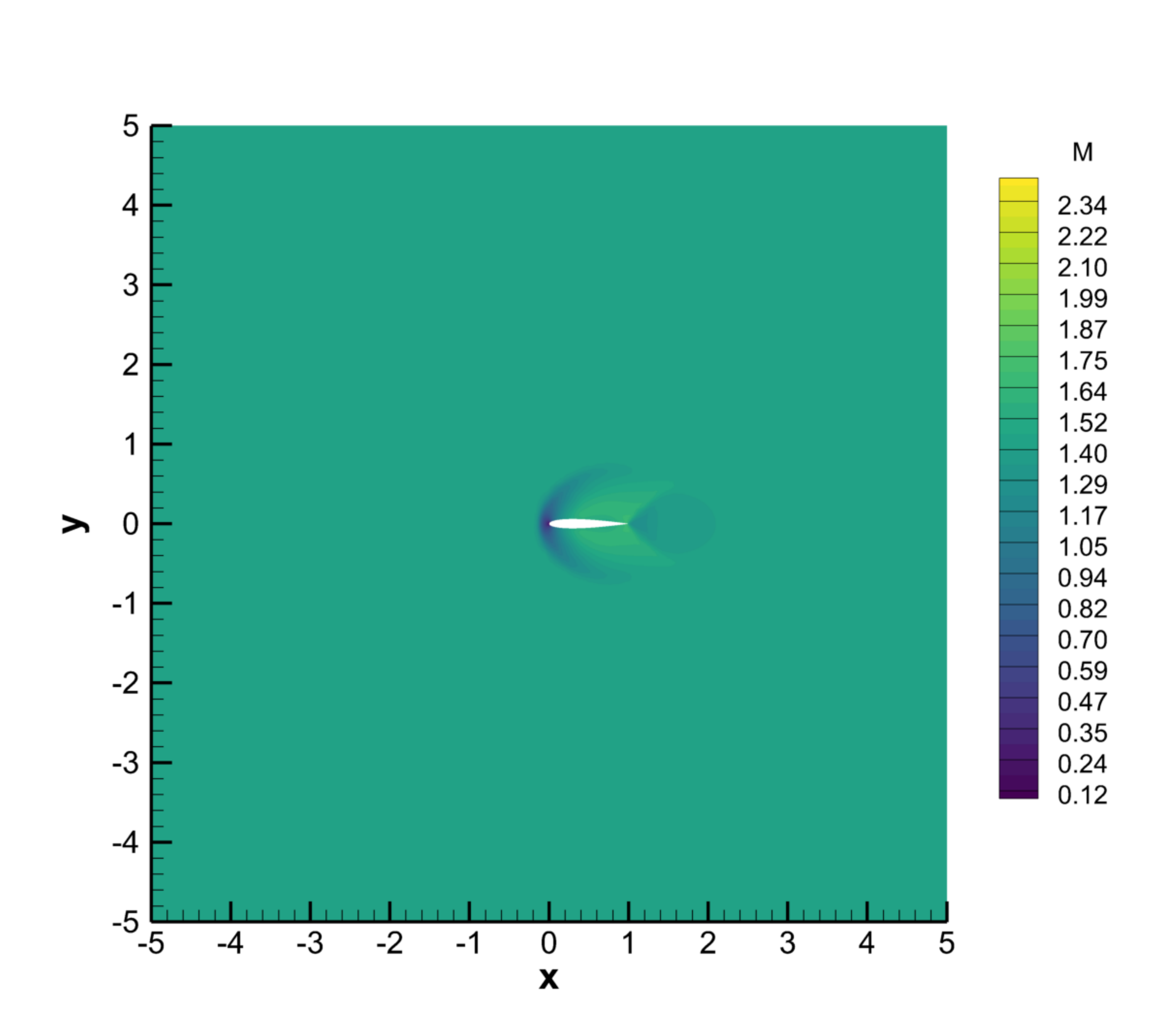}  &           
			\includegraphics[width=0.47\textwidth]{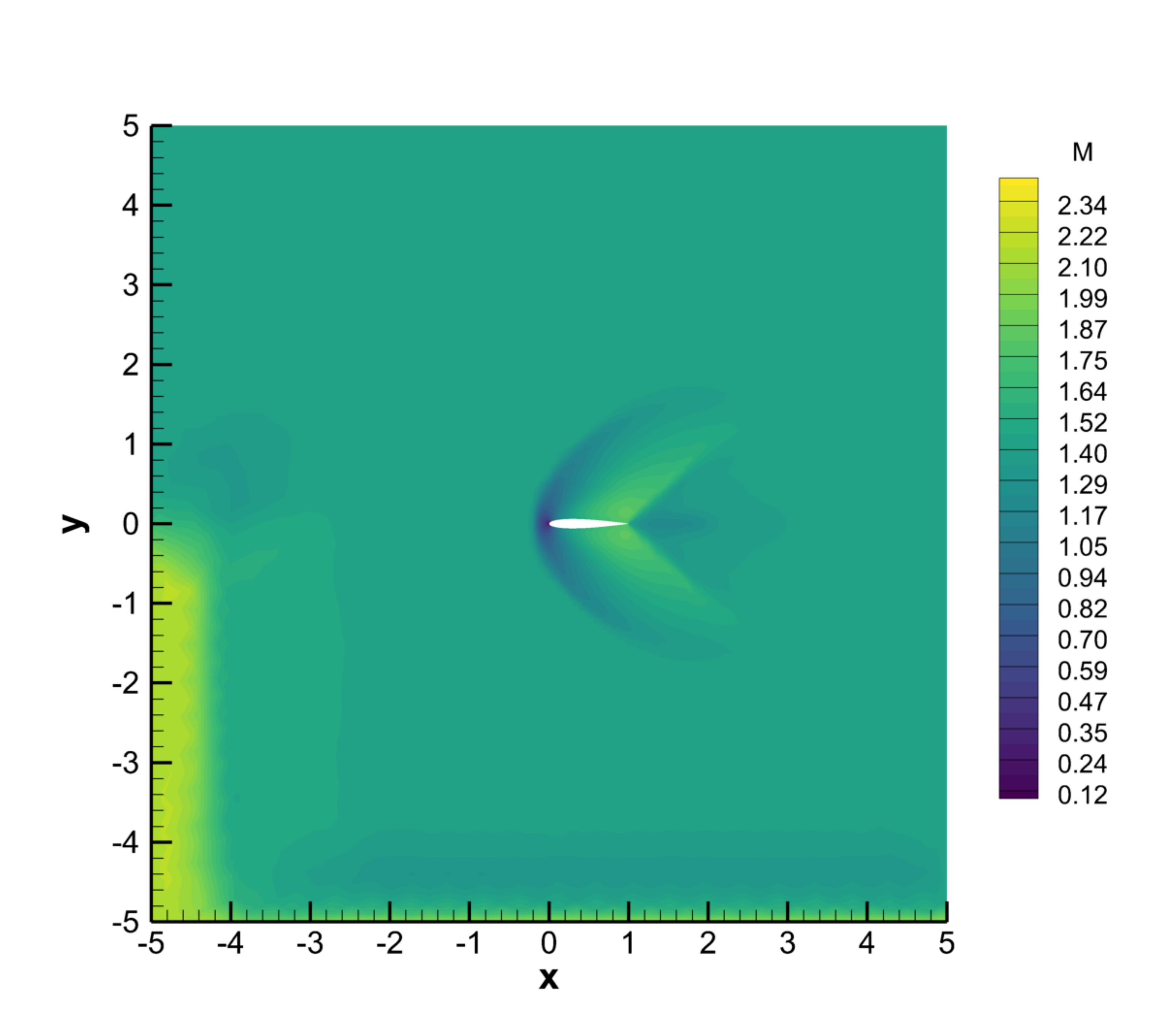}  \\          
			\includegraphics[width=0.47\textwidth]{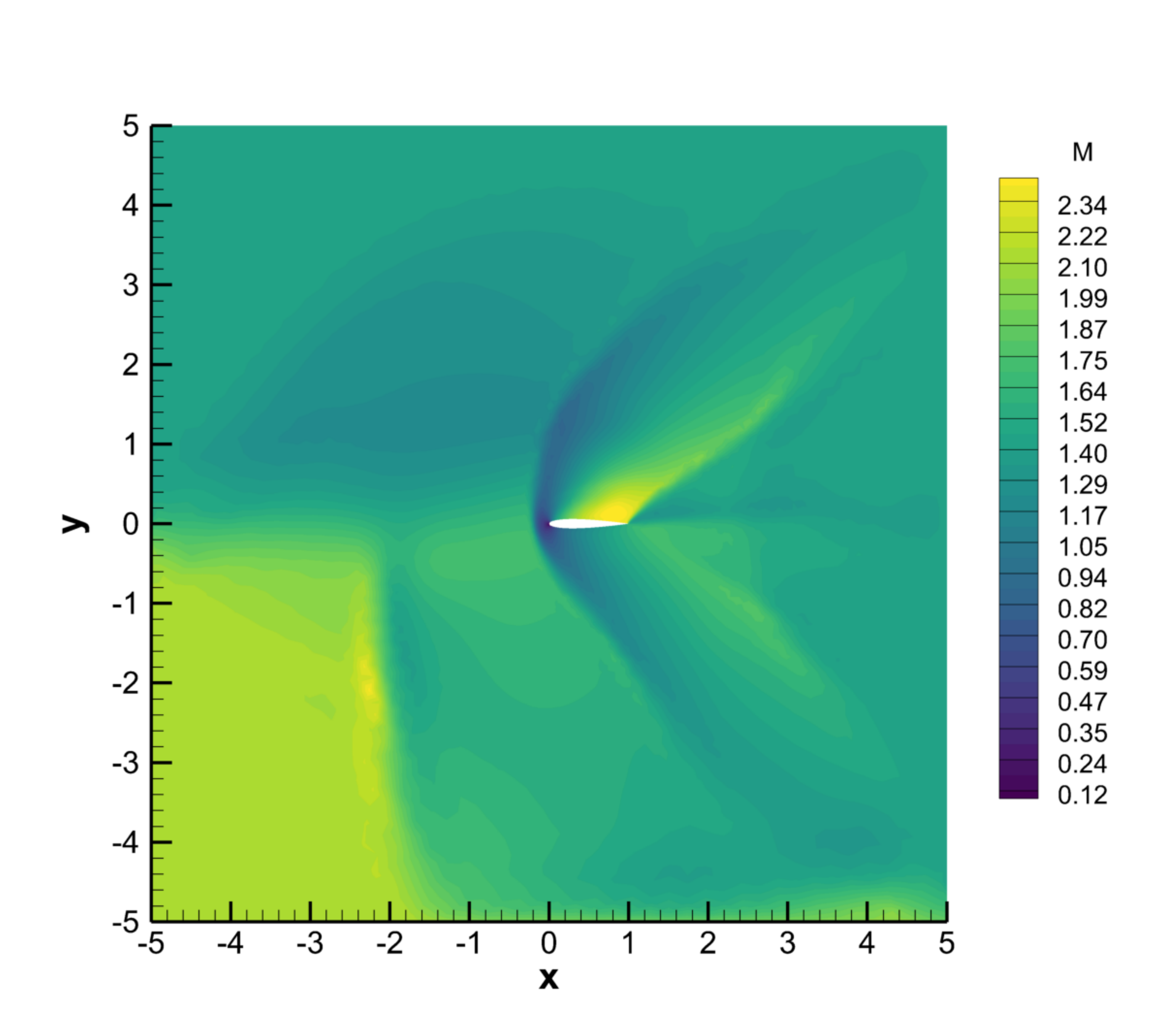}  &           
			\includegraphics[width=0.47\textwidth]{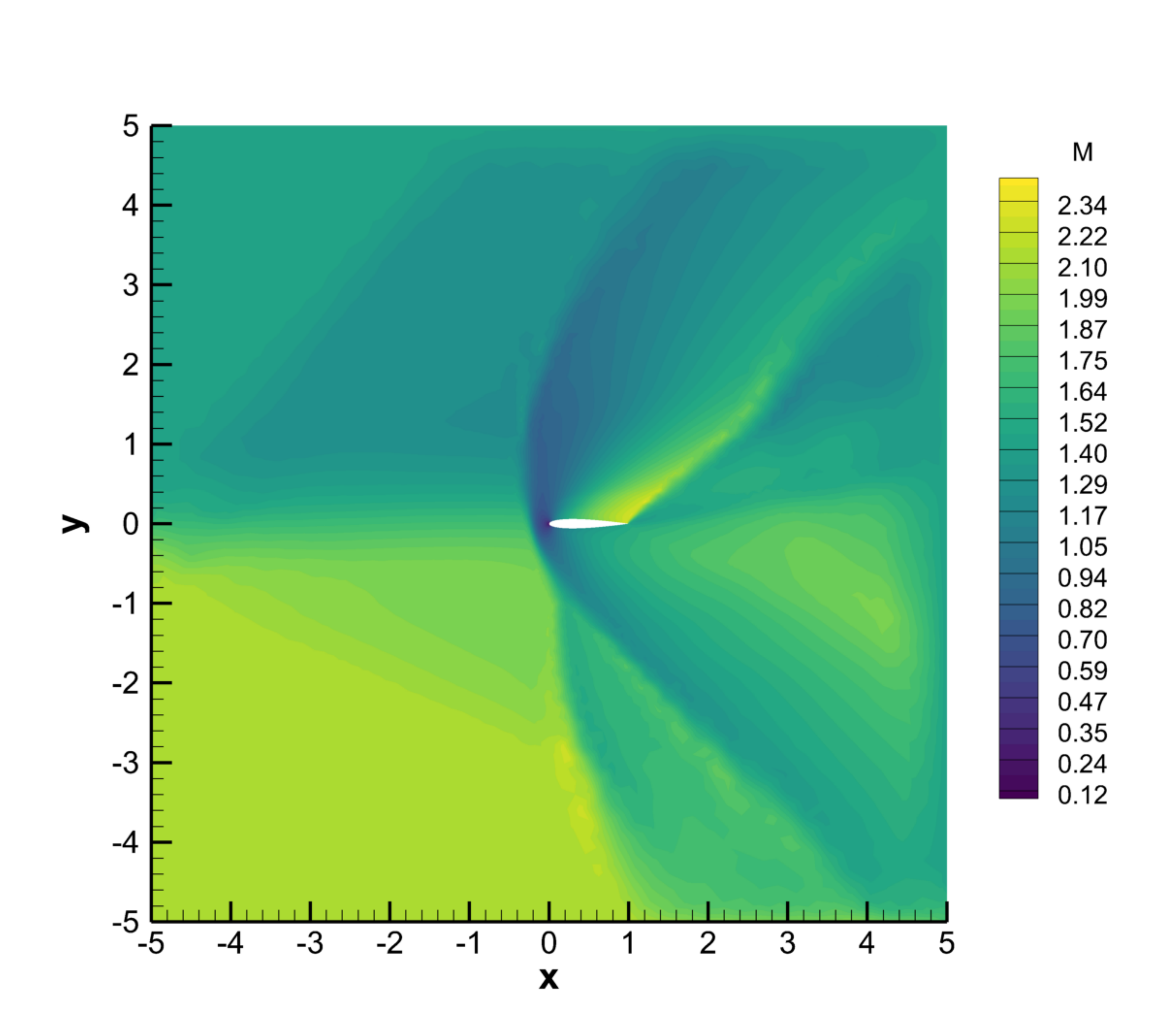}  \\ 
		\end{tabular}
		\caption{Flow around NACA 0012 airfoil with $\varepsilon=5\cdot 10^{-4}$ and time-dependent inflow angle of attack. From top left to bottom right: 40 Mach number contours in the interval $[0.1;2.4]$ at time $t=0.4$, $t=1$, $t=3$ and $t=5$. }
		\label{fig.NACA_Kn5-4_M}
	\end{center}
\end{figure}

\begin{figure}[!htbp]
	\begin{center}
		\begin{tabular}{cc} 
			\includegraphics[width=0.47\textwidth]{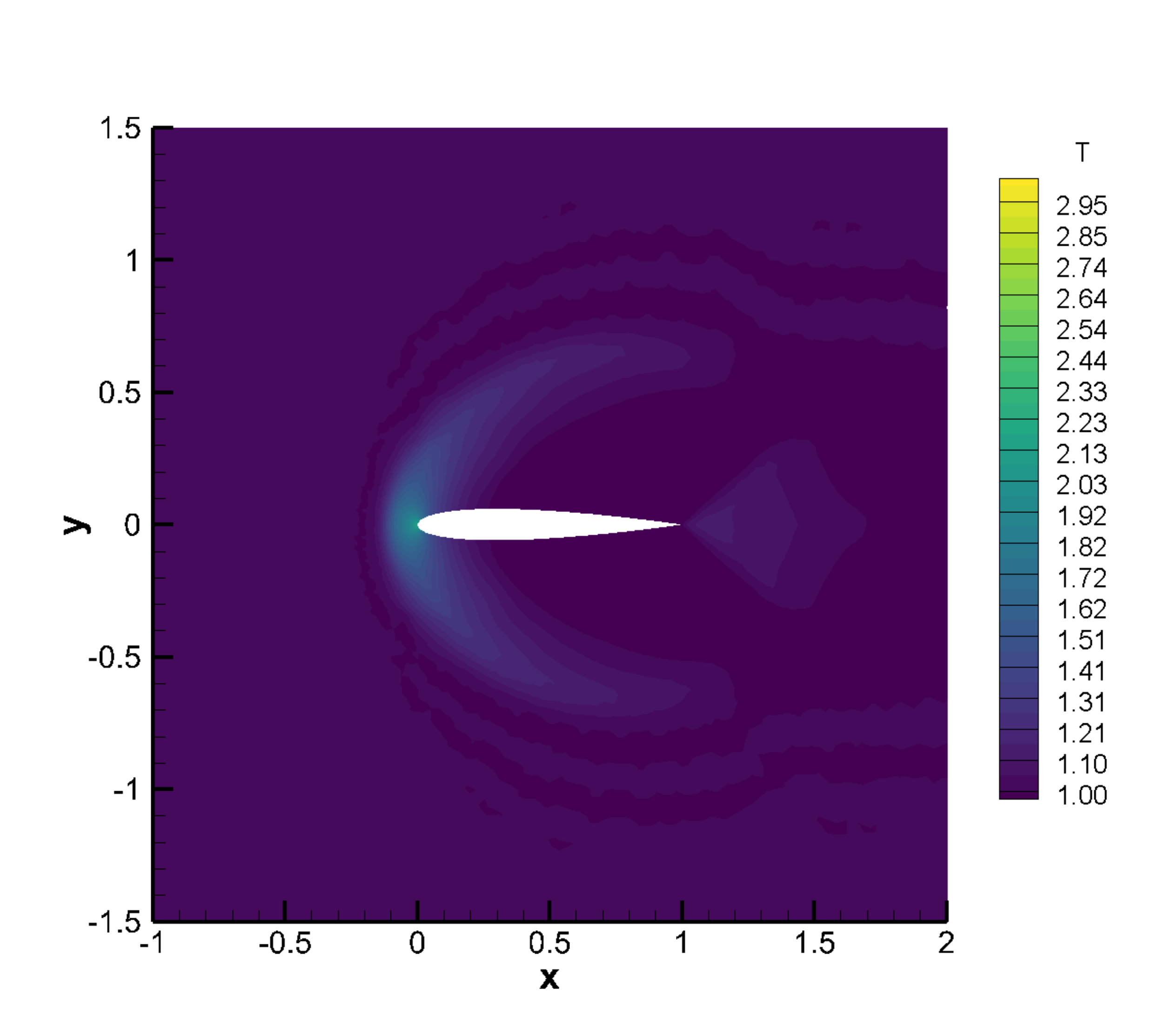}  &           
			\includegraphics[width=0.47\textwidth]{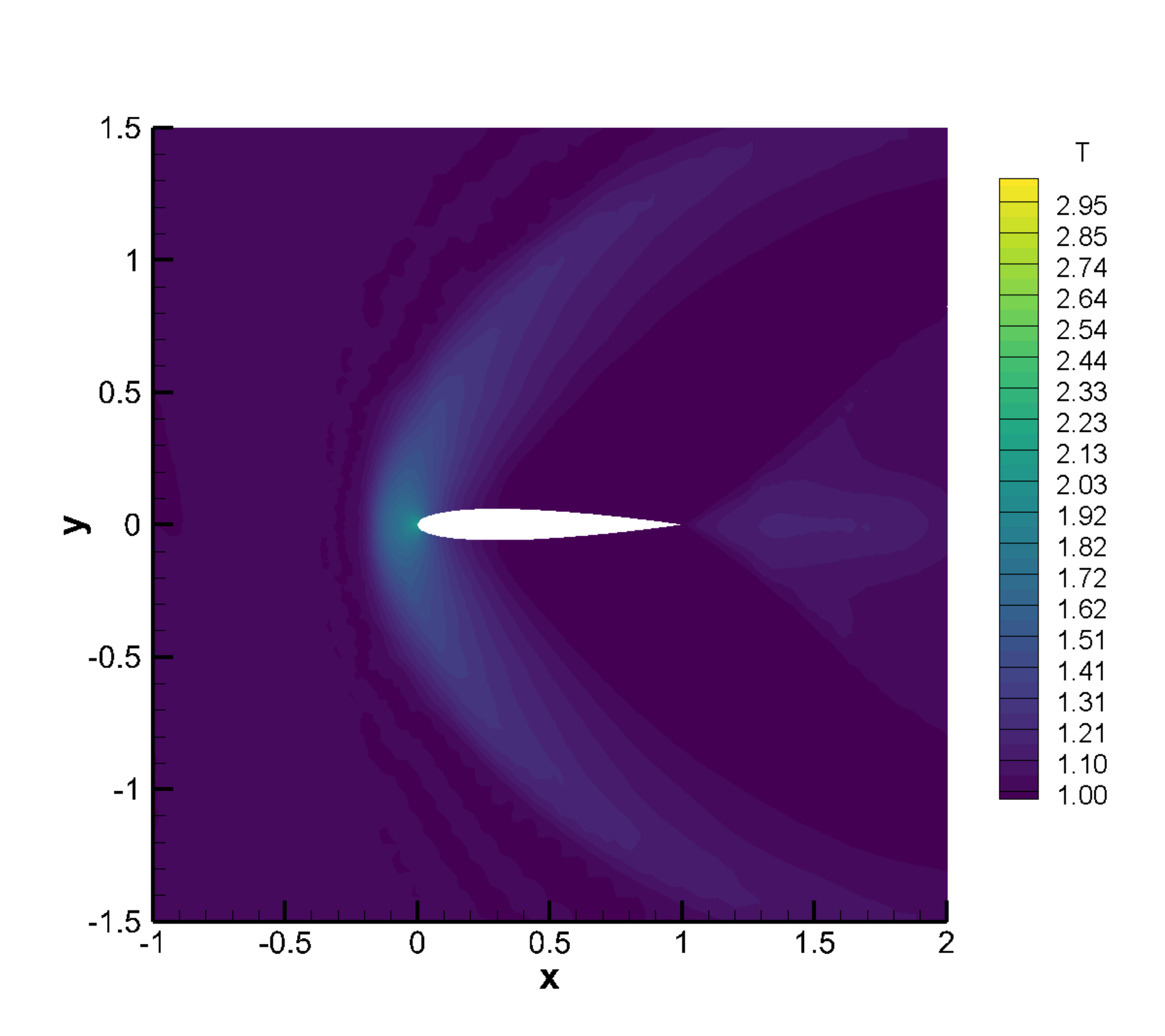}  \\          
			\includegraphics[width=0.47\textwidth]{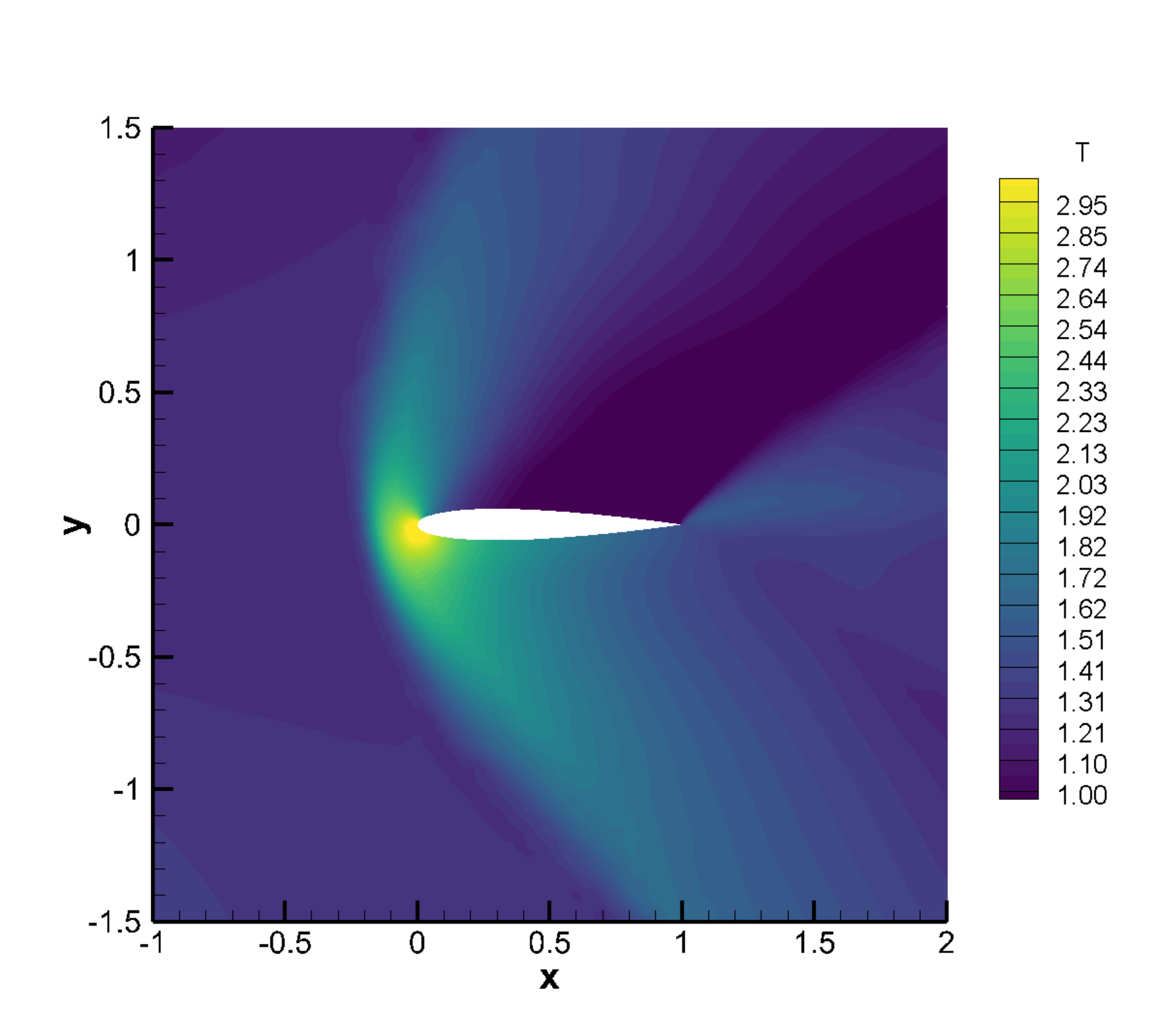}  &           
			\includegraphics[width=0.47\textwidth]{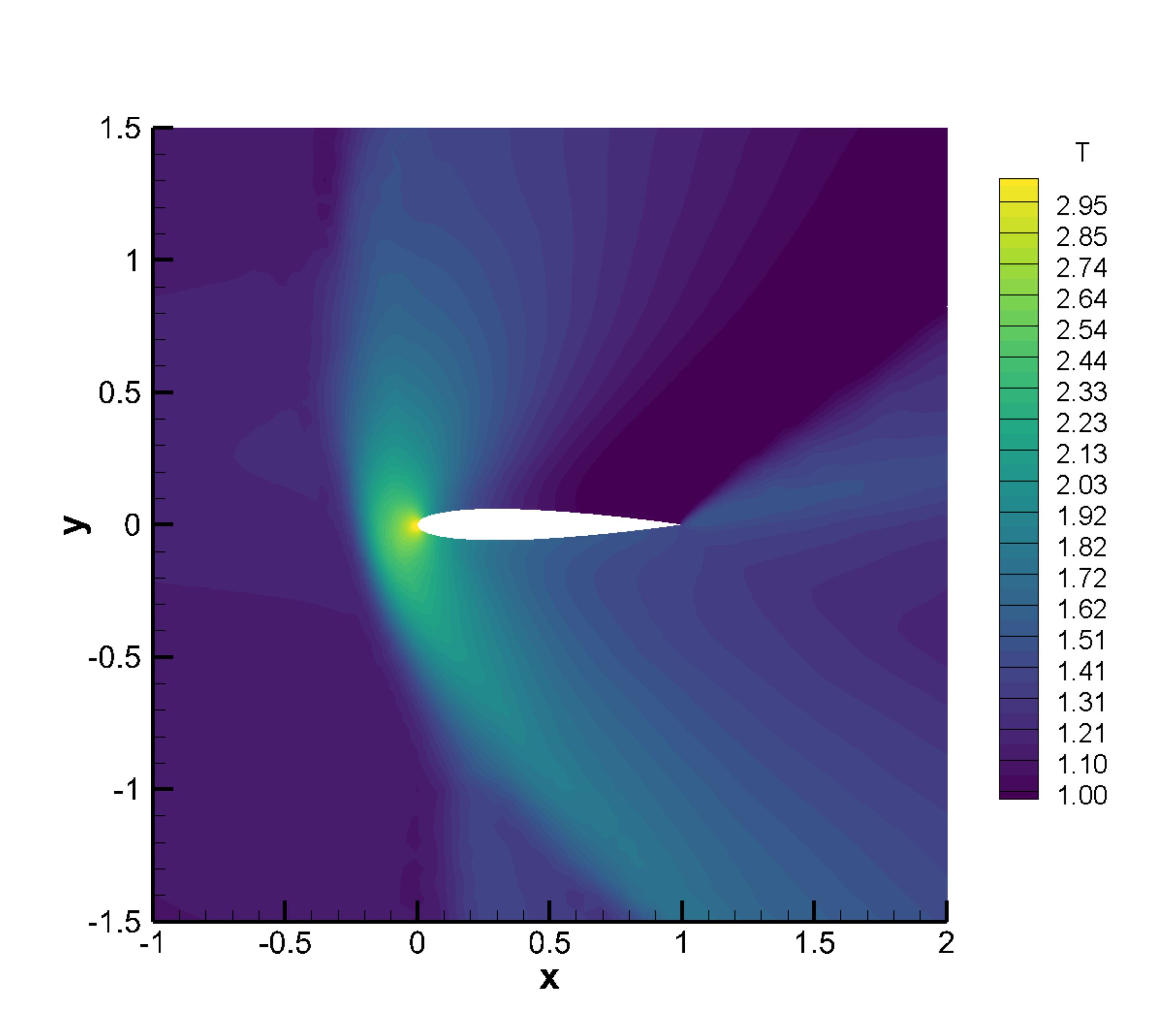}  \\ 
		\end{tabular}
		\caption{Flow around NACA 0012 airfoil with $\varepsilon=5\cdot 10^{-4}$ and time-dependent inflow angle of attack. From top left to bottom right: 40 temperature contours in the interval $[1;3]$ at time $t=0.4$, $t=1$, $t=3$ and $t=5$.}
		\label{fig.NACA_Kn5-4_T}
	\end{center}
\end{figure}

Similarly to a test case proposed in \cite{FKS_DD}, here we want to mimic a change of flying conditions for the object, thus the inflow conditions on the left side of the domain are modified at a certain time of the simulation. This change affects both the inlet velocity and the angle of attack of the incoming flow. The modification of the boundary condition is sharply introduced at a prescribed time, hence no smooth transition occurs. In this way, the scheme is checked for the capability of following these rapid changes of flying conditions which give rise to complicated transitional wave patterns.
With this numerical setting, we consider a supersonic regime with an inflow Mach number of $M=2$ and an angle of attack of $\alpha=0^{\circ}$. The free-stream conditions are as follows:
\begin{equation*}
	(u,v)_\infty=(2,0), \qquad p_\infty=1, \qquad \rho_\infty=1,
\end{equation*}
thus the fluid pressure is set to $p=p_\infty/\gamma$. At time $t=0.4$, the inflow boundary conditions ($x=-5$) undergo an instantaneous change only for the portion of the left boundary defined by $y\leq 0$. Specifically, the inflow velocity is incremented to $(u,v)^\prime=(3,0)$ and the angle of attack is modified so that $\alpha^\prime=4^{\circ}$, hence simulating an increasing slope of the airfoil with respect to the free-stream flow. We consider two different Knudsen numbers, namely $\varepsilon=5 \cdot 10^{-4}$ and $\varepsilon=5 \cdot 10^{-3}$, to highlight the differences which arise from the rarefied level of the gas. The second order DG-IMEX-BDF schemes are used to carry out the simulations up to the final time $t_f=5$.
Figure \ref{fig.NACA_Kn5-4_M} shows the results with $\varepsilon=5 \cdot 10^{-4}$ at output times $t=0.4$ (when the change in the inflow boundary conditions takes place), $t=1$, $t=2$ and $t=5$. The Mach number distribution is depicted, from which it is clearly visible the shock waves that arise in front of the airfoil and departing from the tail. Notice that another shock wave is generated by the change of the inflow velocity, which propagates towards the airfoil until intersection with front shock wave occurs, as visible in the panel at time $t=5$. 
Figure \ref{fig.NACA_Kn5-4_T} plots the temperature distribution at the same output times used for the Mach number with a zoom on the region around the airfoil, where one can notice that high temperature values are located across the shock wave patterns, as expected.

\begin{figure}[!htbp]
	\begin{center}
		\begin{tabular}{ccc} 
			\includegraphics[width=0.33\textwidth]{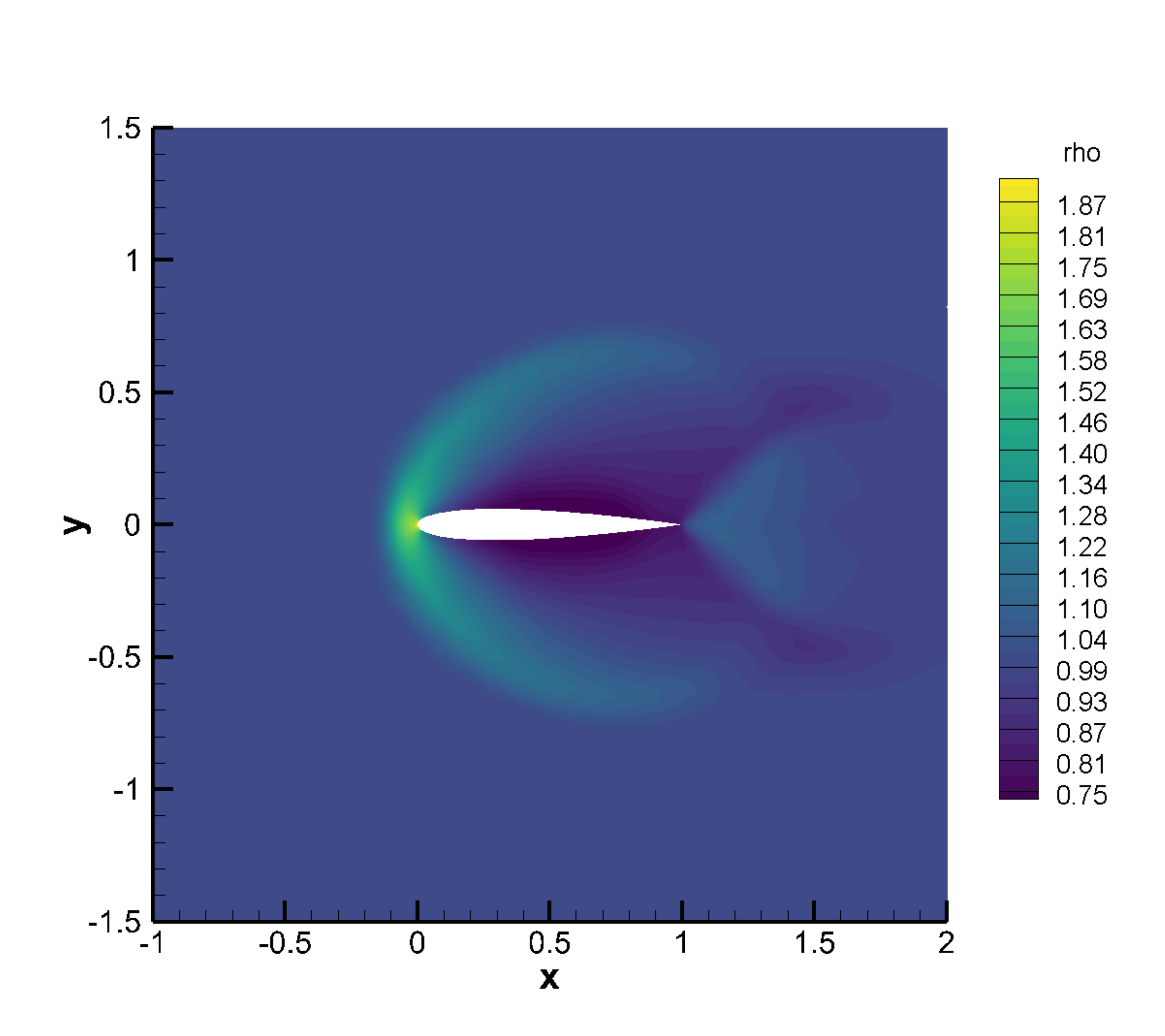}  &           
			\includegraphics[width=0.33\textwidth]{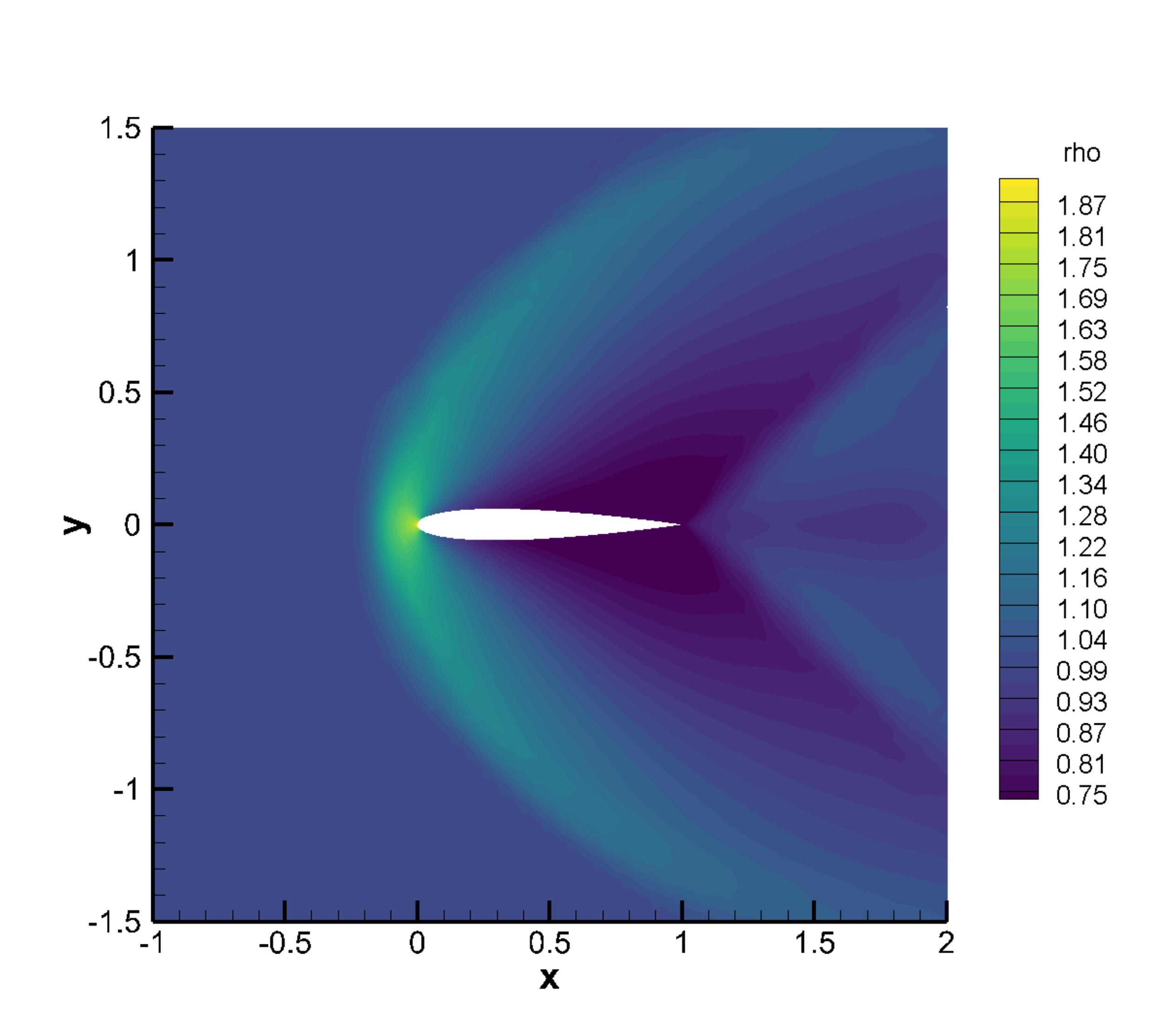}  & 
			\includegraphics[width=0.33\textwidth]{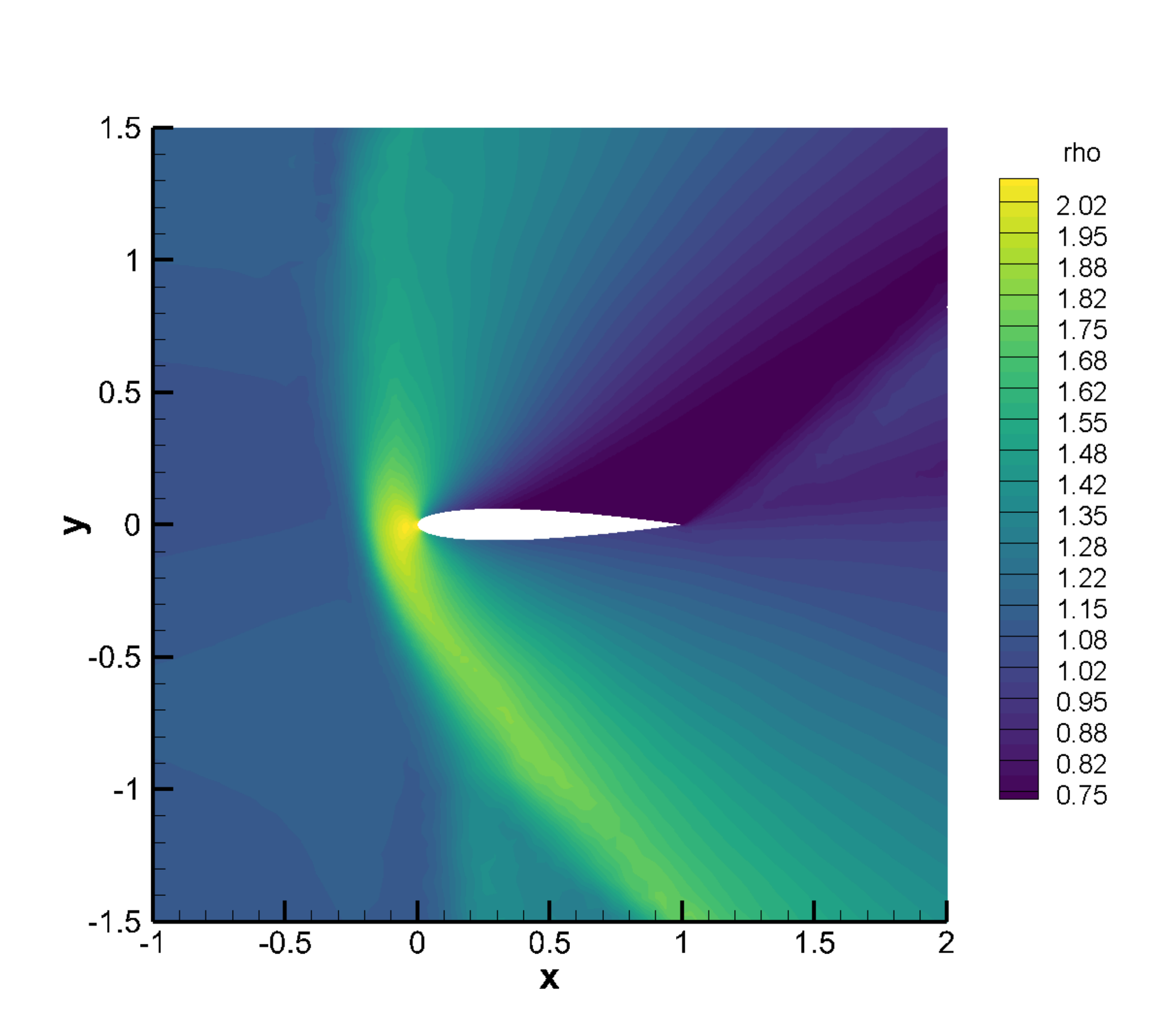}  \\         
			\includegraphics[width=0.33\textwidth]{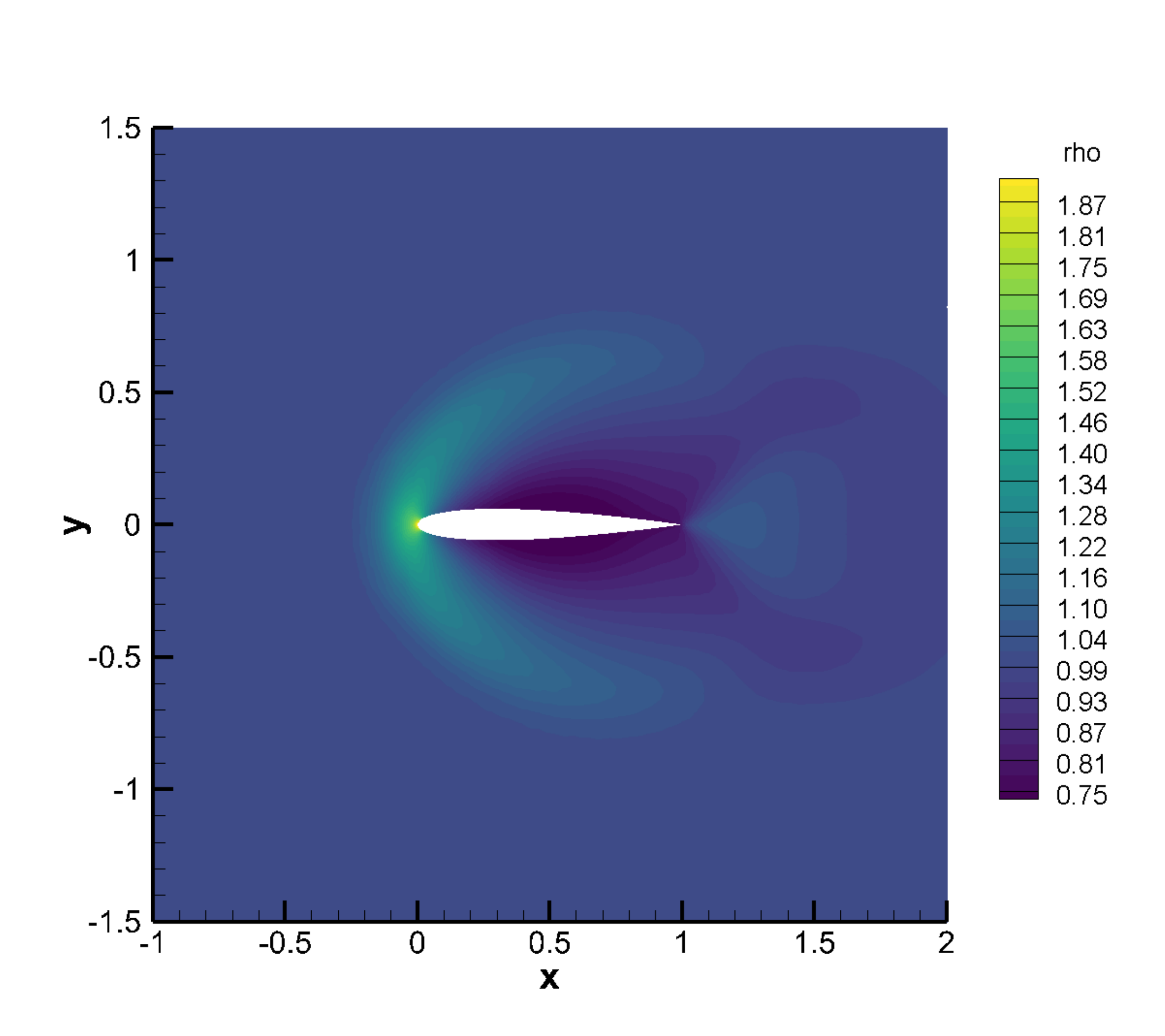}  &           
			\includegraphics[width=0.33\textwidth]{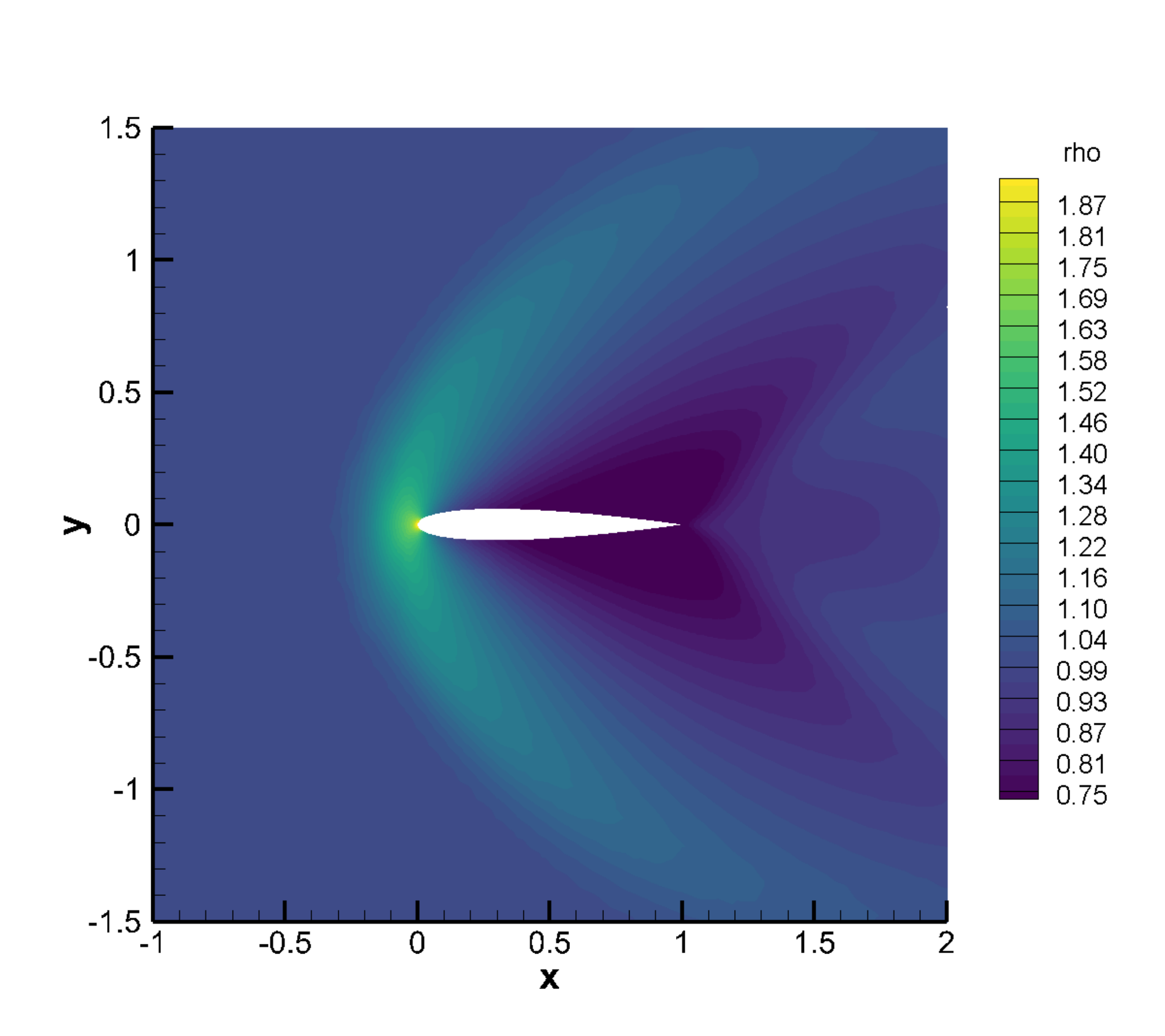}  & 
			\includegraphics[width=0.33\textwidth]{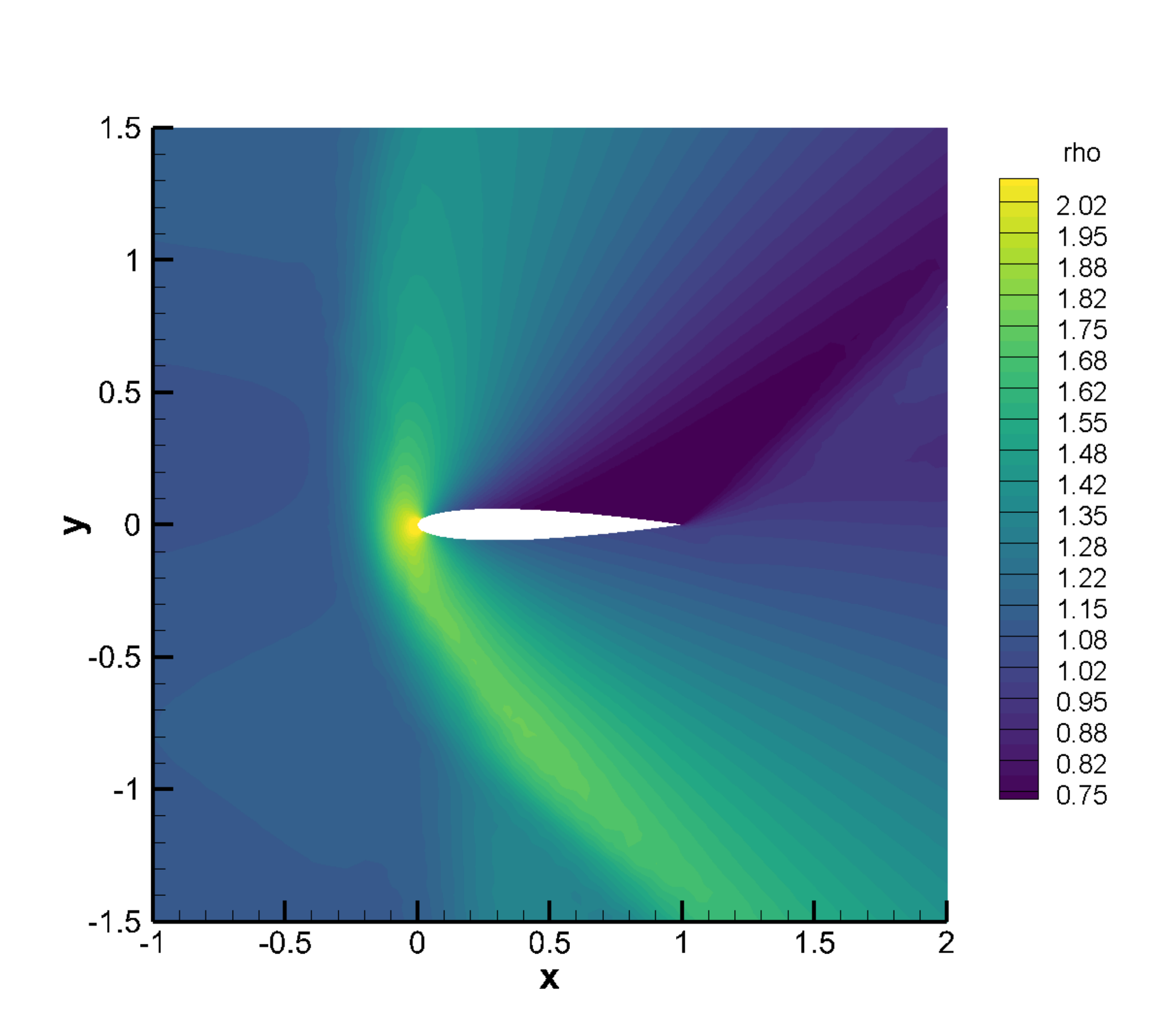}  \\
		\end{tabular}
		\caption{Flow around NACA 0012 airfoil with $\varepsilon=5\cdot 10^{-4}$ (top row) and $\varepsilon=5\cdot 10^{-3}$ (bottom row). Density distribution (40 contours in the interval $[0.65;2.05]$) at time $t=0.4$ (left), $t=1$ (middle) and $t=5$ (right).}
		\label{fig.NACA_Kn_rho}
	\end{center}
\end{figure}

\begin{figure}[!htbp]
	\begin{center}
		\begin{tabular}{cc} 
			\includegraphics[width=0.47\textwidth]{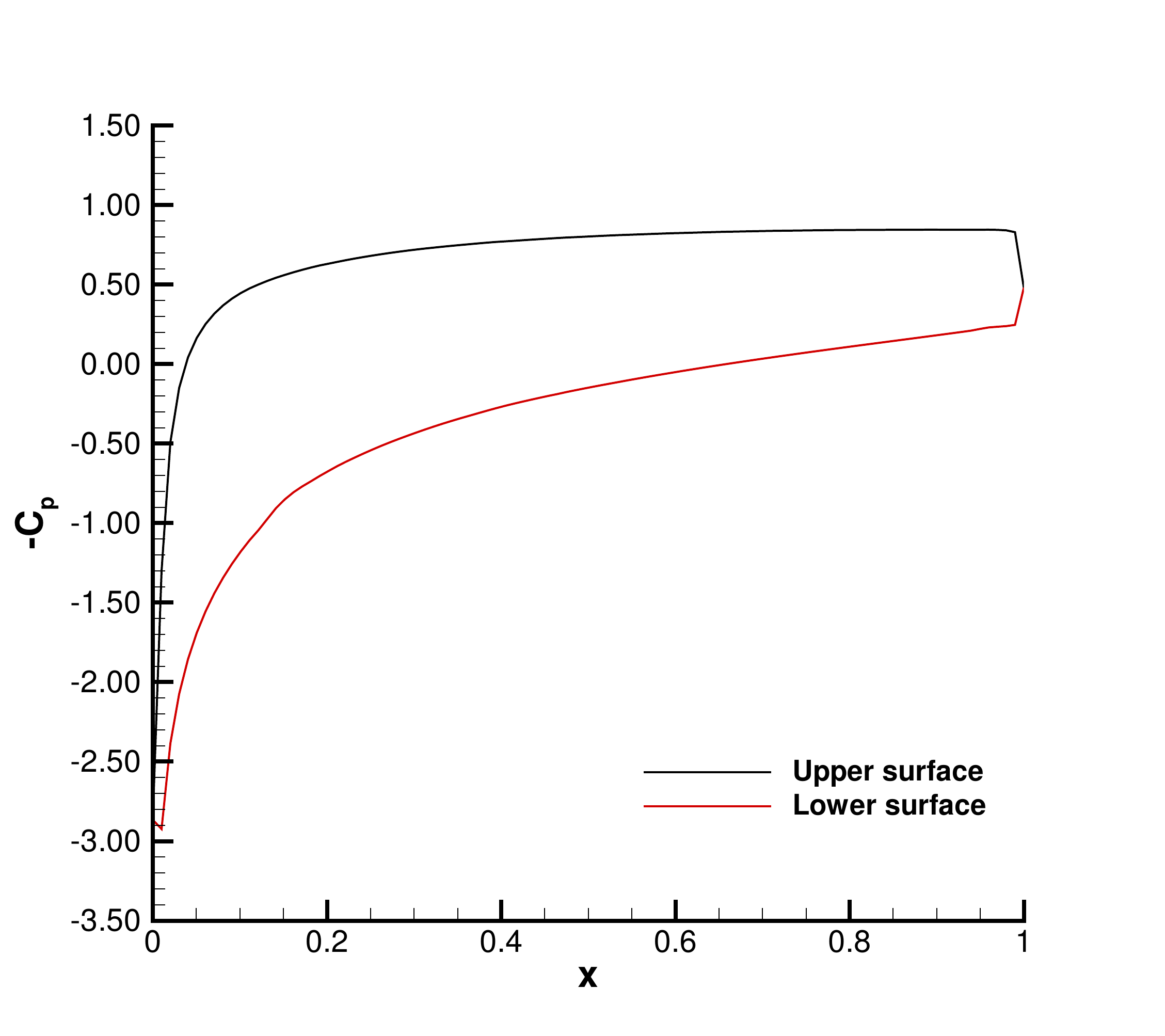}  &           
			\includegraphics[width=0.47\textwidth]{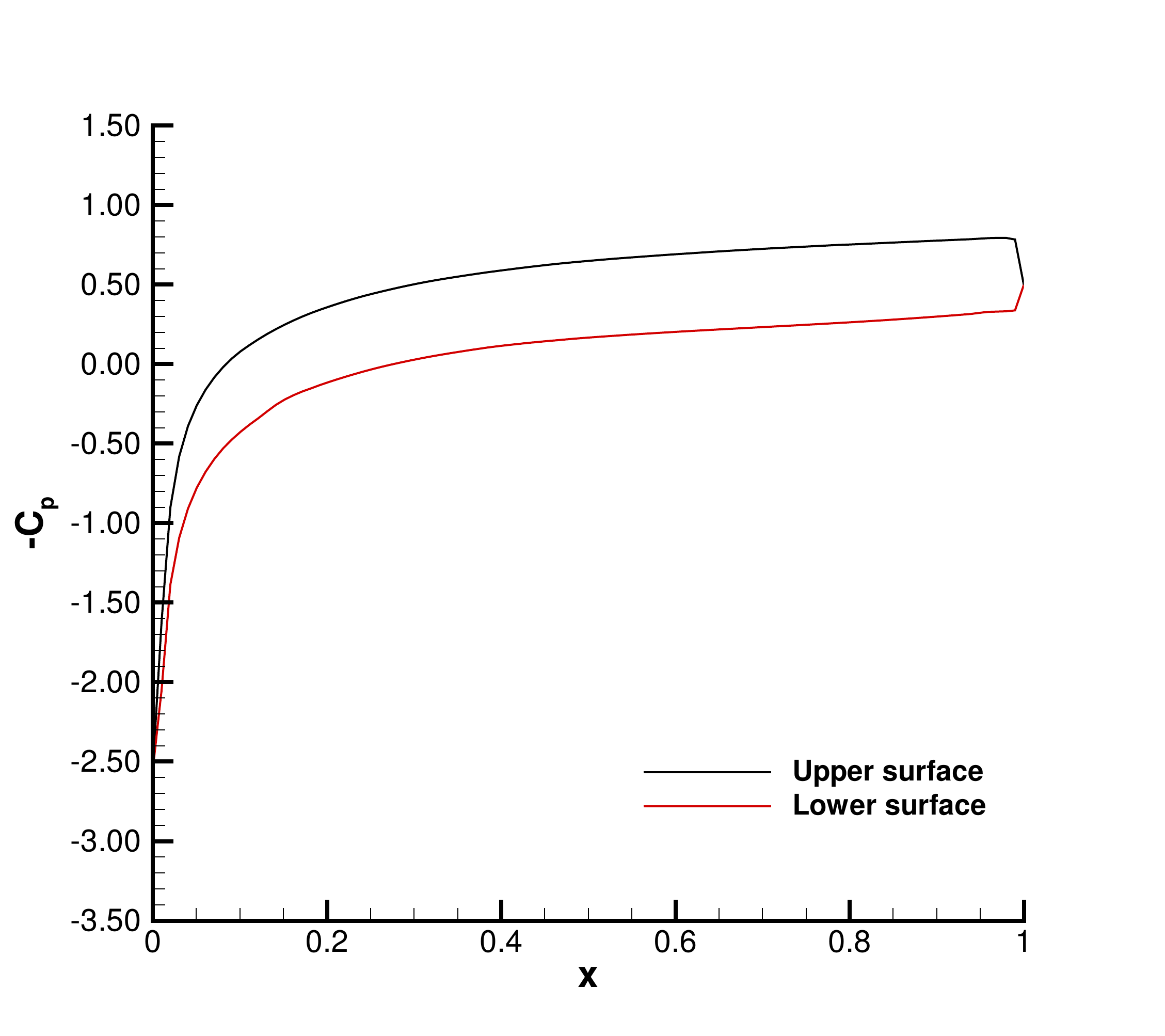}  \\ 
			\includegraphics[width=0.47\textwidth]{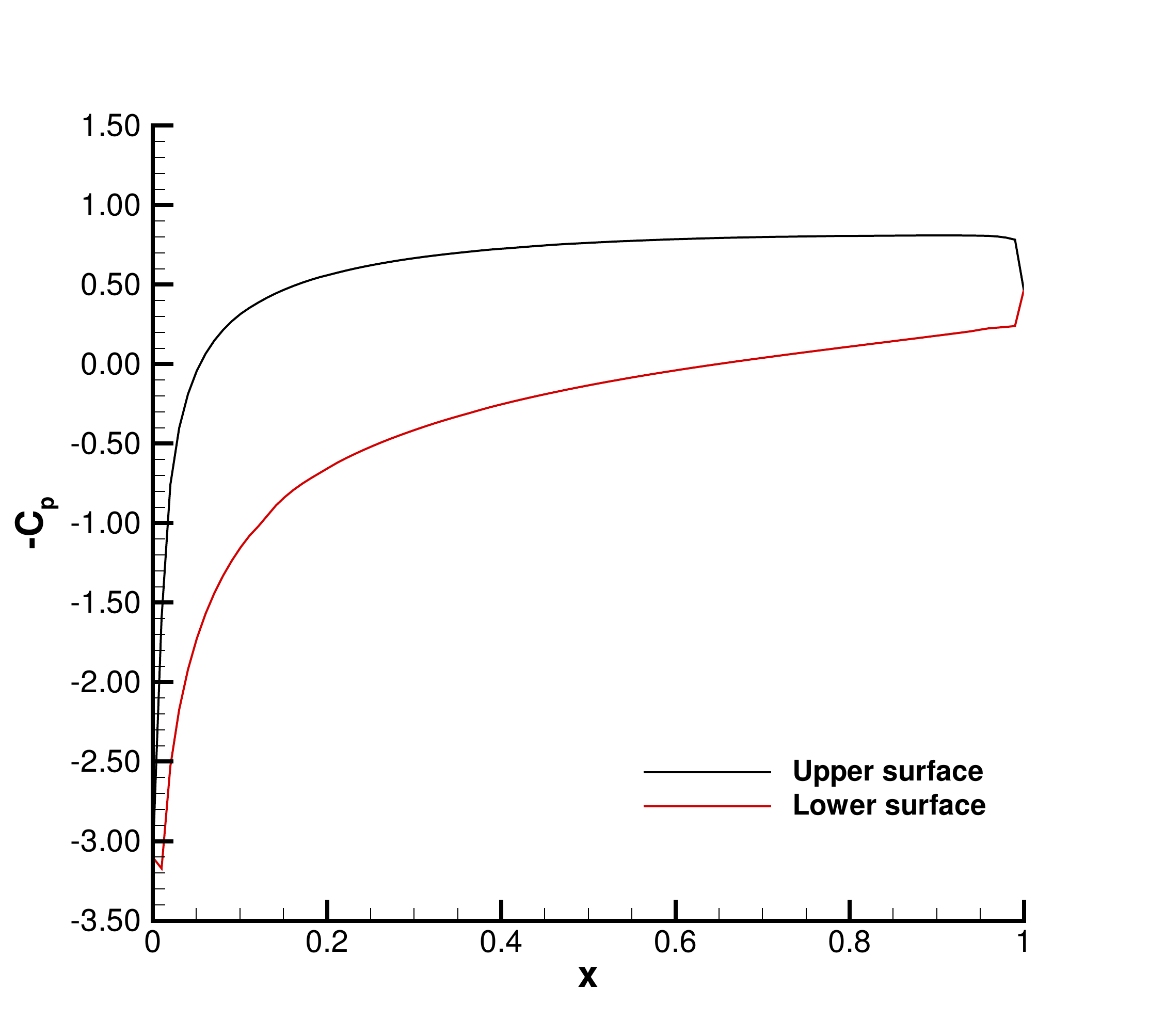}  &           
			\includegraphics[width=0.47\textwidth]{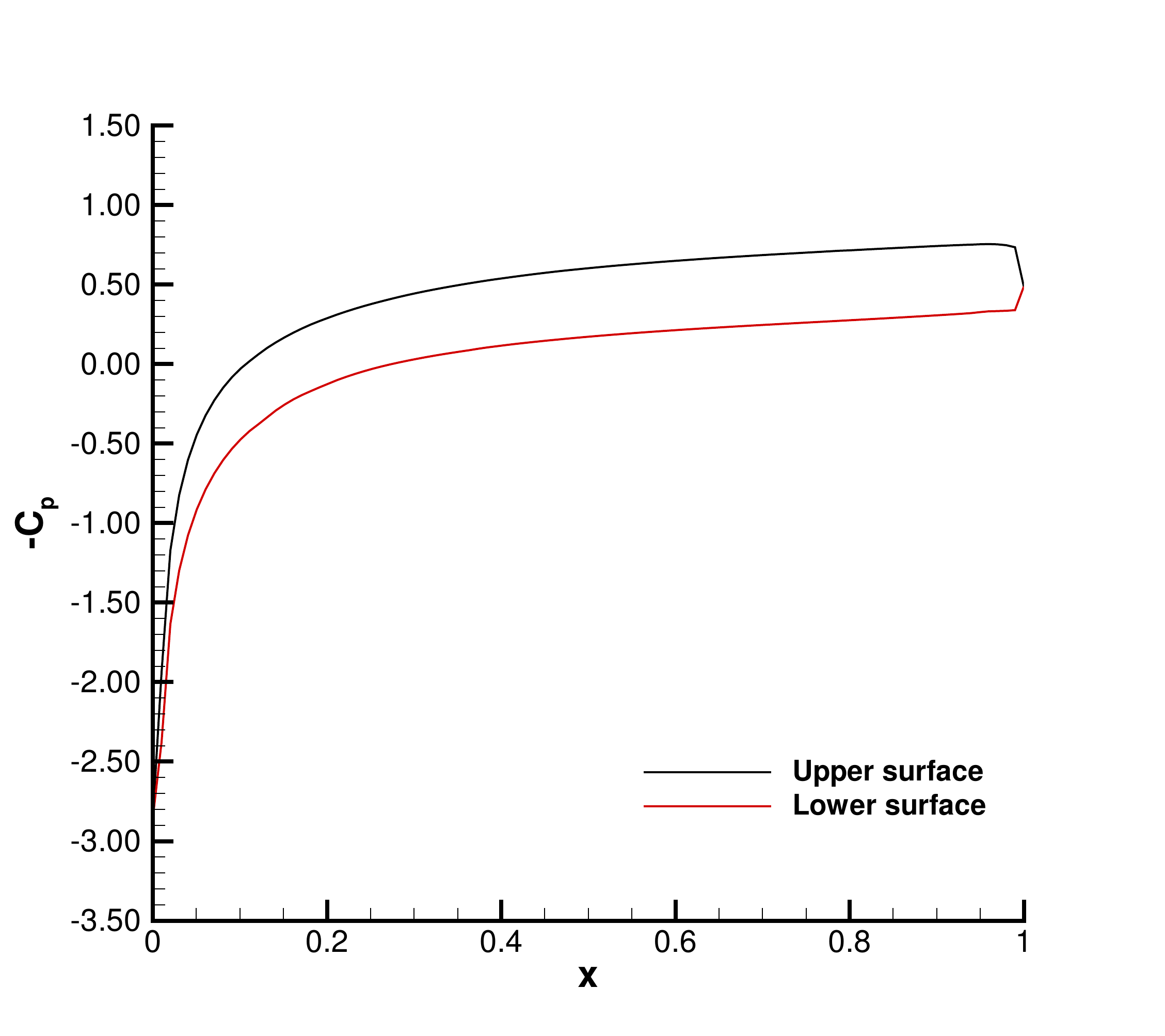}  \\ 
		\end{tabular}
		\caption{Flow around NACA 0012 airfoil with $\varepsilon=5\cdot 10^{-4}$ (top row) and $\varepsilon=5\cdot 10^{-3}$ (bottom row). Pressure coefficient distribution on the upper and lower surfaces of the airfoil at time $t=3$ (left) and $t=5$ (right).}
		\label{fig.NACA-Cp}
	\end{center}
\end{figure}

The same simulation is also computed for the case $\varepsilon=5 \cdot 10^{-3}$, thus considering a more rarefied gas. Figure \ref{fig.NACA_Kn_rho} plots the density distribution around the wing at three different output times for both Knudsen numbers, showing that a denser fluid is correctly approximated with $\varepsilon=5 \cdot 10^{-4}$. In the case $\varepsilon=5 \cdot 10^{-4}$, shocks are also stronger as evident from the higher compression of the gas especially at time $t=5$. However, shocks are visible at the same space-time locations, hence exhibiting the same flow structure due to the supersonic flow conditions. Finally, for this problem we also compute the pressure coefficient $C_p$ on the upper and the lower surface of the airfoil at time $t=3$ and $t=5$ for both Knudsen numbers. This is obtained through the following relation
\begin{equation}
	C_p=\frac{p-p_\infty}{\frac{1}{2} \rho_\infty u_\infty^2}.
\end{equation}
Figure \ref{fig.NACA-Cp} depicts the distribution of $C_p$ at time $t=3$ and $t=5$. The denser fluid generates a higher pressure coefficient at the head of the airfoil and the area bounded between the upper and lower surfaces is slightly larger in the case $\varepsilon=5\cdot 10^{-4}$.
These simulations demonstrate that the new class of DG-IMEX schemes can successfully be applied to more realistic settings involving wing profiles or other objects in atmosphere as for instance the simulation of a vehicle re-entry approaching the soil surface.

\section{Conclusions} \label{sec_conclu}

In this work, we have presented a novel Discontinuous Galerkin Implicit-Explicit spectral numerical method for solving the Boltzmann equation. The class of schemes we have developed enjoys high order accuracy in velocity, in space and in time. On this regard, we have presented up to third order in time and space methods while spectral accuracy has been guaranteed for the solution of the Boltzmann collision integral. Time step limitations typically encountered when solving multiscale kinetic models are overcome thanks to the use of implicit methods for the stiff collisional part. Although being formally implicit, these new DG methods are explicitly solvable, i.e. the inversion of highly nonlinear systems is completely avoided. Moreover, for the first time linear multistep methods have been used in applications involving the resolution of the Boltzmann equation showing their better properties in terms of accuracy and their reduced computational cost in comparison with Runge-Kutta methods. For the space discretization, the modal DG approach has been constructed on arbitrarily polygonal shaped meshes and special care has been taken to avoid the onset of spurious oscillations without compromising the accuracy. 

Next, the proposed approach has been numerically validated. Extensive convergence studies have been performed for all the discretizations employed in this work, considering different rarefied regimes of the gas. Classical fluid dynamics problems involving shock waves and other discontinuities have been reproduced, and comparisons against other discretization techniques as well as other simplified relaxation models have been analyzed. A last part has been dedicated to an application showing that the methods here proposed are ready to be employed in more realistic settings.

In the future, we aim at extending our results to the full six-dimensional case, to work on the mesh adaptation in velocity space to be able to handle very different physical settings, and to develop efficient integration techniques in space in the case of arbitrarily polygonal control volumes with the scope of reducing the computational cost.

\section*{Acknowledgments}
The Authors would like to thank the Italian Ministry of Instruction, University and Research (MIUR) to support this research with funds coming from PRIN Project 2017 (No. 2017KKJP4X entitled "Innovative numerical methods for evolutionary partial differential equations and applications"). 

\bibliographystyle{plain}
\bibliography{biblio}

\end{document}